 \documentclass[draft]{article}

\usepackage{amsmath,amsfonts,amsthm,amssymb,amscd,color}

\renewcommand{\Re}{\mathop{\rm Re}\nolimits} \renewcommand{\Im}{\mathop{\rm
Im}\nolimits}   
  \newcommand{\Ph}{{\mathcal P}}

\theoremstyle{plain} \newtheorem{theorem}{Theorem}[section]
\newtheorem{lemma}[theorem]{Lemma}
\newtheorem{proposition}[theorem]{Proposition}
 \theoremstyle{definition}
\newtheorem{definition}[theorem]{Definition} \theoremstyle{remark}
\newtheorem{remark}[theorem]{Remark} 
 
\newcommand{\R}{{\mathbb R}} \newcommand{\U}{{\mathcal U}}

\newcommand{\Z}{{\mathbb Z}}

\newcommand{\N}{{\mathbb N}}

\newcommand{\resto}{{\mathcal R}} \def\im{{\rm i}}

\newcommand{\C}{\mathbb{C}} 

\def\({\left(}
\def\){\right)}
\def\<{\left\langle}
\def\>{\right\rangle}

\numberwithin{equation}{section}

\begin{document}

  \title{On weak interaction between  a ground state and a non-trapping  potential}

\author {Scipio Cuccagna , Masaya Maeda}

 \maketitle
\begin{abstract}  We  show that ground states of the NLS  moving at nonzero speed are asymptotically stable  if    they either stay far from the potential, or the potential is small, or the
ground state has large speed.
We search an effective Hamiltonian using the Birkhoff normal forms argument in \cite{Cu0},  treating the potential as a perturbation. The so called Fermi Golden Rule, which is used to describe
the decay to 0 of the internal discrete modes of  the ground state, is similar
to that in  \cite{Cu2}. The continuous modes dispersion requires the theory  in
\cite{RSS1} on charge transfer models.

\end{abstract}

\section{Introduction}

We consider the   nonlinear Schr\"odinger equation with a potential (NLSP)

\begin{equation}\label{NLSP}
 \im w_{t }=-\Delta  w +\gamma V(x)w  +\beta  (|w|^2) w  \ , \quad
 (t,x)\in\mathbb{ R}\times
 \mathbb{ R}^3\ .
\end{equation} We have the following   first set of hypotheses.

\begin{itemize} \item[(H1)]      $V\in {\mathcal S}(\R
^3, \R)$ is a   Schwartz function,  $\gamma >0 $,  $  -\Delta +\gamma V  $ has no eigenvalues  and 0 is not a resonance (that is, if $(-\Delta +\gamma V)u=0$  with $u\in  C^\infty$ and     $|u(x)|\le C|x|^{-1} $ for a fixed $C$, then $u=0$).

 \item[(H2)]
$\beta  (0)=0$, $\beta\in C^\infty(\R,\R)$. \item[(H3)] There exists a
$p\in(1,5)$ such that for every $k\ge 0$ there is a fixed $C_k$ with $$\left|
\frac{d^k}{dv^k}\beta(v^2)\right|\le C_k |v|^{p-k-1} \quad\text{if $|v|\ge
1$}.$$ \end{itemize}

We think of $\gamma V(x)w$  as a perturbation  of

\begin{equation}\label{NLS}
 \im w_{t }=-\Delta w   +\beta  (|w|^2) w  .
\end{equation}

A ground state is a solution  $e^{\im\omega t}\phi(x)$ of \eqref{NLS} in $H^2 (\R ^3)$       with $\omega\in \R$ and $\phi(x)>0$.
It   is {\it orbitally stable} if $\forall$ $\epsilon>0$  $\exists$ $\delta>0$ s.t.  if $\|\phi-u_0\|_{H^1}<\delta$
  the solution $u$ of \eqref{NLS} with $u(0)=u_0$ exists globally in time and satisfies
\begin{equation*}
\sup_{t>0}\inf_{s\in\R, y\in\R^3}\|e^{\im s }\phi(\cdot-y)-u(t)\|_{H^1}<\epsilon.
\end{equation*}
By the following hypotheses our  ground states are orbitally stable,  \cite{W2,GSS1,GSS2}:

\begin{itemize} \item[(H4)] on an open interval $\mathcal{O}\subseteq (0, \infty )$ for any $n\in \N$  there is  a $C^\infty $-family $\mathcal{O}\ni\omega\mapsto \phi_{\omega}\in \Sigma _{n}$, see \eqref{eq:sigma} below,  s.t.  $\phi_{\omega}$ is a positive radial solution of \begin{equation}
  \label{eq:B}
  -\Delta u + \omega u+\beta(|u|^2)u=0\quad\text{for $x\in \R^3$.}
\end{equation}

\item [(H5)] $
\frac d {d\omega } \| \phi _ {\omega }\|^2_{L^2}>0 \quad\text{for
$\omega\in\mathcal{O}$.}$

\item [(H6)] Let $L_+:=-\Delta   +\omega +\beta (\phi _\omega ^2 )+2\beta
    '(\phi _\omega ^2) \phi_\omega^2$ with domain   $H^2
 (\R^3)$. We assume   $L_+$ has   one negative eigenvalue
 and
  $\ker L_+  =\text{Span}\{ \partial_{x_j}\phi_{\omega}:  j=1,2,3\}$.
\end{itemize}

By the symmetries of \eqref{NLS}, the following  is a family of solutions  of \eqref{NLS}
\begin{equation}\label{eq:gsset}
\mathcal{G}_{\omega}:=\{e^{\frac \im 2 \textbf{v} \cdot x-\frac \im  4 |\textbf{v} |^2 t +\im  t\omega +\im  \gamma _0 } \phi _\omega (x-\textbf{v} t- D_0) \ |\ \gamma_0 \in \R,\ \textbf{v}, D_0 \in \R^3 \}.
\end{equation}
    We study the effect of the potential  $\gamma V$
on
   a solution $u(t)$ initially close to  $\mathcal{G}_{\omega_0}$.     Setting  $u(t,x):= e^{-\frac \im 2 \textbf{v} \cdot x-\frac \im 4 t |\textbf{v} |^2  }    w(t,x+\textbf{v} t)$
we
can equivalently rewrite   \eqref{NLSP} as
\begin{equation}\label{eq:NLSP}
\begin{aligned}
 &\im \dot u =-\Delta   u   +\gamma V(x+\textbf{v} t) u  +\beta  (|u|^2) u\ ,\quad
 u(0,x)=u_0(x).
\end{aligned}
\end{equation}

\subsection{Statement of the main result} \label{subsec:statement}

We assume  the following hypotheses  on the linearization  $\mathcal{H}_\omega$
defined in \eqref{eq:linearization1}. \begin{itemize}

\item [(H7)]  $\exists$ $\textbf{n}$  and $0<\textbf{e} _1(\omega )\le \textbf{e} _2(\omega )
\le ...\le \textbf{e} _{\textbf{n}}(\omega )$, s.t.
$\sigma_p(\mathcal{H}_\omega)$ consists    of  $\pm   \textbf{e} _j(\omega )$ and  $0$ for $j=1,\cdots, \textbf{n}$.
We   assume $0<N_j\textbf{e} _j(\omega )<
\omega < (N_j+1)\textbf{e} _j(\omega )$ with $N_j \in \N$. We set $N=N_1$. Here
 each  eigenvalue  is repeated a number of times equal to its  multiplicity.   Multiplicities  and $\textbf{n}$ are  constant in $\omega$.

 \item [(H8)]  There is no multi
index $\mu \in \mathbb{Z}^{\textbf{n}}$ with $|\mu|:=|\mu_1|+...+|\mu_k|\leq 2N_1+3$
such that $\mu \cdot \textbf{e}(\omega) =\omega $, where
 $\textbf{e}(\omega)=(\textbf{e}_1(\omega),\cdots,\textbf{e}_{\textbf{n}}(\omega))$.

\item[(H9)] For $\textbf{e} _{j_1}(\omega)<...<\textbf{e} _{j_k}(\omega)$   and $\mu\in \Z^k$ s.t.
  $|\mu| \leq 2N_1+3$, then we have
\[ \mu _1\textbf{e} _{j_1}(\omega)+\dots +\mu _k\textbf{e} _{j_k}(\omega)=0 \iff \mu=0\ . \]
\item[(H10)]    The points $\pm \omega$ are not resonances of $\mathcal{H}_\omega$,
    see p.5 \cite{CPV}.

\item [(H11)] The Fermi golden rule  Hypothesis (H11)   in Sect.
    \ref{subsec:FGR}, see \eqref{eq:FGR}, holds.

\end{itemize}

Solutions of  \eqref{eq:NLSP}   initially  close to a $\mathcal{G}_{\omega}$, for some time can be written as
 \begin{equation} \label{eq:Ansatz1}\begin{aligned}&  u(t,x) =  e^{\im \left ( \frac
12 v(t)\cdot  x  + \vartheta (t)\right )}
 \phi _{\omega (t)} (x -D(t))+ r(t,x  )   .
\end{aligned}
\end{equation}
We will show   this   persists  for all $t\ge 0$ and that the      $(\omega (t) ,v(t) )$   converge.

\begin{theorem}\label{theorem-1.1}
  Let $\omega_1\in\mathcal{O}$.  Assume $(\mathrm{H1})$--$(\mathrm{H11})$   and   $u_0\in {\mathcal S}(\R ^3)$. Fix $M_0>1$ and  $ \mathbf{v} ,y_0   \in \R ^3   $ with   $|\mathbf{v}| > M_0 ^{-1}$
    and $|\gamma | <M_0$, where $\gamma$ is the constant in $(\mathrm{H1})$.  Fix a $\delta _0>0$ small.
    Set
	\begin{equation*} \epsilon:=\inf_{\theta  \in \R }\|u_0-e^{\im
\theta }  \phi_ {\omega _1 }(\cdot -y_0)
\|_{H^1}   + \gamma  \sup _{\text{dist} _{S^2}( \widehat{e}, \frac{\mathbf{v}}{|\mathbf{v}|} ) \le \delta _0} \int _0^\infty  (1+|  |\mathbf{v}|\widehat{e} t+y_0   |^{ 2} )^{-1} dt .
\end{equation*}
	   $\exists$  $\epsilon_0 =  \epsilon_0 (M_0,\omega_1,\delta _0 )>0$ $\&$
$C=  C (M_0,\omega_1,\delta _0 )>0$ s.t. if
$u $ solves
\eqref{eq:NLSP} and
\begin{equation}
	\label{eq:sizeindata} \epsilon   <\epsilon_0,
\end{equation}
  there exist $\omega _+\in\mathcal{O}$, $v_+ \in
\R^3$, $\theta\in C^1(\R_+;\R)$, $ y\in C^1(\R_ +;\R ^3)$ and $h _+ \in H^1$ with
$\| h_+ \| _{H^1}+|\omega _+ -\omega_1  |+ |v_+  |\le C \epsilon $
such that
\begin{equation}\label{eq:scattering} \lim_{t\to +\infty}\|u(t,x)-e^{\im
\theta(t)+  \frac  \im 2{v_+ \cdot x}  } \phi_{\omega _+} (x -y(t) )-
e^{\im t\Delta }h _+ (x)   \|_{H^1_{x} }=0. \end{equation} In the notation of
\eqref{eq:Ansatz1},
 we have  $    r(t,x  )=A(t,x)+\widetilde{r}(t,x)$  such that
$A(t, \cdot ) \in \mathcal{S}(\R^3, \C)$, $|A(t,x)|\le C (t)$ with  $\lim
_{t  \nearrow \infty}C (t)=0$ and such that  $\|
\widetilde{r} \| _{L^p_t( \mathbb{R}_+,W^{1,q}_x)}\le
 C\epsilon $
for any admissible pair $(p,q)$, by which we mean that \begin{equation}\label{admissiblepair}
2/p+3/q= 3/2\,
 , \quad 6\ge q\ge 2\, , \quad
p\ge 2. \end{equation}
\end{theorem}

\begin{remark}
By    $(\mathrm{H2})$, $\beta(|w|^2)w  =O(w^3)$  at   $0$.
This  excludes bound states with small $H^1$ norm and multi-soliton solutions near the ground state.
\end{remark}

 Theorem \ref{theorem-1.1} for   $\gamma =0$    is in \cite{Cu3}.   Here we extend  \cite{Cu3} by adding a
non-trapping potential. We consider only the regime   when the interaction of the soliton with   the potential
is either too quick or is weak because they are far or because the potential is small.      See \cite{JFGS0,JFGS,GHM,HMZ1, HMZ2,HZ2,HZ1,AS,ALS,perelman2}  for various studies on behavior  of ground states and a potential. See also   \cite{M1,M2}.
A related topic is the analysis of multisolitons, see   \cite{RSS2,perelman3, perelman4,MMT,MM,ASFS,CLC}.

  \noindent Theorem \ref{theorem-1.1}   is expected  since for weak interaction of soliton and potential
one  expects an analogue of    the classical particle scattering by a potential.

 \noindent
For   slower velocities  and  classically trapping   potentials   the   time evolution can involve  long time  oscillatory  motion
and complicate patterns, see \cite{JFGS0,JFGS,ASFS}.  Our    result should    be contrasted with the   case considered in \cite{perelman2},   see also  \cite{HMZ1, HMZ2,perelman1},
where the 1D soliton   hits  a relatively large  potential at a not sufficiently large velocity   to survive    and  collapses  becoming radiation,  see
also  \cite{perelman1}.

Our   interaction    is similar to the      solitons interaction in \cite{RSS2,perelman3, perelman4} but is less restrictive since we allow {\it internal modes},
that is the existence of the eigenvalues $\textbf{e}_j(\omega )$  in (H7). This makes the local dynamics near a soliton
  harder to analyze.

The first results on the decay of internal modes go back to  \cite{BP2,SW3}, which treated only case $N_1=1$.  A  quite general analysis of this decay was obtained  first in \cite{bambusicuccagna},
then in \cite{Cu2,Cu3}.
For a more general and simpler  set up see   \cite{bambusi,Cu0}.

Here we are in a situation similar to that in   \cite{Cu3,Cu0}. Although the soliton in not an exact solution
of   \eqref{eq:NLSP}, we still  represent  solutions  $u(t)$  as a sum of a soliton and a remainder. By the weakness of the interaction with the potential,  this representation is preserved for all times. Furthermore,   up to   a solution of the  constant coefficient linear Schr\"odinger  equation  and up to a phase factor and  to translation, our solution
 becomes a ground state of  \eqref{NLS}.   In   \cite{Cu3,Cu0} this is shown using an appropriate effective
Hamiltonian in a neighborhood of the ground state with the same invariants of motion of $u(t)$.   Since the linear momenta are not
invariants of  \eqref{eq:NLSP},
here we   modify this  argument    by considering an appropriate time dependent effective Hamiltonian.

Once we get     right  time dependent coordinate system  and  effective Hamiltonian,   we
  close estimates    as in   \cite{Cu3,Cu2}.  The proof resembles
the standard analysis of scattering for semilinear equations, see Ch.  6 \cite{strauss}, supplemented  with   the
  \textit{Fermi Golden Rule}  (FGR), which involves showing
  square power structure and positive semi--definiteness for
  some       coefficients of the internal modes equations,
  see \eqref{eq:FGR8}--\eqref{eq:FGR81}.  In  (H11) we assume  they are strictly positive  (proved in Prop.2.2 \cite{bambusicuccagna} for      an easier problem).
 Then, by  nonlinear interaction with continuous modes   the energy of the
internal modes (which   left on their own
would    look like harmonic oscillators)
leaks into the continuous modes where is scattered by
linear dispersion. This discussion is in Sect. \ref{subsec:FGR} and is the same of   \cite{Cu2}.

  The scattering of  continuous modes is in Sect. \ref{subsec:cont mode}--\ref{subsec:coupling}, \ref{sec:dispersion}--\ref{sec:completion}. It
   requires    a perturbation argument which exploits linear dispersion.
	The   linearization is similar
to the   {\it charge transfer models}   in \cite{RSS1}, applied  to   interaction of solitons in   \cite{RSS2}.
While  we follow  \cite{RSS1}  and our nonlinear problem is similar to   \cite{RSS2},   the      slow decay   due to the    internal modes   forces us to    add to the  argument   also
Beceanu's   integrating factor \cite{beceanu},   needed to eliminate the terms       $v  \cdot \nabla _x u $ and $\varphi    \sigma _3u$  in \eqref{eq:strich1},    in \cite{RSS2}     eliminated using   fast decay due to absence of internal modes.

    We work in 3$d$. The arguments extend to higher $d$. In 1$d$ and $2d$
     Beceanu's argument  is unavailable due to the fact that      $\Delta $ has a resonance at 0.

We focus only on   case $ \sigma ( -\Delta +\gamma V)=[0,\infty ).  $
 We leave  to a future paper the case of
 $-\Delta +\gamma V$  with exactly one eigenvalue, where asymptotic behavior of
 small energy solutions of \eqref{NLSP} is known, \cite{GNT}.  Notice that \cite{CM1}
  has extended 
  \cite{GNT} generalizing   \cite{SW4,TY3,NPT}.

\section{Set up} \label{section:set up}  We start with some notation.
 We set $\langle x \rangle = (1+|x|^2)^{\frac{1}{2}}$ and
\begin{equation} \label{eq:hermitian}\begin{aligned}& \langle f,g\rangle =\Re  \int
_{\mathbb{R}^3}f(x)  \overline{g}(x)  dx \text{ for $f,g:\mathbb{R}^3\to \mathbb{C}$ }
 .  \end{aligned}\end{equation}
We  identify $\C =\R^2$ and set $J=\begin{pmatrix}  0 & 1  \\ -1 & 0
 \end{pmatrix}$. Multiplication by $\im $ in $\C$ is $J^{-1}=-J$.
For   $n \ge 1$  and   $K=\R , \C$ then   $\Sigma _n=\Sigma _n(\R ^3, K^2   )$
is the   Banach space   with
 \begin{equation}\label{eq:sigma}
\begin{aligned} &
      \| u \| _{\Sigma _n} ^2:=\sum _{|\alpha |\le n}  (\| x^\alpha  u \| _{L^2(\R ^3   )} ^2  + \| \partial _x^\alpha  u \| _{L^2(\R ^3   )} ^2 )  <\infty .      \end{aligned}
\end{equation}
	 We set $\Sigma _0= L^2(\R ^3, K^2   )$.  We   define $\Sigma _{r}$   by  $
      \| u \| _{\Sigma _r}  :=  \|  ( 1-\Delta +|x|^2)   ^{\frac{r}{2}} u \| _{L^2}   <\infty    $ for $r\in \R$.
 For $r\in \N$ the two definitions are equivalent, see    \cite{Cu3}. Notice that for $\mathcal{S}=\mathcal{S}(\R ^3, K^2   )$ the set of Schwartz functions and for
 $\mathcal{S}'=\mathcal{S}'(\R ^3, K^2   )$ the set of tempered distributions we have
 \begin{equation}\label{eq:schwartz}
\begin{aligned} &
      \mathcal{S}= \cap _{n\in \Z} \Sigma _n   \, , \quad    \mathcal{S}'= \cup _{n\in \Z} \Sigma _n.\end{aligned}
\end{equation}

\noindent  We will      by $H^{k,s} (\R ^3, K^2   )  $ the Banach space  defined by
\begin{equation}
\begin{aligned} &
      \| u \| _{H^{k,s}}  :=  \|   (1+|x|^2) ^{\frac{s}{2}}    ( 1-\Delta  )   ^{\frac{k}{2}} u \| _{L^2}   <\infty    .   \end{aligned}\nonumber
\end{equation}
 We set
 $L^{2,s} =H^{0,s}  $, $L^2=L^{2,0}    $,  $H^k=H^{k,0}    $.
To emphasize  spatial variables, we might
denote spaces by $W^{k,p}_x$, $L^{ p}_x$, $H^k_x$, $H^{ k,s}_x$ and $L^{2,s}_x$.
For $I\subset \R $   and $Y_x$ any of them, we   consider
spaces $L^p_t( I, Y_x)$ with   norm $ \| f\| _{L^p_t( I, Y_x)}:= \| \| f\|
_{Y_x} \| _{L^p_t( I )}.$

 Given two topological vector spaces $X$ and $Y$ we denote by $B (X,Y)$ the set of
 continuous   linear operators from $X$ to  $Y$.

\subsection{The linearization} \label{subsec:linearization}

For $B\in C^\infty(\R,\R)$ s.t. $B(0)=0$ and $ B'(t)=\beta (t) $,
we consider the energy
\begin{equation} \label{eq:energyfunctional}\begin{aligned}&
 \textbf{E}(u)=\textbf{E}_0(u)+ 2 ^{-1} {\gamma }     \langle
    V (\cdot + \textbf{v} t) u , u\rangle    \\&
\textbf{E}_0(u):=2 ^{-1}\|
  \nabla u \| ^2_{L^2} + \textbf{E}_P(u) \ , \   \textbf{E}_P(u):=2 ^{-1}
 \int _{\R ^3}B(| u|^2  ) dx . \end{aligned}\end{equation}

\noindent For    $u\in H^1( \R ^3, \C ) $  we have  charge and momenta: 	
\begin{equation}\label{eq:charge}
\begin{aligned}
 & Q(u)= \Pi _4(u)= 2 ^{-1} \|  u  \|  ^2_{L^2}= 2 ^{-1}\langle   \Diamond _4 u , u \rangle
\, ,  \quad  \Diamond _4:=1 ;\\
&     {\Pi }_a(u)=  2 ^{-1} \Im \langle  u_{x_a},  {u}  \rangle
= 2 ^{-1} \langle   \Diamond _a u ,  {u} \rangle  \, , \quad \Diamond _a:=J\partial _{x_a}\text{ for $a = 1,2, 3$.}
\end{aligned}
\end{equation}
 We set $\Pi  (u)=( \Pi _1(u),...,\Pi _4(u))$. We have
$\textbf{E}\in C^2 ( H^1( \R ^3, \C  ), \R  )$ and $\Pi _j\in C^\infty (  H^1(
\R ^3, \C ) , \R  )$. By straight computation we have
   the following   formulas:
 \begin{align} &Q (e^{-\frac 1 2 Jv\cdot x   }u )= Q (u
) \, ; \nonumber  \quad
   \Pi _a (e^{-\frac 1 2 J   v\cdot x   }u )= \Pi _a(u) +2^{-1} v_a Q (u
   )  \text{ for $a = 1,2, 3$} \, ;
\\& \textbf{E}_0(   e^{-\frac 1 2 J v\cdot x   }u )= \textbf{E}_0(   u) +
v\cdot  \Pi   (u)
 + 4 ^{-1}{v^2 } Q (u )    \, ,  \,   v\cdot  \Pi   (u)=\sum _{a\le 3} v_a \Pi _a
 (u).  \label{eq:charge1}\end{align}
By \eqref{eq:charge1} and by (H5), the map $(\omega , v)\to  p=\Pi (e^{ -\frac 1 2 Jv\cdot x   }\phi _\omega     ) $ is a diffeomorphism  into an open
subset   $\mathcal{P}$ of $ \R ^4$.  For $p=p(\omega, v)\in \mathcal{P}$ set $\Phi
_p=e^{-\frac 1 2 Jv\cdot x   } \phi _\omega    $.

For $F\in C^1( \mathbf{U},\C)$ with  $\mathbf{U}$ an open subset  of $H^1$,  the gradient $\nabla F (U)$
is defined by
$  \langle \nabla F (U), X\rangle = dF(U) (X)$, with $dF(U)$ the  Frech\'et derivative at $U$. If $F\in C^2(\mathbf{U},\C ) $ it remains defined the linear operator  $\nabla ^2F(U):H^1\to H^{-1}$

    The $\Phi _p   $ are
constrained critical points of $\textbf{E}_0$ with associated Lagrange
multipliers $\lambda  (p) \in \R^4$ so that  $
 \nabla \textbf{E}_0(\Phi _p  )= \lambda   (p) \cdot \Diamond \Phi _p$, where
 we have
 \begin{equation}
	 \label{eq:LagrMult} \lambda _4(p) =-\omega  (p) -4 ^{-1}{v^2 (p)}    \, ,    \quad  \lambda _a(p):=v_a (p)
	\, \text{ for $a=1,2, 3$.}
 \end{equation}
We set also
\begin{equation}
	 \label{eq:dp} d(p):=\textbf{E}_0(\Phi _{p  }) -   \lambda  (p) \cdot      \Pi   (\Phi _{p  }).
 \end{equation}
For  any fixed   $\tau _0\in \R^4 $ 		a function    $u(t):=  e^{J( t   \lambda
(p) +\tau _0)\cdot \Diamond}\Phi _p  $ 		 is a  {solitary wave} solution of
$\im u_t=-\Delta u +\beta (|u|^2)u$. 		We now introduce
\begin{equation}\label{eq:linearizationL}
 \mathcal {L}_p:=   J(\nabla ^2
 \textbf{E}_0(\Phi _p  )- \lambda  (p) \cdot \Diamond )  .
 \end{equation}
By an abuse of notation, we set
\begin{equation}\label{eq:linearization0}
   \mathcal {L}_\omega := \mathcal {L}_p  \text{ when $v(p)=0$ and $\omega (p)=\omega$.}
 \end{equation}
We have the  following identity,   elementary to check, see \cite{Cu0} Sect.7,
\begin{equation}\label{eq:conjNLS}
 \begin{aligned}
&   {\mathcal L}_p =e^{-\frac 12 J v (p)\cdot x}
  {\mathcal L}_{\omega (p)} e^{\frac 12 J v (p)\cdot x} .
\end{aligned}
\end{equation}
   \begin{remark}\label{rem:conjNLS}
By \eqref{eq:conjNLS}   the spectrum of ${\mathcal L}_p$   depends only on $\omega(p)$.
\end{remark}

 (H5) implies that $\text{rank}  \left [ \frac{\partial
  \lambda _i}{\partial p _j}   \right ]      _{ \substack{ i\downarrow  \  , \
  j \rightarrow}}= 4$. This  and (H6) imply
\begin{equation} \label{eq:kernel1}\begin{aligned}
&\ker {\mathcal L}_p  =\text{Span}\{ J \Diamond _j    \Phi _p:j=1,..., 4   \} \text{ and}\\
&N_g ( {\mathcal L}_p ) = \text{Span}\{ J \Diamond _j    \Phi _p, \partial _{ \lambda _j}     \Phi _p :j=1,..., 4 \} ,
\end{aligned}
\end{equation}
where   $N_g ( L ) :=\cup _{j=1}^\infty \ker (L^j)$.
Recall that we have a well known decomposition \begin{align} 	
\label{eq:begspectdec2}& L^2= N_g(\mathcal{L}_p)\oplus N_g^\perp
(\mathcal{L}_p^{\ast}) \  ,
   \\& N_g (\mathcal{L}_p^{\ast})  =\text{Span}\{   \Diamond _j    \Phi _p,
   J^{-1}\partial _{ \lambda _j}     \Phi _p :j=1,..., 4  \}   .
   \label{eq:begspectdec3}
\end{align}

\begin{lemma} [Modulation Lemma]
  \label{lem:modulation}
     Fix $\underline{n} \in \Z$ and   $\Psi _1=e^{J \tau _1 \cdot \Diamond}\Phi _{p_1}$. Then there exists a neighborhood $\U _{\underline{n}  } $   of $\Psi _1$    in $ \Sigma _{-\underline{n}}(\R ^3, \R ^{2 }) $
    and   functions $p  \in C^\infty (\U _{\underline{n}}   , \mathcal{P})$
	and $\tau \in C^\infty (\U _{\underline{n}}   , \R ^{4   })$ s.t.  $p(\Psi _1)=p_1 $   and $\tau
(\Psi _1) =\tau _1$
  and s.t.  $\forall u\in \U _{\underline{n}}   $
\begin{equation}\label{eq:ansatz}
\begin{aligned}
  &
 u =   e^{J \tau \cdot \Diamond} (  \Phi _{p } +R)
 \text{  and $R\in N^{\perp}_g (\mathcal{H}_p ^*)$.}
\end{aligned}
\end{equation}

\end{lemma}
\proof  First of all we point out that  $R\in N^{\perp}_g (\mathcal{L}_p ^*)$ makes sense for any tempered distribution $R$ and    if  $u  $  is in $H^k$ or   $ \Sigma _{k} $ with $k\ge - {\underline{n}} $  the same is true for $R$. This because the functions  in \eqref{eq:kernel1}
  are Schwartz functions.

By Sect. 7 \cite{Cu0} the hypotheses (A1)--(A6) of \cite{Cu0} hold here.
Hypotheses (B1)--(B2) of \cite{Cu0}  hold here by the discussion between \eqref{eq:charge1}
and \eqref{eq:LagrMult}. (C1) of \cite{Cu0} is \eqref{eq:kernel1} while (C3) of \cite{Cu0}
can be checked by elementary computation.   So Lemma \ref{lem:modulation} follows by Lemma 2.4  \cite{Cu0}.
\qed

\subsection{Choice of $p_0$ and coordinates} \label{subsec:coordinates}

Let $P_{N_g }(p) =P_{N_g(\mathcal{L}_p)}$ be the projection on
$N_g(\mathcal{L}_p) $ related to \eqref{eq:begspectdec2}. We have
\begin{equation} \label{eq:projNg} \begin{aligned} & P_{N_g }(p)X= -J\Diamond
_j \Phi _p\  \langle X ,J^{-1} \partial _{p_j}\Phi _p\rangle +\partial
_{p_j}\Phi _p \  \langle X ,\Diamond _j \Phi _p \rangle  \quad \forall \, X\in \mathcal{S} ' \end{aligned}
\end{equation} summing
on double indexes.
Indeed   range(rhs)$\subseteq N_g(\mathcal{L}_p) $, $P_{N_g }(p)J\Diamond
_j \Phi _p=J\Diamond
_j \Phi _p$ and  $P_{N_g }(p)\partial _{p_j}\Phi _p= \partial _{p_j}\Phi _p$
which follow by $\langle    \Diamond _j \Phi _p ,   \partial _{p_\ell }\Phi _p\rangle  =  \partial _{p_\ell }  p_j $. By (H4)--(H5) and \eqref{eq:kernel1}   we have  \begin{equation} \label{eq:projreg} \begin{aligned} & P_{N_g }(p)\in C^{\infty}
(\mathcal{P}, B (\mathcal{S}',\mathcal{S})). \end{aligned}
\end{equation}
We set
\begin{equation} \label{eq:projP} \begin{aligned} & P(p):=1-P_{N_g }(p). \end{aligned}
\end{equation}
We consider the vector
\begin{equation*}\label{eq:pt}
 \pi (t)=(\Pi _1 (u(t)),..., \Pi _4(u(t))) = (\Pi _1 (u(t)),\Pi _2 (u(t)),\Pi _3 (u(t)),  \Pi _4(u(0)) )  .
\end{equation*}
We  choose  $p_0, v_0, \omega _0$  such that if
$u_0$ is the initial value in  \eqref{eq:NLSP}, then
\begin{equation}\label{eq:p0}
 \Pi (\Phi _{p_0})=\Pi (u_0),  \text{ $v_0=v(p_0)$ and $\omega _0 =\omega  (p_0)$.}
\end{equation}
Furthermore, for $u_0$ satisfying \eqref{eq:sizeindata} with $\epsilon _0$   small enough,
then by
\begin{equation*}
\begin{aligned}
p_0=\Pi(\Phi_{p_0})=\Pi(u_0)=\Pi(\Phi_{p}+R)=p+\Pi(R)
\end{aligned}
\end{equation*}
we have   $|p_0-p|\lesssim   \epsilon ^2$.
Since $|p-p_1|\lesssim \epsilon$, we have
\begin{equation}\label{eq:sizeindata1}
\begin{aligned}
 & |v_0|+ |\omega _0-\omega _1|\lesssim  \epsilon .
    \end{aligned}
 \end{equation}

By  $ \pi (t) = \Pi (u(t))$,     $ v( \pi (t))$ is non constant while
$ \omega  ( \pi (t))=\omega  _0 $.
Since  $\pi (t)$ is not a constant,
it will be convenient to introduce a parameter $\pi \in \mathcal{P}$.

 By \eqref{eq:projreg}--\eqref{eq:projP} we have $(p,\pi )\to  P(p)P(\pi  )P(p_0) \in C^\infty (\mathcal {P}^2,B (\Sigma _k ,\Sigma _k))
$. Furthermore, since the hypotheses of  {Lemma 2.3} \cite{Cu0} hold here, see Lemma \ref{lem:modulation},
 we conclude that
  there exists an $a>0$ such that
 if
 $|p-p_0|<a$ and $|\pi -p_0|<a$ for  $\pi  ,p \in \mathcal{P}$, the  map
$P(p)P(\pi  )P(p_0)  $   restricts into an isomorphism
from $N_g^\perp (\mathcal{L}_{p_0}^*)\cap   \Sigma _k$  to $   N_g^\perp (\mathcal{L}_{p }^*)  \cap  \Sigma _k$
for any $k\ge - {\underline{n}}$.  The arguments in  {Lemma 2.3} \cite{Cu0} hold also with   $\Sigma _k$ replaced by
  $H^k$.

 For $X_k$ equal either to   $H^k$ or to  $ \Sigma _k$, there exists a
fixed $a>0$ such that    for   $k\ge - {\underline{n}}$ the
map,    dependent on the parameter $\pi$ s.t. $|\pi -p_0|<a$,
\begin{equation} \label{eq:coordinate} \begin{aligned} &
  \R ^4 \times \{ p: |p-p_0  |<a  \} \times ( N_g^\perp(\mathcal{L}_{p_0}^*) \cap  X_k)  \to X_k  , \\&  (\tau , p , r) \to u=   e^{J
\tau \cdot \Diamond }
  (\Phi _{p} + P(p)P(\pi ) r)
\end{aligned} \end{equation}
is for  $\| r\| _{ X _k}<a$  a local homeomorphism in the image.
 Inverting we have:
\begin{lemma} \label{lem:gradient R} We have $  r(\pi ,u) \in  C^{0} ( \{ |\pi -p_0|<a  \}\times (\U_{\underline{n}} \cap X _{k }),  X_k)$    for any $k\ge -  {\underline{n}}$,
with the $\U_{\underline{n}}$ in Lemma \ref{lem:modulation}. Furthermore, for     $k\ge -  {\underline{n}}+1$
we have $r\in
C^1(\{ |\pi -p_0|<a  \}\times (\U_{\underline{n}} \cap X _{k }), X _{k -1})$.

For $\U_{\underline{n}}$
sufficiently small in $\Sigma _{-\underline{n}} $, summing on the repeated index $j$ we have
  \begin{align}  \label{eq:gradient R1}&  \partial _ur  =   (
    P(p,\pi , p_0)  )^{-1}  P (p )
\big [ e^{ - J \tau \cdot \Diamond}    Id-
 J \Diamond _j P (p )  r \, d\tau _j-      \partial _{p_j} P (p )  r \, dp _j
 \big  ]     ,  \\&    \label{eq:gradient R2}  \partial _\pi r  = -  (P(p,\pi , p_0) )^{-1}  P (p ) (\partial _\pi P(\pi ))   r ,
\end{align} with $ (P(p,\pi , p_0)   )^{-1}:N_g^\perp
(\mathcal{L} _{p }^{\ast}  )\cap X _k\to N_g^\perp (\mathcal{L}_{p_0}^{\ast} )\cap X _k$ the inverse
of $  P(p,\pi , p_0)
=P (p )P(\pi ) P(p_0 )  :N_g^\perp (\mathcal{L} _{p _0 }^{\ast})\cap X _k\to N_g^\perp
(\mathcal{L}_{p }^{\ast})\cap X _k$.
\end{lemma}
\proof The proof is like in  Lemma 2.5  \cite{Cu0} but for \eqref{eq:gradient R2} which follows by
 \begin{equation*}
\begin{aligned} &
      0=\partial _{\pi}u=  e^{J
\tau \cdot \Diamond }
    P(p) ((\partial _{\pi}P(\pi )) r   +  P(\pi )  \partial _{\pi}r).  \qed \end{aligned}
\end{equation*}

We replace the functions $p_j$ with the functions $\Pi _j$  in the
coordinates $(\tau , p, r)$, and   move to coordinates $(\tau , \Pi , r)$.
As in  {(34) of \cite{Cu0}}, for    $\varrho =\Pi  ( r) $, we have \begin{equation}\label{eq:variables}
\begin{aligned} &
 \Pi _j =   p_j+  \varrho _j   - \Pi _j( (P _{N_g}(p)-P_{N_g}(\pi   )) r) +\langle r,
 \Diamond _j (P_{N_g}(p)-P_{N_g}(\pi  )) r\rangle   .
\end{aligned} \end{equation}
We introduce now a number of spaces.
\begin{definition}\label{def:PhaseSpace} Set  ${\Ph} _j={\Ph}^{0 }_j$
where,  for $n\in \Z$ , we consider  the spaces
\begin{equation}\label{eq:PhaseSpace}\begin{aligned} & {\Ph}^{n }_0:=  \Sigma
_{n } \cap N_g^\perp ({\mathcal H}_{\pi }) \text{ is the space of the $r$ }\\&
{\Ph}^{n }_1:=  \R ^{4} \times  {\Ph}^{n }_0 \text{ is the space of the
$(\varrho , r )$  or $(\Pi , r )$} \\&  {\Ph}^{n }_2:=  \R ^{4} \times
{\Ph}^{n }_1 \text{ is the space of the $(\tau , \varrho , r )$,  $(\tau , \Pi
, r )$ or $(\pi , \Pi , r )$} \\&  {\Ph}^{n }_3:=  \R ^{4} \times  {\Ph}^{n }_2
\text{ is the space of the $(\tau ,\Pi , \varrho , r )$  or $(\pi ,\Pi ,
\varrho , r )$}
   \\&  {\Ph}^{n }_4:=  \R ^{4} \times  {\Ph}^{n }_3
\text{ is the space of the $(\pi , \tau ,\Pi , \varrho , r )$ }  .
\end{aligned} \end{equation}
\end{definition}
Using \eqref{eq:projreg} the r.h.s. of  \eqref{eq:variables}
is smooth in $(\pi , p, \varrho , r)\in  {\Ph}^{ n }_3$ for any $n$,
with
the last two terms  $O(\| r \| _{\Sigma  _{ n}}^{2})$.  It is then easy to
see that
we can apply   the implicit function theorem
     to \eqref{eq:variables}   and get $ p_j=\Pi _j-\varrho _j + {\Psi}
_j(  \pi ,  \Pi , \varrho ,r)
  $   with  ${\Psi} _j\in C^{\infty }(  \mathcal{A},\R )$ for  $\mathcal{A}
	$ a neighborhood     of $(p_0,p_0,  0,0)$ in   $  {\Ph}^{ n }_3$
 for arbitrary $n \ge -\underline{ {n}} $, $\underline{ {n}}$ fixed in Lemma \ref{lem:modulation}.
	   It is elementary    that
  ${\Psi}_j= \mathcal{R}^{1, 2}_{\infty ,\infty}  (\pi , \Pi ,\varrho , r)
   $  and $p= \mathcal{R}^{1, 0}_{\infty ,\infty}  (\pi , \Pi ,\varrho , r)
   $, where the latter symbols are defined as follows.
				\begin{definition}\label{def:scalSymb} For $I$  an interval with 0 in the
interior,  $\mathcal{A}\subset   {\Ph}^{-n }_3$
  a neighborhood of $(p_0, p_0  , 0 ,0)$,
we  say that   $ F \in C^{m}(I\times \mathcal{A},\R)$
 is $\mathcal{R}^{i,j}_{n, m}$
 if    there exists    a $C>0$   and a smaller neighborhood
  $\mathcal{A}'$ of  $(p_0, p_0  , 0 ,0)$   in  $  {\Ph}^{-n }$
  s.t.
 \begin{equation}\label{eq:scalSymb}
  |F(t, \pi ,  \Pi ,\varrho , r )|\le C \|  r\| _{\Sigma   _{-n}}^j (\|  r\|
  _{\Sigma   _{-n}}+|\varrho | +|\Pi -\pi  |)^{i} \text{  in $I\times
  \mathcal{A}'$}.
\end{equation}
 We  will write also  $F=\mathcal{R}^{i,j}_{ n,m}$ or
 $F=\mathcal{R}^{i, j}_{ n,m} (t,\pi, \Pi ,\varrho , r)$.
We say  $F=\mathcal{R}^{i, j} _{n, \infty}$  if $F=\mathcal{R}^{i,j}_{n, l}$ for all $l\ge m$.
We say    $F=\mathcal{R}^{i, j}_{\infty, m} $       if   for all   $l\ge n$    the above   $F$ is the restriction  of an
$F \in C^{m}(I\times \mathcal{A}_{l },\R)$ with  $\mathcal{A}_l$   a neighborhood of 0 in
$  {\Ph}^{-l }$ and
$F=\mathcal{R}^{i,j}_{l, m}$. If   $F=\mathcal{R}^{i, j}_{\infty, m} $ for any $m$, we set $F=\mathcal{R}^{i, j}_{\infty, \infty} $.
\end{definition}

\begin{definition}\label{def:opSymb}  A    $T \in C^{m}(I\times
\mathcal{A},\Sigma   _{n}  (\R^3, \R ^{2 }))$,  with $I\times \mathcal{A}$
like above,
 is $ \mathbf{{S}}^{i,j}_{n,m} $   and  we  write as above
 $T= \mathbf{{S}}^{i,j}_{n,m}$  or   $T= \mathbf{{S}}^{i,j}_{n,m} (t, \pi, \Pi
 ,\varrho , r  )$,
 if      there exists  a $C>0$   and a smaller neighborhood
  $\mathcal{A}'$ of  $(p_0, p_0  , 0 ,0)$   s.t.
 \begin{equation}\label{eq:opSymb}
  \|T(t,\pi ,  \Pi ,\varrho , r )\| _{\Sigma   _{n}}\le C \|  r \| _{\Sigma
  _{-n}}^j (\|  r\| _{\Sigma   _{-n}}+|\varrho | +|\Pi -\pi  |)^{i}  \text{  in
  $I\times \mathcal{A}'$}.
\end{equation} We  use notation
$T=\mathbf{{S}}^{i,j}_{n,\infty }$, $T=\mathbf{{S}}^{i,j}_{\infty,m}$    and  $T=\mathbf{{S}}^{i,j}_{\infty,\infty}$
as
above.

\end{definition}
					
Consider $u$ s.t. $\Pi$ is close to $p_0$ and let
$\pi $ be close to $p_0$ with  $ Q   (\Phi _{\pi   }) =Q   ( \phi _{\omega _0})
$.   	   For $\omega =\omega (p)$,  $v =v (p)$, with $p=p(u)$ as of Lemma \ref{lem:modulation},
we  set
     \begin{align} \nonumber
 {K}_0(\pi ,u):&=\mathbf{E}_0 (u)- \mathbf{E}_0\left (   \Phi  _{\pi   }\right
 ) +
 \left ( \omega  + 4 ^{-1}{v^2 } \right ) \left ( Q    - Q   ( \Phi  _{\pi  })
 \right )  - v  \cdot \left ( \Pi   - \pi   \right ) \\ \label{eq:K} {K}(\pi
 ,u):&={K}_0(\pi ,u)   +2 ^{-1}{\gamma } \langle    V (\cdot - \textbf{v} t) u, u\rangle +v  \cdot \left ( \Pi
 - \pi \right )+ \mathbf{E}_0\left (    \Phi  _{\pi  }\right ) \\&
 =\mathbf{E}  (u)  +
 \left ( \omega +4 ^{-1}{v^2 }\right ) \left ( Q    - Q   (  \Phi  _{\pi  })  \right ) . \nonumber
\end{align}
 We will set $\pi = \pi (t) =\Pi (u
(t))$ for our solution $u(t)$. But, before doing this we will think of
$\pi $ as a fixed  parameter.  We   will  first  apply to
${K}_0(\pi ,u)$   a slight extension of the results of
\cite{Cu0}.  We first  consider the  \textit{strong}  symplectic form, i.e. $X\to  \Omega (X, \ ) $ defines an isomorphism from $L^2$ to its dual,
\begin{equation} \label{eq:Omega} \Omega (X,Y) =  \langle  J^{-1} X , Y \rangle
 .
\end{equation}
\begin{definition}\label{def:HamField} Given a   differentiable function $F$,
the Hamiltonian vectorfield of  $F$ with respect to a  {strong} symplectic form $\Omega $   is  the field $X_F $ such that $ \Omega (X_F ,Y)=  dF (Y)$  $\forall$
tangent vector   $Y $, with $dF$   the  Frech\'et derivative.

\noindent For    $F,G  $       differentiable functions
their   Poisson bracket is
$
  \{ F,G \}  :=  dF (X_G )
$
if $G$ is scalar valued and  $F$ is either scalar or     has values in a
Banach space $\mathbb{E}$.
\end{definition}
Equation \eqref{NLSP} can be written as $u_t=J \nabla   \mathbf{E}(u) = X_\textbf{E}(u)$.

\begin{lemma} \label{lem:PoissIdent}  Consider the  functions $\Pi _j$.  Then
$ X  _{\Pi _j }=
\frac{\partial}{\partial \tau _j} .$
In particular
  \begin{equation} \label{eq:Ham.VecField1}  \{ \Pi _j,\tau _k   \}  =-\delta _{jk} \ , \quad   \{ \Pi _j,\Pi  _k   \}=0  \ , \quad   \{ r, \Pi _j   \}  = 0.
\end{equation}
Furthermore,  for
 \begin{equation} \label{eq:thetad}  \text{$\vartheta :=-\tau _4 $,    $D_a:=\tau _a  $ for $a\le  3$ and $D:=(D_1,D_2,D_3)$}
\end{equation}
we have  at $\Pi =\pi$,
\begin{equation}\label{eq:PoissIdent1} \begin{aligned}  &  \{ \Pi _j   ,K \}  = \{ \Pi _j   ,\textbf{E} \} \, , \,  \{ r   ,K \}  = \{ r   ,\textbf{E} \} \\&    \{ \vartheta   ,K \}  = \{ \vartheta    ,\textbf{E} \} -\omega - 4^{-1}{v^2 }  \, , \,  \{ D_a   ,K \}  = \{ D_a   ,\textbf{E} \} \text{   and}
\end{aligned}
\end{equation}
\begin{equation}\label{eq:PoissIdent2} \begin{aligned} &
 \{ \Pi _j   ,\textbf{E} \} =  \{ \Pi _j   ,K_0\} + 2 ^{-1} {\gamma }    \{ \Pi _j   ,\langle  V (\cdot + \textbf{v} t) u, u\rangle\}  \, , \\&
\{ r   ,   \textbf{E} \} =  \{ r   ,   K_0\} + 2 ^{-1} {\gamma }   \{ r     ,\langle  V (\cdot + \textbf{v} t) u, u\rangle\}  \, , \\&
\{ \vartheta    ,\textbf{E} \}  =   \omega + 4^{-1}{v^2 }   + \{ \vartheta    ,   K_0\} + 2 ^{-1} {\gamma }    \{ \vartheta      ,\langle  V (\cdot +\textbf{v} t) u, u\rangle\} , \\&   \{ D_a   ,\textbf{E} \}  =v_a + \{ D_a      ,   K_0\} + 2 ^{-1} {\gamma }    \{ D_a        ,\langle  V (\cdot + \textbf{v} t) u, u\rangle\} .
\end{aligned}
\end{equation}
\end{lemma}
 \proof   \eqref{eq:PoissIdent1}--\eqref{eq:PoissIdent2} follow from  \eqref{eq:Ham.VecField1}. Like in Lemma 2.7 \cite{Cu0}  \eqref{eq:Ham.VecField1}
is a consequence of  $ X  _{\Pi _j }=
\frac{\partial}{\partial \tau _j} .$ The latter is  $X  _{\Pi _j } (u) = J\nabla \Pi _j (u)  = J \Diamond _ju=    \frac{\partial}{\partial \tau _j}  . $
\qed

 By Lemma \ref{lem:modulation}  we have $|D(0)-y_0|\lesssim \epsilon$.  We claim that \begin{equation}
\label{eq:weakint} \begin{aligned} &
 \gamma  \sup _{\text{dist} _{S^2}( \overrightarrow{{e}}, \frac{\mathbf{v}}{|\mathbf{v}|} ) \le \delta _0} \int _0^\infty  (1+|  |\mathbf{v}|\overrightarrow{{e}} t+D(0)   |^{ 2} )^{-1} dt  <10\epsilon .
\end{aligned}
\end{equation}
Indeed  let  $I_{\overrightarrow{e}}:=\{ t \ge 0 \text{ s.t. }
   | |\mathbf{v}|\overrightarrow{e}t+y_0|\ge 1/2 \}$.
For  $t\in I_{\overrightarrow{e}}$,   by $|D(0) -y_0 | <C\epsilon $
 \begin{equation*}
  \begin{aligned} &    | |\mathbf{v}|\overrightarrow{e}t+y_0| (1 -2C\epsilon )  \le |  |\mathbf{v}|\overrightarrow{e}t+D(0)  |  \le   | |\mathbf{v}|\overrightarrow{e}  t+y_0| (1 +2C\epsilon )   .
\end{aligned}
\end{equation*}
So $ \gamma \|    \langle  |\mathbf{v}|\overrightarrow{e} t+D(0)  \rangle  ^{-2}  \| _{L^1_t(I _{\overrightarrow{e}} )} \le  \gamma  (1 -2C\epsilon ) ^{-2}    \|    \langle  |\mathbf{v}|\overrightarrow{e}t+y_0  \rangle  ^{-2}  \| _{L^1_t(\R _+  )} < 2 \epsilon    .$

\noindent If     $| |\mathbf{v}|\overrightarrow{e} t+y_0|<  \frac 1 2$,  by  Chebyshev's inequality we have $|\R _+\backslash I_{\overrightarrow{e}}| <  \frac {5 \epsilon }{ 4       \gamma } $.   Then  $\gamma  \|    \langle  |\mathbf{v}|\overrightarrow{e}t+D(0)  \rangle  ^{-2}  \|   _{L^1_t(\R _+\backslash I _{\overrightarrow{e}} )}\le \gamma  |\R _+\backslash I_{\overrightarrow{e}}|   <    \frac{5}{4}    \epsilon   .$
 Adding up   we get \eqref{eq:weakint}.

\section{Darboux Theorem for fixed $\pi$} \label{sect:symplectic}

Our  aim  is to find
 in this section
Darboux coordinates   and in Sect.  \ref{sect:pullback} and \ref{sec:speccoo} an appropriate expansion of $\textbf{E}_0$. This is done near $ e^{ J \tau \cdot \Diamond}   \Phi _{\pi }$ for all $\pi $ near $p_0$
   and  is    the same of \cite{Cu0} except for the fact that in \cite{Cu0} we fix $\pi =p_0$ (uniformity
  in $\pi$ depends  on  regularity on parameters of solutions of ODE's) and that in  \cite{Cu0} we  focus only on the set
  $\Pi =\pi$ (this is dealt here  introducing slightly more general symbols in
  Def. \ref{def:scalSymb} and  \ref{def:opSymb} compared to those in Def. 2.8 and 2.9 in \cite{Cu0}). Up to minor modifications the proofs are the same of \cite{Cu0}.
 In this section, which is time independent,
the variable $t$ should not be confused with the time variable $t$ in \eqref{eq:NLSP}
or  \eqref{NLS}.

For the $P(p)$  in \eqref{eq:projP} and   $dr$ the Frech\'et derivative of $u\to r(\pi , u)$ for fixed $\pi$, see Lemma \ref{lem:gradient R},
   consider the differential forms
\begin{equation} \label{eq:Omega0} \Omega  ^{(\pi )}  (X,Y) :=  d  \tau _j
\wedge d\Pi  _j   (X,Y)+   \langle  J^{-1} P(\pi )dr  X , P(\pi )dr Y   \rangle
 .
\end{equation}

\begin{lemma}
  \label{lem:1forms}
  At the points  $e^{ J \tau \cdot \Diamond}   \Phi _{\pi }$,  for all $\tau
  \in \R  ^{4}$  we have $\Omega   ^{(\pi )}=\Omega   .$

  \noindent We have $ d\mathrm{B_0}^{(\pi )}=\Omega  ^{(\pi  )}$ and $
  d\mathrm{B }^{(\pi )} =\Omega    $ for the following  forms:
\begin{equation*} \label{eq:1forms1}\begin{aligned} &
 \mathrm{ B} _0 ^{(\pi )}:=\tau _j
d\Pi  _j  +2 ^{-1}   \langle  J^{-1} P(\pi )r   , P(\pi )dr \rangle ;
\quad \mathrm{B}^{(\pi )}:=\mathrm{B_0}^{(\pi )}+\alpha ^{(\pi )} \text{ for } \alpha
^{(\pi )}  :=   \\&
 \langle  \Gamma ^{(\pi )}(p) P(\pi )r  +
\beta  _j^{(\pi )}(p,r) P^*(p) \Diamond _j P(p)P(\pi )r ,P(\pi )dr
\rangle -\beta _j^{(\pi )}(p,r)  d\Pi _j  ,
     \\&
 \Gamma ^{(\pi )}(p):=2 ^{-1} J ^{-1}     \left ( P(p) -P(\pi  )  \right ) \ ,
\\& \beta  _j^{(\pi )}(p,r) :=2 ^{-1}\ \frac{    \langle   P^*(p)J ^{-1}P(\pi )r ,
(\partial  _{p_{ j}}P(p)) P(\pi )r    \rangle  }{ 1+   \langle  \Diamond _j
P(p)P(\pi )r , (\partial  _{p_{ j}}P(p))P(\pi )r    \rangle  } \
   .
\end{aligned} \end{equation*}

\end{lemma} \proof  This is exactly Lemma 3.1  \cite{Cu0} if we substitute $r$
with
  $R=P(\pi )r $. \qed

$\alpha ^{(\pi )}$ is used to define in Lemma \ref{lem:vectorfield0} the vector field whose flow yields the Darboux coordinates. See Sect.7 \cite{Cu2} for a brief reminder of Moser's proof of
Darboux theorem, see also below around \eqref{eq:fdarboux}.

 A direct computation,
Lemma 3.3   \cite{Cu0}, yields the following formulas.
\begin{lemma}
  \label{lem:dalpha2}     Summing on repeated indexes  the exterior differential   $d\alpha ^{(\pi  )}$
 equals
	\begin{equation*} \begin{aligned}     & \widehat{\delta } _   k    \partial
_{p_k}{\beta } ^{(\pi )}_ {j } d\Pi _j \wedge  d\Pi _k + \langle \widehat{\Gamma}_j+ (
\widehat{\delta } _   k  \partial _{p_k}{\beta }^{(\pi )} _ {j } -\widehat{\delta } _
j
  \partial _{p_j}{\beta } ^{(\pi )} _ {k } )    \Diamond _k P(p)P(\pi )r,P(\pi )dr
\rangle   \wedge \\& d\Pi _j  +   2 \langle  \Gamma ^{(\pi )}  (p)  P(\pi )dr, P(\pi
)dr\rangle + \langle \widetilde{\beta}  _j   , P(\pi )dr\rangle
 \wedge \langle  P^*(p) \Diamond _j P(p)P(\pi )r ,P(\pi )dr \rangle
  \end{aligned}
\end{equation*} 	where we have  (this time not summing on repeated indexes) 		
  \begin{equation*}
\begin{aligned}& \widehat{\delta} _k :=  (1+\langle \Diamond _k P(p)
P(\pi )r , (\partial _{p_k}P(p))P(\pi )r \rangle ) ^{-1} \  ,  \quad \widehat{\Gamma}_j :=  \\& -P^*(\pi )\nabla _r \beta _j ^{(\pi )}-\widehat{\delta} _j  [
(\partial _{p_j}\Gamma  ^{(\pi )} ) P(\pi )r +\sum _{k\le 4}  \beta ^{(\pi )} _k(\partial _{p_j}
\left (   P^*(p) \Diamond _k P(p)  \right ) )P(\pi )r  ]\\& +  \sum _{k\le 4}(
\widehat{\delta } _   k
 \partial _{p_k}{\beta } _ {j }^{(\pi )} -\widehat{\delta } _
   j    \partial _{p_j}{\beta } _ {k }^{(\pi )} )   ( P^*(p)-1)\Diamond _k P(p)  P(\pi
   )r
\\ & \widetilde{\beta}  _j := P^*(\pi )\nabla _r \beta  _j^{(\pi )}+   \widehat{\delta} _j
(\partial _{p_j}  (  \Gamma  ^{(\pi )}    +\sum _{k\le 4}
  \beta  _k^{(\pi )}  P^*(p) \Diamond _k P(p)   ) )P(\pi )r
   . \end{aligned}
\end{equation*} 	\end{lemma}

For any fixed large $n$ we have  ${\beta }^{(\pi )} _ {j } =\resto ^{0,2}_{n ,\infty}$. Then $\nabla _r \beta  _j^{(\pi )}=\mathbf{S}^{0,1}_{n,\infty} $. Similarly,
$\Gamma ^{(\pi )}(p)r=\mathbf{S}^{1,1}_{n,\infty} $  by $p=\resto ^{0,1}_{n,\infty} $, see
above  Def. \ref{def:scalSymb}. It is then easy to conclude that
$\widetilde{\beta}  _j= \resto ^{0,2}_{n,\infty}\cdot \Diamond r +\mathbf{S}^{0,1}_{n,\infty}$ and
 $\widehat{\Gamma}=\mathbf{S}^{0,1}_{n,\infty}$. Then schematically  we have
\begin{equation*} \begin{aligned}   d\alpha ^{(\pi )}   &=    \resto ^{0,2}_{n,\infty}
 d\Pi _j \wedge  d\Pi _k + \langle \resto ^{0,2}_{n,\infty}\cdot \Diamond r
 +\mathbf{S}^{0,1}_{n,\infty}
 , dr
\rangle \wedge  d\Pi _j \\& +   2 \langle  \Gamma ^{(\pi )} (p)  P(\pi )dr, P(\pi
)dr\rangle + \langle \resto ^{0,2}\cdot \Diamond r +\mathbf{S}^{0,1}_{n,\infty}
 , dr
\rangle \wedge  \langle   \Diamond r +\mathbf{S}^{0,1}_{n,\infty}
 , dr
\rangle .
  \end{aligned}
\end{equation*} We  define now $Y$   by $
 i_{Y} \Omega  ^{(\pi )}=-  \alpha ^{(\pi )}
   $   (where  $i_{Y} \Omega  ^{(\pi )}:=\Omega  ^{(\pi )} (Y, \ ) $), that is
\begin{equation*}    \begin{aligned} &
   (Y ) _{\tau _j} d\Pi _j -  (Y ) _{\Pi  _j} d\tau  _j  +\langle J^{-1} P(\pi
   )(Y ) _{r} , P(\pi ) dr \rangle
   \\ &= \beta _j ^{(\pi )}  d\Pi _j  -     \langle  \Gamma ^{(\pi )}
   P(\pi ) r +
\beta  _j ^{(\pi )} P^*(p) \Diamond _j P(p)P(\pi ) r,P(\pi ) dr   \rangle
.
   \end{aligned}
\end{equation*} This yields
\begin{equation*}      \begin{aligned} &
   (Y ) _{\tau _j} =\beta _j^{(\pi )}(p,r) =\mathcal{R}^{0,2}_{n,\infty}
      \ , \quad  (Y ) _{\Pi  _j}=0\  , \\ &
	(Y ) _{r}= (P(\pi ) P(p_0)) ^{-1}   (-  J \Gamma ^{(\pi )} (p)  P(\pi ) r
- \beta _j^{(\pi )}(p,r)  J P^*(p) \Diamond _j P(p)P(\pi ) r  )
 \\& = \mathbf{S}^{1,1}  _{n,\infty} +J\mathcal{R}^{0,2} _{n,\infty} \cdot    \Diamond r  .
   \end{aligned}
\end{equation*} Define the operator $\mathcal{K}$   by $ i_{X}d
\alpha ^{(\pi )} =
 i_{\mathcal{K} X} \Omega  ^{(\pi )}$. For
 vectors $X$ with $(X ) _{\Pi _j}=0$  (where  $(X)_\kappa :=d \kappa X $  for any coordinate $\kappa$) we have
\begin{equation*} \begin{aligned}&   (\mathcal{K} X)_{\tau _j} =
  \langle       \mathcal{R}^{0,2} _{n,\infty} \cdot    \Diamond r
+ \mathbf{S}^{0,1}_{n,\infty}, (X)_{r} \rangle    , \quad  (\mathcal{K} X)_{\Pi _j} = 0,
     , \quad  (\mathcal{K} X)_{r} =\\& = (P ( \pi )P ( {p}_0))^{-1}J   [   2
\Gamma ^{(\pi )}  (X)_{r} + \langle       \mathcal{R}^{0,2}_{n,\infty}  \cdot    \Diamond r +
\mathbf{S}^{0,1}_{n,\infty}, (X)_{r} \rangle    ] (      \Diamond r +
\mathbf{S}^{0,1}_{n,\infty}  )   . \end{aligned} \end{equation*}

System $
 i_{\mathcal{X}^{t \pi}} \Omega _t ^{(\pi )} =- \alpha ^{(\pi )}
   $ is equivalent to $(1+t \mathcal{K})\mathcal{X}^{t \pi} =Y. $
We have:
\begin{lemma}
  \label{lem:vectorfield0}   For any fixed large $n$ and for $\varepsilon
  _0>0$
	consider the subset  $\U _{1}$ of $H^1$ defined by 	 $| \Pi  -p_0|\le
\varepsilon _0$,  $  \|r\| _{\Sigma  _{-n}}\le   \varepsilon
  _0$ and $| \Pi
(R)|\le   \varepsilon _0$.  Assume   	$| \pi  -p_0|\le   \varepsilon
  _0$.
  Then for $ \varepsilon
  _0 $  small enough
     $\exists$ a unique vectorfield
		$ \mathcal{X}^{t \pi}  :\U  _{ {1}}  \to
 \Ph _2  $       s.t.  $
 i_{\mathcal{X}^{t \pi}} \Omega _t ^{(\pi )} =- \alpha ^{(\pi )}
   $, where $\Omega _t^{(\pi )}:=\Omega  ^{(\pi )}+t(\Omega -\Omega ^{(\pi)})$
   for  $|t|<4$.
   For  $|t|<3$, we have
\begin{equation}\label{eq:quasilin1} \begin{aligned} &   (\mathcal{X}^{ t \pi
})_\Pi =0 \, , \,
 (\mathcal{X}^{ t \pi })_{\tau _j} =  \mathcal{T}_j
 \, , \,  (\mathcal{X}^{ t \pi })_r = J \mathcal{A }_j \Diamond _j r  +
 \mathcal{D}    \text{ with} \end{aligned}
\end{equation}  \begin{equation}\label{eq:quasilin2} \begin{aligned} &
\mathcal{T}_j,
    \mathcal{A }_j=\mathcal{R}^{0,2}_{n,\infty}
    (t, \pi , \Pi , \Pi (r),r)   \text{ and }  \mathcal{D }
    =\mathbf{S}^{1,1}_{n,\infty }(t, \pi , \Pi , \Pi (r),r).   \end{aligned}
\end{equation}
\end{lemma} \proof   The proof is like Sect. 3.1
\cite{Cu0}. It rests in solving $(1+ t\mathcal{K})\mathcal{X}^{ t,\pi }=Y$.
 It is enough to focus on  $(\mathcal{X}^{ t,\pi })_r $. By
 $  \Gamma  ^{(\pi )}=\langle \mathbf{S}^{1,0}_{n,\infty} ,\  \rangle \mathbf{S}^{0,0}_{n,\infty}$
 we have  schematically
 \begin{equation}\label{eq:quasilin3}
\begin{aligned} &    (\mathcal{X}^{ t \pi } )_r +t
 \langle  \mathcal{R}^{0,2} _{n,\infty} \cdot
 \Diamond r + \mathbf{S}^{1,0}_{n,\infty}, (\mathcal{X}^{ t \pi } )_r\rangle
  (J\mathcal{R}^{0,2}_{n,\infty}  \cdot    \Diamond r + \mathbf{S}^{0,0} _{n,\infty})\\& =
  J\mathcal{R}^{0,2} _{n,\infty} \cdot    \Diamond r + \mathbf{S}^{1,1}_{n,\infty} .   \end{aligned}
\end{equation}
 Applying to both sides $ \langle  \mathcal{R}^{0,2}_{n,\infty}  \cdot
 \Diamond r + \mathbf{S}^{1,0}_{n,\infty}, \ \rangle $  we get
$
	\langle  \mathcal{R}^{0,2}_{n,\infty}
 \cdot
 \Diamond r + \mathbf{S}^{1,0}_{n,\infty}, (\mathcal{X}^{ t,\pi } )_r\rangle
 =\mathcal{R}^{1,2}_{n,\infty}$.
 Plugging back in \eqref{eq:quasilin3},
  we get $(\mathcal{X}^{ t,\pi } )_r= J\mathcal{R}^{0,2} _{n,\infty} \cdot    \Diamond r +
  \mathbf{S}^{1,1}_{n,\infty}.$
 \qed

  We   need   information on the
flows generated by systems like   \eqref{eq:quasilin1}--\eqref{eq:quasilin2}.
  Lemma \ref{lem:ODE} is
 a simple extension of   the standard  arguments in Lemma 3.8 \cite{Cu0}.

\begin{lemma} \label{lem:ODE} For
$n,M,M_0,\textbf{s},\textbf{s}',k,l\in \N\cup \{ 0 \}$ with $1\le l\le M$,
consider
\begin{equation} \label{eq:ODE}\begin{aligned} &
  \dot \tau   _j (t)= T_j (t,\pi ,\Pi , \Pi (r), r  ) \  , \quad \dot \Pi  _j
  (t)   =0 \  ,  \\&  \dot R (t)
	= J \mathcal{A }_j(t,\pi ,\Pi , \Pi (r), r   )  \Diamond _j r  +
\mathcal{D}(t,\pi ,\Pi , \Pi (r), r   ) , \end{aligned}   \end{equation}
  where  we assume what follows:
			\begin{itemize} \item[(i)] $P_{N  _g ( p_0  )} (\mathcal{A }_j(t,\pi ,\Pi ,\varrho ,r) J\Diamond _j
r  + \mathcal{D}(t,\pi ,\Pi ,\varrho ,r)) \equiv 0$;

\item[(ii)]   for $|s|<5$,  $|\pi -p_0| <\widetilde{\varepsilon}
  _0$,   $|\Pi -p_0|
    <\widetilde{\varepsilon} _0$ and $\| r\| _{\Sigma _{-n} \cap N^\perp _g ({\mathcal
    H}^*_{p_0} )} <\widetilde{\varepsilon}
  _0$
  we have
$\mathcal{A }, T=\resto ^{0,M_0+1}_{n,M}(t,\pi ,\Pi , \Pi (r), r  ) $
   and
  $ \mathcal{D}  =\textbf{S} ^{i,M_0 }_{n,M}(t,\pi ,\Pi , \Pi (r), r  )$.

\end{itemize} Let  $k\in \Z\cap [0,n-(l+1) ]$ and set for
$\textbf{s}^{\prime \prime}\ge 1$ \begin{equation} \label{eq:domain0}
\begin{aligned}   &   \U _{\varepsilon _1,k}^{\textbf{s}^{\prime \prime}}
:= \{ (\tau, \Pi  ,r) \in  { {\Ph}}^{\textbf{s}^{\prime \prime}}_{2} \
 : \  |\Pi -p_0|\le \varepsilon _1 \  ,  \  \|  r \| _{\Sigma _{-k }} + |\Pi
 (r)|
  \le \varepsilon _1\}.
\end{aligned}    \end{equation} Then for $\varepsilon _1>0$ small enough,  \eqref{eq:ODE} 	 defines a      flow
$\mathfrak{F} _{t \pi } =(  \mathfrak{F}
  _{t \pi } ^{\tau    } , \mathfrak{F}  _{t \pi } ^{\Pi   }
   ,  \mathfrak{F}  _{t \pi } ^{r} )$
	  for $t\in [-2, 2]$   in   $\U _{\varepsilon _1,k}^{1}$
 with $\mathfrak{F}  _{t \pi }  ^{\Pi   } =\mathfrak{F}  _{0 \pi }  ^{\Pi   }$
 such that, see \eqref{eq:PhaseSpace},
\begin{equation} \label{eq:ODE-1}\begin{aligned} & \mathfrak{F}   _{t\pi}   ^{\tau } (
  \Pi ,\varrho,  r) =\tau +  \resto  ^{0,M_0+1 }_{  n- l-1 , l}(t,\pi , \Pi , \varrho , r ) \text{ and}     \end{aligned}   \end{equation}
\begin{equation} \label{eq:ODE1}\begin{aligned} & \mathfrak{F} _{t\pi}   ^{r} (
  \Pi ,\Pi (r), r) = e^{J q(t,\pi , \Pi , \Pi (r), r)\cdot \Diamond } ( r+
\textbf{S} (t,\pi , \Pi ,  \Pi (r),r    ))        \end{aligned}   \end{equation}
 in $(-2,2)\times  \mathcal{V}$, $\mathcal{V}\subset \Ph _3 ^{l+1 -n}= \{  (\pi,  \Pi, \varrho ,r)\}$
 a neighborhood of $(p_0,p_0, 0,0 )$,
\begin{equation} \label{eq:ODEpr210}\begin{aligned}   &  q  =
 \resto  ^{0,M_0+1 }_{  n-  l-1, l}(t,\pi , \Pi , \varrho , r )
    \, , \quad  \textbf{S}  =  \textbf{S} ^{i,M_0 }_{  n- l-1,l}
    (t,\pi , \Pi , \varrho , r ) .
\end{aligned}    \end{equation}

\noindent We have also  $\textbf{S}= \textbf{S}_1+ \textbf{S} _2$  with
\begin{equation} \label{eq:duhamel}   \begin{aligned}   &  \textbf{S}_1 (t, \pi
, \Pi , \Pi (r), r    ) = \int _0^t \mathcal{D }(t',\pi , \Pi ,\Pi (r (t') ) ,r
(t')    ) dt'\\&
    \textbf{S}_2 (t, \pi , \Pi , \varrho, r   )
=\textbf{S} ^{i,2M_0+1 }_{  n- l-1, l}(t,\pi , \Pi , \varrho , r ).
\end{aligned}    \end{equation}
\begin{equation} \label{eq:ODEindex}   \begin{aligned}   &   \text{For $n- l-1 \ge  \textbf{s}' \ge  \textbf{s}+l \ge l $ and    $k\in
\Z\cap [0,n- l-1 ]$}
\end{aligned}    \end{equation}
and
     for  $ \varepsilon _1>0$ sufficiently small, we have for $ {\Ph}  ^{\textbf{s} }_2
     =\{ ( \tau ,\Pi ,r)\}$
 \begin{equation} \label{eq:reg1}\begin{aligned} &\mathfrak{F}  _{t \pi}
  \in C^l((-4,4)\times \{ \pi : |\pi -p_0| <\widetilde{\varepsilon} _0\} \times \U _{\varepsilon _1,k}^{\textbf{s}'}
  ,  {\Ph}  ^{\textbf{s} }_2
 )
. \end{aligned}   \end{equation} Furthermore, $\exists$ $ \varepsilon _2>0$
s.t. $\mathfrak{F}  _{t \pi}( \U _{\varepsilon _2,k}^{\textbf{s}'})
\subset  \U _{\varepsilon _1,k}^{\textbf{s}'}$ for all $|t|< 4$ and $|\pi -p_0| <\widetilde{\varepsilon} _0$.

\noindent  The constants  implicit in \eqref{eq:ODE-1}, \eqref{eq:ODEpr210}--\eqref{eq:duhamel}
 and the constants $\varepsilon _1$ and $\varepsilon _2$
can be chosen to depend only on the constants from (ii).
   Finally, we have      \begin{equation}
\label{eq:ODE11}\begin{aligned} & \mathfrak{F}  _{t \pi }   (e^{J \tau \cdot
\Diamond } U) = e^{J \tau \cdot \Diamond }\mathfrak{F}  _{t \pi}   (  U) .
\end{aligned}   \end{equation} 	\end{lemma}

\proof
 Set $S= e^{-Jq \cdot \Diamond }r$ for $q\in \R ^{4}$.  Then  consider the following system:
 \begin{align} \nonumber &\dot \tau    = T  (t,\pi ,\Pi , \varrho , P(p_0)e^{ Jq \cdot \Diamond }S  )    , \  \dot \Pi    =0 \  ,
    \dot S
	= e^{-Jq \cdot \Diamond } \mathcal{D}(  t,\pi ,\Pi ,\varrho ,P(p_0)e^{ Jq \cdot \Diamond }S   )   ,
\\& \label{eq:ODEpr0}\dot q = \mathcal{A } (t,\pi ,\Pi ,\varrho ,P(p_0)e^{ Jq \cdot \Diamond }S    ) \quad , \quad q(0)=0  ,  \\& \dot \varrho _j= \langle S,  e^{-Jq \cdot \Diamond } \Diamond _j\mathcal{D}( t,\pi ,\Pi , \varrho ,P(p_0) e^{ Jq \cdot \Diamond }S   ) \rangle  \ . \nonumber
\end{align}
 We have  $P_{N  _g ( p_0  )}e^{ Jq \cdot \Diamond }S=0 $ $\forall$ $t$.
Indeed this holds  for $t=0$ by  $q=0$ and  $S=r\in N  _g ^\perp ( p_0  )$.
For $t>0$  it holds  by Gronwall inequality since, by
Hypothesis   (i), we have
\begin{equation*}  \begin{aligned} &
   \partial _t  P_{N  _g ( p_0  )}e^{ Jq \cdot \Diamond }S= P_{N  _g ( p_0  )} [ J \mathcal{A }       \cdot \Diamond  P_{N  _g ( p_0  )}  e^{ Jq \cdot \Diamond }S  +       J \mathcal{A }       \Diamond P(p_0) e^{ Jq \cdot \Diamond }S +   \mathcal{D}    ]
        \\& = P_{N  _g ( p_0  )}  J \mathcal{A }     \cdot \Diamond  P_{N  _g ( p_0  )}  e^{ Jq \cdot \Diamond }S  \, ,  \text{ where }  \mathcal{A }=\mathcal{A } (t,\pi ,\Pi ,\varrho ,P(p_0)e^{ Jq \cdot \Diamond }S)   .
\end{aligned}   \end{equation*}

 The field in \eqref{eq:ODEpr0} is $C^l( (-3,3)\times \U _{-k},  \Sigma _{ \mathbf{s}^{\prime \prime}}
\times \R ^{12})$  for
  $    k, \mathbf{s}^{\prime \prime}\in [0, n-(l+1) ]$  and    $ l\le M$, and    with
$\U _{-k}\subset \R ^{16}\times \Sigma _{-k}  =\{ (\pi    ,\Pi, q , \varrho , S ) \}$ a neighborhood of  $(p_0, p_0, 0, 0, 0)$.
 For $l\ge 1$ we can apply  to \eqref{eq:ODEpr0} standard
theory of ODE's to conclude that  there is a $\U _{-k}'$ like $\U _{-k}$
such that   the flow  is of the form
\begin{equation} \label{eq:ODEpr1}\begin{aligned}
   &  \tau (t)
	=  \tau + \textbf{T} (t, \pi ,\Pi ,\varrho , r )  \, , \quad \Pi (t)=\Pi  \, ,   \\&  S(t)
	=  r+ \textbf{S} (t, \pi ,\Pi ,\varrho , r )  \ , \quad  \textbf{S} (0, \varrho ,r) =0\ ,  \\& q(t)=  {q} (t,\pi ,\Pi ,\varrho , r ) \ , \quad  {q}(0,\pi ,\Pi , \varrho ,r) =0\ ,  \\& \varrho (t)=\varrho + \overline{\varrho} (t, \pi ,\Pi ,\varrho , r)  \   ,   \quad   \overline{\varrho} (0,\pi ,\Pi ,\varrho , r ) =0\ ,
\end{aligned}   \end{equation}
\begin{equation} \label{eq:ODEpr2}\begin{aligned} \text{with   }  &
     \textbf{S}  \in C^l((-2,2)
\times \U _{-k}' , \Sigma _{n- (l+1) }
 )    \\&  \textbf{T}, \overline{\varrho},{q} (t,\pi ,\Pi ,\varrho , R )  \in C^l((-2,2)  \times   \U _{-k}', \R ^{ 4}
 ).
\end{aligned}   \end{equation}

 \noindent For $S\in \Sigma  _{ 1}  \cap B_{\Sigma _{-k}}$ and $S(0)=S$,
choosing $\textbf{s}^{\prime \prime}\ge 1$  we have $S(t)\in\Sigma  _ {1}  $
with $\Pi (S(t))=\varrho (t)$  for  $\varrho (0)=\varrho =\Pi (S)$.  Then \eqref{eq:ODEpr2} yields
\eqref{eq:ODEpr210}.

 \noindent  $\exists$ fixed $C$ s.t. for all $|t|\le 2$  and for $r(0)=r$ and $\varrho (0) = \varrho  $ we have
\begin{align} \label{eq:gronwall0} & \| r(t) \| _{\Sigma _{ s ^{\prime\prime}} }\le   C\| r  \| _{\Sigma _{ s ^{\prime\prime}} } \, , \\&    | \varrho (t) - \varrho    |    \le  C  \| r \| _{\Sigma _{ l+1  -n}}^{ M_0+1 }
(|\pi - \Pi| +|\varrho | +\| r \| _{\Sigma _{ l+1   -n}}) ^i  . \label{eq:gronwall1}
\end{align}
 Indeed \eqref{eq:gronwall0}--\eqref{eq:gronwall1}  can be obtained as in \cite{Cu0} by a standard
use of  Hypothesis (ii) and Gronwall inequality.   Using     \eqref{eq:gronwall0}--\eqref{eq:gronwall1}
it is elementary to get the bounds \eqref{eq:scalSymb}--\eqref{eq:opSymb}
 of Definitions \ref{def:scalSymb} and \ref{def:opSymb} needed to get \eqref{eq:ODE-1} and
\eqref{eq:ODEpr210}.  For this and the routine  proof of the rest of this lemma, see Lemma 3.8 \cite{Cu0}.\qed

\noindent Classically     Darboux Theorem  follows by $ i_{\mathcal{X}
^{t \pi }} \Omega  _t^{(\pi  )}=-\alpha ^{(\pi )}$  and by
\begin{equation}\label{eq:fdarboux}  \begin{aligned} & \partial _t
 \mathfrak{F}_{ t  \pi}^*\Omega  _t^{(\pi  )}        =   \mathfrak{F}_{ t
 \pi}^*   (L_{\mathcal{ X} ^{ t  \pi}}\Omega  _t^{(\pi
)}+\partial _t\Omega  _t^{(\pi  )}  )     =
 \mathfrak{F}_{ t  \pi}^*
  ( d i_{\mathcal{X}  ^{ t  \pi}}\Omega  _t^{(\pi  )}+d\alpha ^{(\pi )}  )
=0  \end{aligned} \end{equation}
with $L_X$ the Lie derivative, whose definition is   not needed   here.
Since this  $\mathfrak{F}_{ t  \pi}$ is not a differentiable flow on any given manifold\footnote{For this reason (3.44)   \cite{bambusi}  is formal.   Also (3.42)   \cite{bambusi} is formal, see  Remark 6.6 \cite{Cu0}. Notice that the concise   Remark 17  \cite{bambusi} is in \cite{bambusi}     a surrogate  for
     part of the proof of Darboux Theorem.},  \eqref{eq:fdarboux} is
  formal. Still,     Sect. 3.3 and   Sect. 7
\cite{Cu0} (i.e.   a regularization and a limit argument for $  \mathfrak{F}  _{t \pi}$) yield:

\begin{lemma}[Darboux Theorem]\label{lem:darboux} Consider \eqref{eq:ODE} defined by the field
${\mathcal X}^{t \pi }$ and consider indexes and notation of Lemma
\ref{lem:ODE} (in particular   $M_0=1$ and $i=1$; $n$ and $M$ can be
arbitrary). Consider $\textbf{s}'$,$\textbf{s}$ and $k$ as in \ref{lem:ODE}.
Then  for   $\mathfrak{F}  _{1 \pi} \in C^l( \{ \pi :|\pi -p_0|<\varepsilon
_1 \} \times  \U _ {\varepsilon
_1,k}^{\textbf{s}^{\prime  }}
  ,  {\Ph}  ^{\textbf{s} }_3)$
derived from \eqref{eq:reg1},
 we have $ \mathfrak{F}   _{1 \pi  }^* \Omega =\Omega  ^{(\pi )}$  for all fixed $\pi$. \qed
\end{lemma}

The  symbols  $\resto ^{i,j}_{k,
m}$ and   $\textbf{S} ^{i,j}_{k,
m}$ are preserved by the flows  in Lemma \ref{lem:ODE}:
\begin{lemma}
  \label{lem:ODEdomains} Assume Lemma
  \ref{lem:ODE}.  	Consider the  flow of system    \eqref{eq:ODEpr0}.
	  Then for $  \textbf{S}  =  \textbf{S} ^{i,M_0 }_{  n- l-1 ,l}$   and
$q,\overline{\varrho}= \resto  ^{i,M_0+1 }_{  n-(l+1), l}$  	   \begin{equation} \label{eq:ODEdomains1} \begin{aligned} &
\mathfrak{F}  _{t \pi}  ^{r} ( \Pi , \varrho  ,r ) = e^{Jq(t
  , \pi , \Pi ,\varrho  ,r)\cdot \Diamond } ( r+ \textbf{S} (t
    , \pi , \Pi ,\varrho  ,r    )) \\&  \mathfrak{F}  _{t \pi}  ^{\varrho}
    (\Pi,   \varrho  ,r ) =
 \varrho + \overline{\varrho} (t  , \pi , \Pi ,\varrho  ,r)
\text{ with:} \end{aligned}   \end{equation}

\begin{itemize} \item[(1)]  $\mathfrak{F}  _{t \pi} \in C^l((-2,2)\times \U _{-k  }
  , \Ph ^{-h }_3
 )
$ for $\U _{-k } $ a neighborhood  of  $(p_0, p_0, 0, 0)$ in $  \Ph ^{-k
}_3=\{ (\pi , \Pi , \varrho , r)\}$,
  $k \in \Z \cap [0,n-(l+1) ]$ and $h\ge \max\{ k
+l ,(2 l+1) -n\}$;  \item[(2)]   $\forall$  $\resto ^{a,b}_{h, l}(
\pi , \Pi ,\varrho  ,r)$  (resp. $\textbf{S}^{a,b}_{h, l}(   \pi , \Pi ,\varrho  ,r)$)
 and $(h,k)$ as in (1)  with $h\le n-  l-1  $,
$\resto ^{a,b}_{h, l}\circ \mathfrak{F}
_{s \pi}=\resto ^{a,b}_{k, l}( s  , \pi , \Pi ,\varrho  ,r) $ (resp.
$\textbf{S} ^{a,b}_{h, l}\circ \mathfrak{F}  _{s \pi}=\textbf{S} ^{a,b}_{k,
l} ( s  , \pi , \Pi ,\varrho  ,r)$).

\end{itemize} 	\end{lemma} \proof   The proof  is  the same of that
 of Lemma 3.9  in \cite{Cu0}. \qed

\section{Pullback   of $K_0(\pi , u)$} \label{sect:pullback}

  We seek  a  coordinate  system
	where   $u\to K_0(\pi , u)$    has a  simple effective part. We will then treat
$2 ^{-1} \langle V u, u\rangle $ as a perturbation.   The first step is the
following.
 \begin{lemma} \label{lem:ODE1}   Consider
  $\mathfrak{F}^\pi =\mathfrak{F}_1^\pi \circ \cdots \circ \mathfrak{F}_L^\pi$
  with  $   \mathfrak{F}_j^\pi = \mathfrak{F}_ {j,t=1,\pi}$ transformations as
  of Lemma \ref{lem:ODE}. Suppose that for $j$
      we have $M_0=m_j$, with given numbers  $1\le m_1\le ...\le m_L$.
      Suppose also that all the $j$   have the same pair   $ n$
      and $ M $, which we assume sufficiently large.
      Let $i_j=1$ if $m_j=1$. Fix $0<m'<M$.

      \begin{itemize}
\item[(1)]
      Let  $ n >  2L(m'+1)    + s'_L > 4L(m'+1)  + s_1$, $s_1\ge 1$.
	Then,  for any $\varepsilon >0$ there exists a $\delta  >0$ such that
$\mathfrak{F}^\pi \in C^{m'} (\{|\pi -p_0|< \varepsilon
  _0  \} \times \U
^{s'_L}  _{\delta  ,a}, \U ^{s _1}  _{\varepsilon ,h})$    for $0\le a\le h
$ and $0\le h < r-(m'+1)  $ for a fixed $\varepsilon
  _0>0$ small enough.

\item[(2)]
      Let    $ n >    2L(m'+1)  )    +h>  4L(m'+1)     + a$,
       $a\ge 0$. The above composition, interpreting
       the $\mathfrak{F}_j^\pi$'s as maps in the $(\Pi ,\varrho , r)$
				variables  as in Lemma \ref{lem:ODEdomains}, 				 yields also
$\mathfrak{F}^\pi    \in C^{m'}( \{  \U _{-a}
  , \Ph ^{-h }_2
 )
$ for  $\U_  {-a} $ a sufficiently small neighborhood of  $(p_0, p_0,0,0)$
in $  \Ph ^{-a}_3=\{(\pi ,\tau , \Pi , \varrho , r)\} $ and with $  \Ph
^{-h}_2=\{( \tau , \Pi , \varrho , r)\} $.

\item[(3)] $\exists$ a neighborhood   $\U_  {-a}  $ of  $(p_0, p_0,0,0)$ in  $
    \Ph ^{-a}_2 =\{( \pi , \Pi , \varrho , r)\} $
 and   functions
$\mathcal{R}^{i, j}  _{a,m'}\in C^{m'} ( \U_  {-a}  ,\R ) $ and
 $\mathbf{S}^{i, j}  _{a,m'}\in C ^{m'}(  \U_  {-a} ,\Sigma _{ a}) $ s.t.
\begin{equation} \label{eq:ODE2}\begin{aligned} & \Pi (r'):= \Pi (r)\circ
  \mathfrak{F}^\pi = \Pi (r)
 +\mathcal{R}^{i_1, m_1+1}  _{a,m'}(\pi , \Pi, \Pi (r), r)  , \\&
  p':=p\circ \mathfrak{F}^\pi    = p +\mathcal{R}^{i_1, m_1+1}
   _{a,m'}(\pi , \Pi, \Pi (r), r),\\& \Phi  _{p'} =
     \Phi  _{p}  +\mathbf{S}^{i_1, m_1+1}  _{a,m'}(\pi , \Pi, \Pi (r), r)
     .
\end{aligned}   \end{equation} \item[(4)] If
 $F(e^{J\tau \cdot \Diamond }U)\equiv F(U)$   $\exists$ functions
$ \mathbf{S}^{i,j }_{k  ,l} =\mathbf{S}^{i,j }_{k  ,l}
  (\pi , \Pi, \Pi (r), r)$ s.t.
\end{itemize} \begin{equation*} \label{eq:ODE4}\begin{aligned} &
  F\circ \mathfrak{F}^\pi (U)= F\left ( \Phi  _{p}+P(p)
  (r+  \mathbf{S}^{i_1,m_1 }_{b+1  ,m '})
  + \textbf{S}^{i_1,m_1 +1 }_{b-1,m'} \right )   , \, b= n-4L(m'+2) .
\end{aligned}   \end{equation*}
\end{lemma} \proof The proof is  almost  the same of Lemma 4.1 \cite{Cu0}.
\qed

The following lemma  is  partially proved inside  Lemma 4.1 \cite{Cu0}.
\begin{lemma}
  \label{lem:transf} Consider the setup of Lemma \ref{lem:ODE1}.  In particular,  let $0<m'<M$.
	Set    $a'_{L }:=n- 2L m' $ .    Consider the coordinates  decomposition
		$\mathfrak{F}^\pi= ( \mathfrak{F}^\pi _\tau , \mathfrak{F}^\pi _\Pi , \mathfrak{F}^\pi _\varrho,\mathfrak{F}^\pi _r)$.
		 Then, we have
	\begin{equation}\label{eq:tranR}\begin{aligned} & \mathfrak{F}^\pi _r (\pi , \Pi , \varrho ,r)=e^{Jq (\pi , \Pi , \varrho ,r) \cdot \Diamond } ( r+ \textbf{S}^{i_1,m_1  }_{a'_{L } ,m'} (\pi , \Pi ,\varrho ,r) )  \ , \\&  q (\pi , \Pi ,\varrho ,r)  =\resto ^{0,m_1 +1  }_{a'_{L } ,m'} (\pi , \Pi ,\varrho ,r) \ ,   \\&  \mathfrak{F}^\pi _\tau  (\pi , \Pi , \varrho ,r) -\tau =\resto ^{0,m_1 +1  }_{a'_{L } ,m'} (\pi , \Pi ,\varrho ,r)    .\end{aligned}
\end{equation}
	\end{lemma}\proof   Here $\mathfrak{F}^\pi  =\mathfrak{F}^\pi _1 \circ \cdots \circ \mathfrak{F}^\pi _L$.
Lemma \ref{lem:transf}  follows by  Lemma \ref{lem:ODE}  when $L=1$. By induction let us suppose that Lemma \ref{lem:transf}  is true for  $\mathfrak{F} ':= \mathfrak{F}^\pi _1 \circ \cdots \circ \mathfrak{F}^\pi  _{L-1} $ with $a'_{L-1 }  $.
Then, omitting the arguments  $(\pi , \Pi , \varrho ,r)$,
we have
\begin{equation*}\begin{aligned} &
 \mathfrak{F}^\pi _r=e^{J(q' \circ\mathfrak{F}^\pi  _{L} )  \cdot \Diamond } \left ( e^{Jq_L  \cdot \Diamond } ( r+ \textbf{S}^{i_L,m_L  }_{n-(m'+1) ,m'}  )+ \textbf{S}^{i_1,m_1  }_{a'_{L-1} ,m'}\circ\mathfrak{F}^\pi  _{L}  \right ) \\&  = e^{J ({q'\circ\mathfrak{F}^\pi   _{L}+q_L} )  \cdot \Diamond }  \left  ( r+ \textbf{S}^{i_L,m_L  }_{n-(m'+1) ,m'}  )+e^{-Jq_L  \cdot \Diamond }\textbf{S}^{i_1,m_1  }_{a'_{L-1}-  m'  ,m'}  \right ),
\end{aligned}
\end{equation*}
where $\textbf{S}^{i_1,m_1  }_{a'_{L-1} ,m'}\circ\mathfrak{F}^\pi  _{L} =\textbf{S}^{i_1,m_1  }_{a'_{L-1}-  m'  ,m'} $
by  the last claim in Lemma \ref{lem:ODEdomains}. Since   $e^{-Jq_L  \cdot \Diamond }\textbf{S}^{i_1,m_1  }_{a'_{L-1}- m'  ,m'} =\textbf{S}^{i_1,m_1  }_{a'_{L-1}-2m' ,m'}  $  we conclude the first equality  in  \eqref{eq:tranR} for $q:=q' \circ\mathfrak{F} _{L}+
q_L$.  We have $q_L=\resto ^{0,m_L+1  }_{n-(m'+1) ,m'}$  by \eqref{eq:ODEpr210}.  We have  by Lemma
\ref{lem:ODEdomains}
\begin{equation*}\begin{aligned} &
 q' \circ\mathfrak{F}^\pi _{L}= \resto ^{0,m_1 +1  }_{a'_{L-1 } ,m'} \circ\mathfrak{F}^\pi  _{L}=
 \resto ^{0,m_1 +1  }_{a'_{L-1 } -m' ,m'}
\end{aligned}
\end{equation*}
and the latter is    an $ \resto ^{0,m_1 +1  }_{a'_{L  }  ,m'}$.  Similarly, by induction and   Lemmas
\ref{lem:ODE} and \ref{lem:ODEdomains} we have the following, which completes the proof:
\begin{equation*}\begin{aligned} &  \mathfrak{F}_\tau ^\pi:=
 \tau  \circ \mathfrak{F}^\pi  =  (\tau + \resto ^{0,m_1 +1  }_{a'_{L-1 } ,m'}) \circ
\mathfrak{F}^\pi   _{L}= \tau +  \underbrace{\resto ^{0,m_L +1  }_{n-  m'-1   ,m'}  +\resto ^{0,m_1 +1  }_{a'_{L-1 } -m' ,m'}}_{\resto ^{0,m_1 +1  }_{a'_{L  }  ,m'}} . \qed
\end{aligned}
\end{equation*}

The following lemma is a slight generalization   of   Lemmas 4.3 and
4.4 \cite{Cu0}  (it reduces to these lemmas for $\Pi =\pi =p_0$) and can be proved similarly.

\begin{lemma}
  \label{lem:back}  Consider a  transformation
  $\mathfrak{F} =\mathfrak{F}^\pi  _1 \circ \cdots \circ \mathfrak{F}^\pi  _L$
  like in Lemma \ref{lem:ODE1}  and with $m_1=1$,
   with same notations, hypotheses and conclusions.
	In particular we suppose $n$ and $M$ sufficiently large that
 the conclusions of  Lemma \ref{lem:ODE1}
	hold for preassigned sufficiently large $s=s'_L$, $k'=a$ and $m' $.
   Then
	   there are  a  $k$ as large as we want but $k\ll k'$, an $m$ as large as we
want but $m\ll m'$ 		and an  $\varepsilon >0$ sufficiently small,
 such that in
			$\{ |\pi -p_0|<  \varepsilon \}\times \U^{s}_{\varepsilon ,k}$   we have the expansion
\begin{align*}  \label{eq:back1}   & K _0\circ \mathfrak{F}=
 \psi (\pi , \Pi ,\Pi  (P(\pi )r))
 +     {2}^{-1}\Omega (  \mathcal{L} _{\pi   } P(\pi )r, P(\pi )r )
  +  \resto ^{1,2} _{k,m}   +\mathbf{E}_P(r)+\textbf{R}',
\\&   \textbf{R} ^{\prime \prime } =     \sum _{j=2}^4 \langle B_{j } (\pi ,
\Pi,   r,\Pi (P(\pi )r) ), r  ^{   j} \rangle
      +\int _{\mathbb{R}^3}
B_5 (x,  \pi , \Pi , r, r(x),\Pi (P(\pi )r) )   r^{   5}(x) dx \nonumber \end{align*}
with:
\begin{itemize}
  \item  for $d(p) $ as in  \eqref{eq:dp}, $\psi  (\pi , \Pi ,\Pi  (P(\pi )r))$ is defined by
 \begin{equation} \label{eq:psi}\begin{aligned}  & \psi  (\pi , \Pi ,\Pi  (P(\pi )r)) = d (\Pi -  \Pi  (P(\pi )r) )-
	d(\pi ) \\& + \left ( \lambda  (\Pi -  \Pi  (P(\pi )r) )-
	\lambda  (\pi )  \right )  \cdot   \pi  ;
\end{aligned}\end{equation}

\item $  {\resto} ^{1,2} _{k,m} =\resto ^{1,2} _{k,m}   (\pi , \Pi ,\Pi
    (r), r) $;

\item $ r ^j(x)$  represent $j-$products of components of $r$;

  \item
$B_{j }( \pi , \Pi, r,\varrho  ) \in C^{m } ( \U _{-k}, \Sigma _k
(\mathbb{R}^3, B   (
 (\mathbb{R}^{2   })^{\otimes d},\mathbb{R} ))) $
   for $2\le d \le 4$ with $\U _{-k} \subset \Ph ^{-k}_3$ a
   neighborhood of  $(p_0, p_0, 0, 0)$;

   \item  we have  $ B_2(\pi ,  \pi ,0 , 0  )=0;
$

 \item  for all \footnote{ By an oversight, in the statement of Lemma 4.3 \cite{Cu0}  it is required also   $|\zeta  |\le \varepsilon$ for some small $\varepsilon >0$, but in fact the proof of Lemma 4.3 \cite{Cu0} shows that this constraint is unnecessary.}
$ \zeta \in \mathbb{R}^{2   }$
 and $(\pi , \Pi, \varrho ,r) \in\U _{-k} $
 we have  for $i\le m$ and for    fixed constants  $C_{i,L}$
\begin{equation} \label{eq:B5}\begin{aligned} &  \| \nabla _{\pi , \Pi
,r,\zeta, \varrho  } ^iB_5(  \pi , \Pi,  r,\zeta  ,\varrho  ) \| _{\Sigma
_k(\mathbb{R}^3,   B   (
 (\mathbb{R}^{2   })^{\otimes 5},\mathbb{R} )} \le C_{i,L}  .
 \end{aligned}  \end{equation}
\end{itemize}

\end{lemma}
\qed

\section{Spectral coordinates   associated to $\mathcal{H} _{\pi}$}
\label{sec:speccoo}
We now consider the complexification of $L^2(\R ^3, \R ^2)$ and think of     $\mathcal{L}_p$, defined in \eqref{eq:linearizationL}, and $J$
as   operators  in  $L^2(\R ^3, \C ^2)$. We consider
the bilinear map      $\langle \cdot , \cdot \rangle $  defined in  $L^2(\R ^3, \R ^2) $  by
formula   \eqref{eq:hermitian} and its natural extension as bilinear form defined in
 $L^2(\R ^3, \C ^2)$.  We similarly extend   $\Omega   (\cdot , \cdot ) = \langle J^{-1}\cdot , \cdot \rangle$. We denote
the extension by the    same notation.
 We     set
\begin{equation}\label{eq:linearization1}     \text{   $\mathcal{H}_p:=\im \mathcal{L}_p$  with   $ \mathcal {H}_\omega := \mathcal {H}_p$   when $v(p)=0$ and $\omega (p)=\omega$.}
\end{equation}
By  \eqref{eq:conjNLS} we have
\begin{equation}\label{eq:conjNLS1}
 \begin{aligned}
\ {\mathcal H}_p =e^{-\frac 12 J v (p)\cdot x}
 {\mathcal H}_{\omega (p)} e^{\frac 12 J v (p)\cdot x}
\end{aligned}
\end{equation}
In the course of the proof $\pi$ is a parameter that is later set to be equal to $\pi =\pi (t)=\Pi (t)$.
Since $ \Pi (0)=p_0$  and  $ \Pi _4(t)=Q(u(t))=\Pi _4(0)$ for all $t$,  in the notation of
\eqref{eq:LagrMult}  and of \eqref{eq:p0} we have $\omega (\Pi (t))=\omega _0$. By \eqref{eq:conjNLS1}
we have  $ {\mathcal H}_{\pi } =e^{-\frac 12 J v(\pi )\cdot x}
 {\mathcal H}_{\omega _0} e^{\frac 12 J v(\pi )\cdot x}$ for all $\pi =\pi (t)$ so that in particular the spectrum is constant in $t$ with $\sigma ( {\mathcal H}_{\pi }) =  \sigma ( {\mathcal H}_{\omega _0 })  $.   For ${\mathcal H}_{\omega _0 }$ we assume the hypotheses (H7)--(H10) in Sect. \ref{subsec:statement}.
  We have
\begin{equation}\label{eq:Homega}
 \begin{aligned} &
 {\mathcal H}_{\omega  }  =   \im  J  (-\Delta  + \omega ) +  \im  J
  \begin{pmatrix}     \beta  (\phi _\omega ^2 )+2 \beta  '(\phi _\omega ^2 ) \phi_\omega^2  &
0  \\
0 &     \beta  (\phi _\omega ^2 )
 \end{pmatrix}  .
\end{aligned}
\end{equation}

\noindent  We can see the  $\sigma_e ({\mathcal H}_{\omega  }) =(-\infty , -\omega ]\cup [\omega , \infty )$ by the following:
\begin{align} &  \label{eq:Homega1}M^{-1}
 {\mathcal H}_{\omega  } M =  {\mathcal K}_{\omega  } ,
 \\&
 {\mathcal K}_{\omega  }:=
   \sigma _3  (-\Delta  + \omega ) +
 \begin{pmatrix}     \beta  (\phi _\omega ^2 )+  \beta  '(\phi _\omega ^2 ) \phi_\omega^2  &
   \beta  '(\phi _\omega ^2 ) \phi_\omega^2   \\
  - \beta  '(\phi _\omega ^2 ) \phi_\omega^2  &  -\beta  (\phi _\omega ^2 )-  \beta  '(\phi _\omega ^2 ) \phi_\omega^2
 \end{pmatrix} \nonumber \\&
 M:=
  \begin{pmatrix}   1  &
1  \\
-\im  &   \im
 \end{pmatrix}   \, , \quad   M^{-1} =\frac{1}{2}
  \begin{pmatrix}   1  &
\im   \\
1  &   -\im
 \end{pmatrix}   \, , \quad   \sigma _3=\begin{pmatrix} 1 & 0\\0 & -1 \end{pmatrix}
   .  \nonumber
\end{align}

\noindent Notice that  we have $M^T=2 \overline{M}^{-1}$, with $M^T$ the
transpose of $M$ and $ \overline{M}$ obtained taking complex
conjugation of the entries of $M$. We will later use the equality
 \begin{equation} \label{eq:Homega2}\begin{aligned}&   {M}^{-1}J   {M} =
  - \overline{M}^{-1}J  \overline{M}=
    \overline{M}^{-1}J ^{-1} \overline{M} =-\im \sigma _3.   \end{aligned}
\end{equation}
 \begin{equation} \label{eq:defe}\begin{aligned}&   \text{Set  $\textbf{e} =(\textbf{e} _{1}, ... ,\textbf{e} _{\textbf{n}} )$ with
  $\textbf{e} _{j}:=\textbf{e} _{j}(\omega _0)$ for the $\textbf{e} _{j}(\omega  )$ in (H7)
  Sect. \ref{subsec:statement}.}   \end{aligned}
\end{equation}
We state without proof the following elementary well known  lemma. \begin{lemma}
  \label{lem:Specdec} $\phi _{\omega _0}(x)>0$ $\forall$ $x$ and (H4)--(H7)
  yield the spectral decomposition
\begin{align}  \label{eq:spectraldecomp} & N_g^\perp (\mathcal{H}_{\pi }
^*)  =  \big (\oplus _{\mu  \in \sigma _p\backslash \{ 0\}}
\ker (\mathcal{H}_{\pi}  - \mu
 ) \big) \oplus L^2_c ( {\pi} )\\& \nonumber  L^2_c ( {\pi} ):=
\left\{N_g(\mathcal{H}_{\pi } ^\ast)\oplus \big (\oplus _{\mu    \in \sigma
_p\backslash \{ 0\}}   \ker (\mathcal{H}_{\pi } ^*- \mu
 ) \big)\right\} ^\perp ,
\end{align}
 where orthogonality is with respect to the  bilinear form
$\langle \cdot , \cdot \rangle $     in  $L^2(\R ^3, \C ^2)$.
\end{lemma}
\noindent Notice that the $N_g^\perp (\mathcal{H}_{p }
^*) $ in \eqref{eq:spectraldecomp} equals $N_g^\perp
(\mathcal{L}_p^{\ast}) \otimes _\R  \C$ for the  $N_g^\perp
(\mathcal{L}_p^{\ast}) $ in \eqref{eq:begspectdec2}.
\begin{lemma}
  \label{lem:basis}  It is possible to choose   eigenfunctions
   $\xi  ^{(\pi )}_j\in  \ker (\mathcal{H}_{\pi} -\textbf{e} _j )$
so that  $ \Omega ( \xi _j ^{(\pi )} ,\overline{\xi} _k ^{(\pi )})=0$ for $j\neq k$ and  $ \Omega
( \xi _j  ^{(\pi )},\overline{\xi}  _j^{(\pi )})= \im s_j $ with  $s_j \in \{ 1,-1\}$.  We have $
\Omega ( \xi _j ^{(\pi )} , {\xi} _k ^{(\pi )})=0$ for all $j$ and $k$.
 We have $ \Omega (\xi ^{(\pi )}_j , F )=  \Omega (\overline{\xi} ^{(\pi )}_j , F )=   0$ for any
$\xi ^{(\pi )}_j   $ and any $F\in L^2_c ( {\pi} )$.
  Furthermore,  we have $s_j=1$ for all $j$. \end{lemma} \proof
 Everything is elementary and, except for the last sentence, is proved in
 Lemma 5.2 \cite{Cu0}.
Hypotheses (H5) and (H6), Lemma \ref{lem:back} and formula \eqref{eq:H2}
below,
 imply that $s_j=1$ for all $j$.
Notice that by \eqref{eq:conjNLS} it is sufficient to consider the case $v(\pi )=0$
getting functions $\xi _j$, and then set $\xi _j  ^{(\pi )}= e^{-\frac 12 J v(\pi )\cdot x}\xi _j   $.

\qed

For any $r\in N^\perp _g({\mathcal H} _{p_0 }  ^*)$
with $r=\overline{r}$
we consider     $z\in \C ^{\textbf{n}}$ and   $ f\in L_c^2 ({p_0 } )$ s. t.  \begin{equation}
  \label{eq:decomp2}
  P(\pi )  r  =\sum _{j=1}^{\mathbf{n}}z_j \xi
  _j  ^{(\pi )} +
\sum _{j=1}^{\mathbf{n}}\overline{z}_j\overline{\xi  } _{  j}^{(\pi )}  + P_c(\pi )f , \end{equation}
where $ P_c(\pi )$ is the projection on $L^2_c(\pi )$ associated to \eqref{eq:spectraldecomp}.
Set
\begin{equation}
  \label{eq:Pd}
  P_d(\pi ) :=1- P_c(\pi ). \end{equation}
We write  also for $\ell =c,d$
 \begin{equation}\label{eq:decomp22}
  P_\ell (\omega  )    =  P_\ell (p )     \text{  if $\omega =\omega (p)$ and $v(p)=0$}, \end{equation}
  with an abuse of notation.
	The following equality holds true:
 \begin{equation}\label{eq:decomp23}
  e^{ J \frac{v (p) \cdot x }2}    P_c(p )    e^{ -J \frac{v (p) \cdot x }2}  =
P_c(\omega  (p) )    . \end{equation}

\noindent We have  $f=\overline{f}$. By the first part of
Lemma \ref{lem:basis}  we have \begin{equation}
  \label{eq:H2}  \begin{aligned} &
	 \frac  1 {2}\Omega (  \mathcal{L} _{\pi   } P(\pi )r, P(\pi )r ) =
	-\frac \im 2
   \Omega  (  \mathcal{H} _\pi P(\pi )r  , P(\pi )r)   \\& = \sum _{j=1}^{\mathbf{n}}s_j \textbf{e}_j  |z_j|^2   -\frac \im 2
   \Omega  (   \mathcal{H} _\pi  P_c(\pi ) f  , P_c(\pi )  f )=:H_2(\pi , z,f). \end{aligned}
\end{equation} Notice that $  \Omega  ( -\im  \mathcal{H} _\pi  P_c(\pi ) f  , P_c(\pi ) f )>c \|  f\| _{H^1}^2$ for a
$c>0$
 by Theorem 2.11 \cite{CPV}.   So   the condition $s_j=1$ for all $j$
is related to the fact that,   as shown  in   \cite{GSS1,GSS2,W2},  under  Hypotheses (H5)--(H6) the ground
 states $ \Phi _p$  are local minimizers of the functional $\textbf{E}_0$
under the constraint $\Pi =p$.

The  Fr\'echet derivatives of  the functions $r $  and  $f$ are mutually related  by
\begin{equation*}
    P_c(\pi )dr =\sum _{1\le j\le \mathbf{n}}(dz_j \xi _j+d\overline{z}_j \overline{\xi} _j)
    +P_c(\pi )d f .
\end{equation*}
  In these coordinates we have
	\begin{equation}
  \label{eq:OmegaCoo}\begin{aligned} & \Omega ^{(\pi )}= \sum _{1\le j\le 4}d\tau _j \wedge  d\Pi
  _j +
   \Omega  (   P_c(\pi )dr, P_c(\pi )dr)  \\& =   \sum _{1\le j\le 4}d\tau _j \wedge  d\Pi  _j  +   \im
   \sum _{1\le j\le \mathbf{n}} dz_j\wedge d \overline{z}_j
+  \Omega (  P_c(\pi )  d f , P_c(\pi ) d f )  .  \end{aligned}\end{equation} 			 For a function $F$ and a fixed $\pi $  the Hamiltonian  field
of $F$ w.r.t. $\Omega ^{(\pi )}
 $  is
\begin{equation*}\begin{aligned}   &
    X_F ^{\pi}= \sum _{1\le j\le 4}  ( (X_F^{\pi})_{\Pi _j} \frac{\partial}{\partial \Pi _j}+ (X_F^{\pi})_{\tau _j} \frac{\partial}{\partial \tau _j}  )   \\&  +\sum _{1\le j\le \mathbf{n}}  ( (X_F^{\pi})_{z_j} \xi _j (x) +
 (X_F^{\pi})_{\overline{z}_j}\overline{\xi  }_j( x )  )
+ P_c(\pi ) (X_F^{\pi})_{f}, \quad   (X_F^{\pi})_{f} \in L^2_c (p_0),\end{aligned}
\end{equation*}
where $ \xi _j=\frac{\partial}{\partial z _j} $ and $ \overline{\xi } _j=\frac{\partial}{\partial \overline{z} _j} $.
By $i_{X_F ^{(\pi)}}\Omega ^{(\pi)} =dF$ and by
\begin{equation*} \begin{aligned} &
 dF= \partial _{\Pi _j}F d\Pi _j+\partial _{\tau_j}F d\tau_j+\partial _{z_j}F dz_j+\partial _{\overline{z}_j}F d\overline{z}_j+   \langle \nabla _fF, d f\ \rangle  \\&
i_{X_F ^{(\pi)}}\Omega ^{(\pi)}=    (X_F^{\pi })_{\tau_j}  d\Pi _j-    (X_F^{\pi})_{\Pi _j}  d\tau_j + \im  (X_F^{\pi})_{z_j}  d\overline{z}_j\\& - \im  (X_F^{\pi})_{\overline{z}_j}  dz_j+   \langle  J ^{-1} P_c(\pi ) (X_F^{\pi})_{f}, P_c(\pi )d f\ \rangle ,
\end{aligned}\end{equation*}
we get
\begin{equation*}  \label{eq:hamvf1}
\begin{aligned} &   (X_F^{\pi})_{\tau_j}=\partial _{\Pi _j}F\ , \  (X_F^{\pi})_{\Pi _j}=-\partial _{\tau _j}F \ , \
(X_F^{\pi})_{z_j}=-\im  \partial _{\overline{z}_j}F\ , \  (X_F^{t})_{\overline{z}_j}= \im    \partial _{z_j}F
\end{aligned}\end{equation*}
and by $ P^*_c(p_0 ) P^*_c(\pi )J ^{-1} P_c(\pi ) (X_F^{\pi})_{f}=\nabla _fF $ and $P^*_c(p)J ^{-1}=J ^{-1}P_c (p)$
we get
\begin{equation*}  \label{eq:hamvf2}   (X_F^{\pi})_{f}= (P _c(p_0 ) P_c (\pi )P _c (p_0 ) )^{-1}    J \nabla _fF
. \end{equation*}

\noindent This implies that the Poisson bracket w.r.t. $\Omega ^{(\pi )}$ can be expressed as
\begin{equation} \label{eq:poiss}\begin{aligned} & \{  F,G  \} ^{(\pi)} :=dF(X_G^\pi) = \partial _{\tau  _j}F \partial _{\Pi _j}G- \partial _{\Pi _j}F \partial _{\tau  _j}G\\& -\im \partial _{z_j}F\partial _{\overline{z}_j}G
+\im \partial _{\overline{z}_j}F\partial _{z_j}G +  \langle  \nabla _fF, (P _c(p_0 ) P_c (\pi )P _c (p_0 ) )^{-1}J \nabla _fG \rangle.
\end{aligned}\end{equation}
We have the following result, which, apart from Claim (3)  follows from Lemma \ref{lem:ODE} (Claim (3) follows from Sect. 3.3 \cite{Cu0}, see Lemma 5.3  \cite{Cu0}).
\begin{lemma}
  \label{lem:chi} For  $i\in \{  0,1\}$ fixed and $n,M\in \N$
sufficiently large let
	  \begin{equation*} \label{eq:chi1}\chi   =\sum _{|\mu +\nu |=M_0 +1}
b_{\mu\nu} (\pi ,\Pi , \Pi (f))  z^{\mu} \overline{z}^{\nu} +\im \sum _{|\mu +\nu
|=M_0 } z^{\mu} \overline{z}^{\nu}
 \langle J    B_{\mu   \nu
}(\pi ,\Pi ,\Pi (f))
  , f \rangle
\end{equation*} with $ b_{\mu\nu}(\pi , \Pi , \varrho)= \resto   ^{i,0}
_{n,M}(\pi , \Pi ,\varrho)$ and $ B_{\mu\nu}(\pi , \Pi ,\varrho)=  \textbf{S}
^{i,0}  _{n,M}(\pi , \Pi ,\varrho)$ and with   \begin{equation} \label{eq:symm}
\overline{b}_{\mu\nu}  = {b}_{\nu\mu}   \   , \quad   \overline{B}_{\mu   \nu }
=-B_{\nu\mu} , \end{equation} (so that $\chi$ is real valued for
$f=\overline{f}$).     Then we have what follows.
\begin{itemize} \item[(1)]
 Consider the vectorfield $X_\chi ^\pi  $ defined
with respect to $\Omega ^{(\pi )}$.
 Then, summing on repeated indexes (with the equalities defining the
 field $X_\chi ^{cl}$), we have:
 \begin{equation*}   \begin{aligned} &(X_\chi ^\pi ) _{z_j}= -\im \partial
 _{\overline{z}_j}
 \chi =: (X_\chi ^{\pi cl} ) _{z_j} \,  , \quad (X_\chi ^\pi ) _{\overline{z}_j}=
  \im \partial _{z_j}
 \chi  =: (X_\chi ^{\pi  cl} ) _{\overline{z}_j} \, , \\&  (X_\chi ^\pi ) _{f}
=  (P _c(p_0 ) P_c (\pi )P _c (p_0 ) )^{-1}\partial  _{\Pi _j(f)}\chi \, P_c^*(\pi )
     J \Diamond _j f +
     (X_\chi ^{\pi cl} ) _{f}  \\&\text{ where }  (X_\chi ^{\pi cl} ) _{f}:= (P _c(p_0 ) P_c (\pi )P _c (p_0 ) )^{-1}z^{\mu}
     \overline{z}^{\nu}B_{\mu   \nu
}(\pi , \Pi ,\Pi (f)).
 \end{aligned}
  \end{equation*}
\item[(2)] For $\phi ^{t,\pi }$   the flow of  $X_\chi ^\pi$, see
    Lemma \ref{lem:ODE}, and      $(z^t,f^t)=  (z,f)\circ \phi ^{t,\pi
    }$,   \begin{equation} \label{eq:quasilin51}
    \begin{aligned} &    z^t  =   z  +
 \mathcal{Z}(t,\pi , \Pi ,\Pi (f),  z, f)  \, , \\
  & f^t  =e^{Jq(t,\pi , \Pi ,\Pi (f),  z, f )\cdot \Diamond } ( f+ \textbf{S}(t,\pi , \Pi ,\Pi (f),  z, f ) )
\end{aligned} \end{equation} where, for a preassigned pair $(k,m)$  and sufficiently small neighborhoods of 0,  $
B_{\Sigma  _{-k}}$ in $ \Sigma
_{-k}\cap X_c(\pi ) $  and   $B_{\C ^{\mathbf{n} }} $ (resp.$B_{\R
^{4}}$)
   in $\C ^{\mathbf{n} }$ (resp.$ {\R ^{4}}$)
and for
$(\pi , \Pi ,\varrho ,  z, f)\in \mathcal{A}:= (B_{\R^{4}} +p_0)^2\times B_{\R ^{4}}\times
     B_{\C ^{\mathbf{n} }}  \times    B_{\Sigma  _{-k}}$, we have
\begin{equation} \label{eq:ODEpr21}\begin{aligned} &
     \textbf{S} \in C^m((-2,2)\times
     \mathcal{A} \, , \,   {q}   \in C^m((-2,2)\times    \mathcal{A} , \R ^{ 4}
 ) \, , \\&   \mathcal{Z }  \in C^m((-2,2) \times \mathcal{A} ,\C ^{\mathbf{n} })
\end{aligned}   \end{equation} with for fixed $C$
 \begin{equation} \label{eq:symbol11}   \begin{aligned}   &   | q (t, \pi , \Pi
  ,\varrho , z,f  ) |\le C (|z|+\| f\| _{\Sigma _{-k}}) ^{M_0+1} \\&
  |\mathcal{Z} (t, \pi , \Pi
  ,\varrho , z,f  ) |+ \| \textbf{S} (t, \pi , \Pi
  ,\varrho , z,f  ) \| _{\Sigma _{k}} \le C (|z|+\| f\| _{\Sigma _{-k}}) ^{M_0
  } .
\end{aligned}    \end{equation} We have $ \textbf{S} (t, \pi , \Pi
  ,\varrho , z,f  )= \textbf{S}_1 (t, \pi , \Pi
  ,\varrho , z,f  )+ \textbf{S}_2 (t, \pi , \Pi
  ,\varrho , z,f  )$ with
\begin{align}  \nonumber    \textbf{S}_1 (t, \pi , &  \Pi
  ,\varrho , z,f ) =\int _0^s
 (X_\chi ^{\pi st} ) _{f}\circ \phi ^{t',\pi } dt' ,  \
  \| \textbf{S}_2 (t, \pi , \Pi
  ,\varrho , z,f  ) \| _{\Sigma _{k}}\le  \\& \le C
  (|z|+\| f\| _{\Sigma _{-k}}) ^{2M_0+1 }(|z|+\| f\| _{\Sigma
  _{-k}}+|\varrho | +|\pi - \Pi |)^i . \label{eq:symbol12}
\end{align}

\item[(3)] For   $s,s',k$ as in
    Lemma \ref{lem:ODE},     $\phi  ^{t ,\pi } \in C^l(  \U
    ^{s'}_{\varepsilon _1,k}
  ,  { \Ph}^{s }_1)$
 satisfies
   $\phi ^{t,\pi \ast  }\Omega  ^{(\pi)}=\Omega ^{(\pi)}$
  { in } $C^{\infty}(   \U ^{s'}_{\varepsilon _2,k}
  , B  (  { \Ph}^{s' }_1 \otimes  { \Ph}^{s' }_1, \R )
 )  $ for $\varepsilon  _2>0$ sufficiently small.
\end{itemize}
\end{lemma}
Lemma \ref{lem:back} expressed in terms of the splitting \eqref{eq:decomp2}
yields the following result. The proof is the same of  that of Lemma 5.4 \cite{Cu0}.
\begin{lemma}
  \label{lem:ExpH11}  Consider
  $\mathfrak{F}  =\mathfrak{F}_1 ^{\pi}\circ \cdots \circ \mathfrak{F}_L^{\pi}$
  like in Lemma \ref{lem:ODE1},   with $m_1=2$ and for fixed $n$ and $M$
  sufficiently large.
	Let $(k,m)$ be a fixed pair.
     Consider decomposition
  \eqref{eq:decomp2}.  Then on a domain $\U^s_{\varepsilon , k} $
  like
  \eqref{eq:domain0}  and for $|\pi -p_0|< \varepsilon $
\begin{equation}  \label{eq:ExpH11} \begin{aligned} &  {K}_0\circ \mathfrak{F}=  {
\psi} (\pi, \Pi ,\Pi ( f)) +H_2'   +\textbf{R} \end{aligned} \end{equation}
where  $  {{{\psi}}} ( \pi, \Pi ,\varrho )$  is $C^{m}$
with  $\psi $ defined by \eqref{eq:psi} and with   what follows.\begin{itemize}
\item[(1)]
 We have
\begin{equation}  \label{eq:ExpH2} H_2 '= \sum _{\substack{ |\mu +\nu |=2\\
\mathbf{e}  \cdot (\mu -\nu )=0}}
 g_{\mu \nu} ( \pi , \Pi , \Pi  (f) )  z^\mu
\overline{z}^\nu   - \frac{\im }{2} \langle J^{-1} \mathcal{H}_{\pi
 } P_c(\pi ) f,  P_c(\pi )  f\rangle .
\end{equation}

\item[(2)] We have $\textbf{R}  = \sum _{j=-1}^{4} {\textbf{R} _{ j }} + \resto ^{1,2}_{k,m}(\pi , \Pi ,\Pi (f), f) $, with  for $\varrho =\Pi (f) $:
\begin{equation*}   \begin{aligned} &{ {\mathbf{R}} _{  - 1 }}=
 \sum _{\substack{ |\mu +\nu |=2\\
\mathbf{e} \cdot (\mu -\nu )\neq 0  }} g_{\mu \nu } (\pi, \Pi ,\varrho )z^\mu
\overline{z}^\nu  +\im \sum _{|\mu +\nu |  = 1} z^\mu \overline{z}^\nu  \langle
J  G_{\mu \nu }(\pi ,\Pi , \varrho) ),f\rangle  ;\end{aligned}
\end{equation*} For $ {N}$ as  in (H7)  and  for $\varrho =\Pi (f) $,
    \begin{equation*}    \begin{aligned} &   {\textbf{R} _0}=   \sum _{|\mu
    +\nu |= 3}^{2 {N}+1} z^\mu
\overline{z}^\nu   g_{\mu \nu }( \pi , \Pi , \varrho ) ;  \quad {\textbf{R} _1}= \im
\sum _{|\mu +\nu |= 2}^{2 {N} } z^\mu \overline{z}^\nu \langle J
G_{\mu \nu }(\pi , \Pi , \varrho ), f\rangle ; \\ &   {\textbf{R} _2}=
\langle \mathbf{B}_{2 } ( \pi , \Pi   ,\varrho ),   f  ^{   2} \rangle
     \text{ with $\mathbf{B}_{2 } (\pi , \pi , 0)=0$}
\end{aligned}\nonumber \end{equation*}
  where   $f^d(x)$    represents schematically $d-$products of components of
  $f$;
\begin{equation*}    \begin{aligned} &   {\textbf{R} _3}=   \sum _{ \substack{
|\mu +\nu |=\\=  2N+2}} z^\mu \overline{z}^\nu   g_{\mu \nu }( \pi , \Pi   ,
z,f,\varrho )        +\sum _{ \substack{ |\mu +\nu |=\\= 2N+1}} z^\mu
\overline{z}^\nu  \langle J    G_{\mu \nu }(\pi , \Pi   , z,f,  \varrho ),
f\rangle  ;\\ &   {\textbf{R} _4}=  \sum
_{d=2}^4 \langle B_{d } ( \pi , \Pi   ,  z ,f,\varrho ),   f  ^{   d} \rangle
       +\int _{\mathbb{R}^3}
B_5 (x, \pi , \Pi   , z ,f, f(x),\varrho )  f^{   5}(x) dx\\& +
\widehat{\textbf{R}} _2(   z ,\pi , \Pi   , f,\varrho )+   E_P (  f) \text{ with
$B_2(0,0,\varrho )=0$.} \end{aligned}\nonumber \end{equation*}

\item[(3)] For $\delta _j:=( \delta _{1j},
    ..., \delta _{mj}),$ \begin{equation}  \label{eq:ExpHcoeff1}
    \begin{aligned} & g_{\mu \nu } (\pi , \pi   , 0 ) =0 \text{ for $|\mu +\nu
    | = 2$  with $(\mu , \nu )\neq (\delta _j, \delta _j)$ for all $j$,} \\&
    g_{\delta _j \delta _j } ( \pi , \pi   ,0 ) =\mathbf{e} _j (\omega (\pi ))  , \quad
    G_{\mu \nu }( \pi , \pi   , 0 ) =0 \text{ for $|\mu +\nu | = 1$ }.
    \end{aligned} \end{equation} The  $g_{\mu \nu } ( \pi , \Pi   ,\varrho )$
    and $G_{\mu \nu }(  \pi , \Pi   ,\varrho )$ are $C ^{m }$ in $( \pi , \Pi   ,\varrho )$
    with $G_{\mu \nu }(  \pi , \Pi   ,\varrho ) \in C^m  (
    \mathrm{V},\Sigma _k (\mathbb{R}^3,\mathbb{C}^{2N  }))$,
  for   a
small neighborhood $\mathrm{V}$ of $( p_0,p_ 0,0 )$ in $\R ^{12}   $
 and they satisfy
symmetries analogous to \eqref{eq:symm}.

\item[(4)] We have $g_{\mu \nu }( \pi ,\Pi , z,  \varrho   ), g_{\mu \nu }( \pi ,\Pi , z, f, \varrho   )  \in C^{ m }(  \mathbb{{U}},
 \mathbb{C} ) $   for   a
small neighborhood $\mathbb{{U}}$ of $( p_0,p_ 0,0, 0, 0)$ in $\R ^{12} \times  \mathbb{C}^{
\mathbf{n}}\times (\Sigma _{-k}\cap X_c(\pi))  =\{ (\pi , \Pi , \varrho , z,f  )\}$.

\item[(5)] $G_{\mu \nu }( \cdot , \pi , \Pi   , z,  \varrho ) , \in C^{ m } (
    \mathbb{{U}}, \Sigma _k (\mathbb{R}^3,\mathbb{C}^{2   }))) $.

\item[(6)] $B_{d }( \cdot ,\pi , \Pi   ,  z,f,\varrho  ) \in C^{m } (
    \mathbb{{U}}, \Sigma _k (\mathbb{R}^3, B   (
 (\mathbb{C}^{ 2})^{\otimes d},\mathbb{R} ))) $, for $2\le d \le 4$.  We have also
$\mathbf{B}_{2 }( \pi , \Pi ,   \varrho  )\in C^{m } (
    \mathbb{{U}}, \Sigma _k (\mathbb{R}^3, B   (
 (\mathbb{C}^{ 2})^{\otimes 2},\mathbb{R} ))) $.

\item[(7)] Let $ \zeta \in \mathbb{C}^{2   }$. Then for
  $B_5(\cdot ,  \pi , \Pi   ,z ,f, \zeta  ,\varrho )$   we have   for fixed constants $C_l$ (the
  derivatives are not in the holomorphic sense)
\begin{equation} \label{H5power2}\begin{aligned} &\text{for  $|l|\le  m $ , }
\| \nabla _{  z ,f,\zeta, \varrho  } ^lB_5(  \pi , \Pi   ,z,f,\zeta  ,\varrho
) \| _{\Sigma _k(\mathbb{R}^3,   B   (
 (\mathbb{R}^{2   })^{\otimes 5},\mathbb{R} )} \le C_l .
 \end{aligned}\nonumber \end{equation}
 \item[(8)] We have $\widehat{\textbf{R}} _2
 \in C^{m} ( \mathbb{{U}} ,\C
 )$ and
\begin{equation*} \begin{aligned} &    | \widehat{\textbf{R}} _2 (\pi , \Pi   ,z ,f, \varrho
)|  \le C (|\Pi - \pi |+ |z|+  \| f \| _{\Sigma _{-k}}) \| f \| _{\Sigma
_{-k}}^2. \qquad \qed\end{aligned}\end{equation*}
\end{itemize}
\end{lemma}

\begin{definition}[Normal Forms]
\label{def:normal form} A function $Z( z,f, \varrho ,\pi , \Pi )$ is in normal form if  $
  Z=Z_0+Z_1
$
where $Z_0$ and $Z_1$ are finite sums of the following type:
\begin{equation}
\label{e.12a}Z_1=  \im \sum _{\mathbf{e}    \cdot(\mu-\nu)\in \sigma _e(\mathcal{H} _{p_0})}
z^\mu \overline{z}^\nu \langle  J  G_{\mu \nu}( \pi , \Pi , \varrho   ),f\rangle
\end{equation}
with $G_{\mu \nu}( \cdot  ,\pi , \Pi ,\varrho )\in  C^{m} (  \widetilde{U},\Sigma _{k }(\R ^3, \C ^{2 }))$
 for  fixed $k,m\in \N$, with $\widetilde{U} =\{ p: |p-p_0|< a  \}^2 \times U$   and
 $U\subseteq \R ^{4}$ a neighborhood of 0;
\begin{equation}
\label{e.12c}Z_0= \sum _{  \mathbf{e}  \cdot(\mu-\nu)=0} g_{\mu   \nu}
( \pi , \Pi ,  \varrho  )z^\mu \overline{z}^\nu
\end{equation}
and $g_{\mu   \nu}  \in  C^{m} (   \widetilde{U},
\mathbb{C})$. For $\mathbf{e}$ see  \eqref{eq:defe}.
We assume furthermore that the above coefficients
satisfy the symmetries in \eqref{eq:symm}: that is $\overline{g}_{\mu \nu}=g_{\nu\mu  }$ and $\overline{G}_{\mu \nu}=-G_{\nu \mu }$.
 \end{definition}

 \noindent

  \begin{remark}
\label{lem:normfor} Notice that it is easy to see that  there exists
an $\delta _0>0$ such that for $|\pi - p_0|<\delta _0$  the formula $\omega (\pi )     \cdot(\mu-\nu)\in \sigma _e(\mathcal{H} _{\pi })$ for  $|\mu
    +\nu |\le   2N+2$  is independent of $\pi $. Similarly   independent of $\pi $
    is  the equality   $ \mathbf{e} (\omega (\pi ))  \cdot(\mu-\nu)=0$ for  $|\mu
    +\nu |\le  2N+2$.   To see this  pick   $l_1=1<l_ 2<...<l_{k  }\le \textbf{n}$  and $\l _{k+1}=\textbf{n}+1$
       so that  $\mathbf{e}_j (\omega )=\mathbf{e}_{i} (\omega ) $ if and only if $j,i \in [l_a, l_{a+1})$  for some $a$. Notice that by Hypothesis (H7) the $l_j$ do not depend on $\omega$.
    Now we have
   \begin{equation*}  \begin{aligned} &
  \textbf{e} (\omega ) \cdot (\mu -\nu ) = \sum _{\ell =1}^{k }\widehat{\mu} _\ell \textbf{e}_{l_\ell } (\omega ) =0 \text{  with }
  \widehat{\mu} _\ell:=
   \sum  _{ l_\ell \le j <l_{\ell  +1}   }
	(\mu _j- \nu _j)    .
\end{aligned}
\end{equation*}
       By Hypothesis (H9) we conclude that $\widehat{\mu} _\ell =0$  for all $\ell$. But this is a
       condition which does not depend on $\omega$.
 \end{remark}

 In view of Remark \ref{lem:normfor}, Lemma \ref{th:main} below is almost the same of Theorem 6.4
 \cite{Cu0} except for the  two points discussed at the beginning of Sect. \ref{sect:symplectic}: the fact that the expansion is  for all $\pi $ near $p_0$
 rather than just for one such $\pi$; the fact that we are not restricting $\Pi$ to $\Pi =\pi$. With minor adjustments, in view also of the fact that   $\Pi$ is invariant by our coordinate changes,
 the proof  of  Lemma \ref{th:main}, which we only sketch,  is almost the same of the proof   of Theorem 6.4
 \cite{Cu0}.
\begin{lemma} \label{th:main}   For any integer $2\le \ell \le  2N_1+1$
we have  transformations $\mathfrak{F} ^{(\ell )} = \mathfrak{F}_1 \circ \phi
_2\circ ...\circ \phi _\ell   $, with   $\mathfrak{F}_1$  the transformation in
Lemma \ref{lem:darboux} the $\phi _j$'s
  like in Lemma \ref{lem:chi}, such that
  the conclusions of  Lemma \ref{lem:ExpH11} hold
  but with $\textbf{R} _{-1}
^{(\ell )} =0$. We have the following expansion with $H_2^{(\ell )}$ of the form \eqref{eq:ExpH2}
\begin{equation*}
    K_0\circ  \mathfrak{F} ^{(\ell )}   = \psi(\pi, \Pi, \Pi(f)) + H_2^{(\ell )}   +
    \resto^{1,2}_{k,m}(\pi , \Pi  ,\Pi (f) , f)+ \sum _{j=0}^4
    \textbf{R}_{j}^{(\ell )} \text{ with:}
\end{equation*}   \begin{itemize}
\item[(i)]  all the nonzero terms in $\textbf{R} _0 ^{(\ell )} $ with $|\mu
    +\nu |\le \ell $  are in normal form, that is $ \mathbf{{e}}  \cdot (\mu -\nu
    )=0$;

\item[(ii)] all the nonzero terms in $\textbf{R} _1^{(\ell )}  $ with $|\mu
    +\nu |\le \ell -1$  are in normal form, that is $ \mathbf{{e}} \cdot (\mu -\nu
    )\in \sigma _e(\mathcal{H} _{\omega _0})$. \end{itemize} \end{lemma}

\proof (Sketch).  The proof consists of one part where the statement is
assumed for $\ell >1$ and then proved for $\ell+1$ and another part where the
statement is proved for $\ell =2$. Here we will briefly sketch the proof
of the first.

 We   assume that 
 $H^{(\ell)}:= K_0\circ  \mathfrak{F} ^{(\ell )}   $   for  an  $\ell >1$ has the desired properties for  indexes $(k',m')$ (instead of $(k,m)$) arbitrarily large.
We consider the representation \eqref{eq:ExpH11} for $H^{(\ell)}$ and we set $\mathbf{h}=H^{(\ell)}(\pi , \Pi , z,f,\varrho )$    replacing
  $\Pi (f)$ with $\varrho$ in \eqref{eq:ExpH11}.  We have equalities
\begin{align}
   \label{eq:derivmain1}  & g_{\mu \nu}(\pi , \Pi ,\varrho )
=\frac{1}{\mu !\nu ! }  \partial _z ^\mu \partial _{\overline{z}}^\nu
 \mathbf{h}_{|(z,f,\varrho )= (0,0,\varrho )}   \  ,  \quad |\mu +\nu |\le 2 {N}+1 ,\\&  \label{eq:derivmain2} \im J G_{\mu \nu}(\pi , \Pi ,\varrho )
=\frac{1}{\mu !\nu ! } \partial _z ^\mu \partial _{\overline{z}}^\nu \nabla _f
 \mathbf{h} _{|(z,f,\varrho )= (0,0,\varrho )} \  ,  \quad |\mu +\nu |\le 2 {N}
\end{align}
where we are using  that  $\mathbf{h} $ is
	$C^{2 {N}+2}$
	near $(p_0,p_0,0,0)$  in
$\Ph ^{s _0}_{3}=\{  ( \pi , \Pi ,\varrho ,  r ) \}$ for $m '\ge    2 {N}+2 $  for $s_0>   2 $
by
   Lemma \ref{lem:ExpH11}.   In particular, Lemma \ref{th:main} is a statement on
    the fact that some of the derivatives in \eqref{eq:derivmain1}--\eqref{eq:derivmain2} are 0.

   We consider now  a yet unknown $\chi$ as in  \eqref{eq:chi1}   with $M_0=\ell$,
 $i=0$,  $M=m'$ and $r=k'$.
Set $\phi :=\phi ^1$,   $\phi ^t$  the flow
of Lemma \ref{lem:chi}. We are seeking $\chi$  such that $H^{(\ell )}
\circ\phi $ satisfies the conclusions of Lemma \ref{th:main} for $\ell +1$.  Like in \cite{Cu0} we have \begin{equation}  \label{eq:backH}    \begin{aligned} &   H_2 ^{(\ell )}\circ
\phi  = H_2^{(\ell )} +  \int _0^1 \{  H_2^{(\ell )}  , \chi \}  ^{\pi cl}\circ \phi ^t dt  + \resto ^{0, \ell+2} _{\widetilde{k},m'} (\pi , \Pi , \Pi (f), r)
\end{aligned}
\end{equation}
with  $ \widetilde{{k}}=    k^{\prime   } - 2m'-1   $ and with $\{  H_2^{(\ell )}  , \chi \}  ^{\pi cl}=\{  H_2^{(\ell )}  , \chi \}  ^{(\pi ) }| _{\Pi (f) =\varrho} $ for  $\varrho$
treated like a constant. Furthermore, for  $\ell >1$ like in \cite{Cu0}  we have
\begin{equation} \label{eq:key001}\begin{aligned} &   H_2^{(\ell )}  \circ
\phi  = H_2^{(\ell )}  +    \{  H_2^{(\ell )}  , \chi \}  ^{\pi cl}+    \resto ^{0, \ell+2} _{\widetilde{k},m'} (\pi , \Pi ,\Pi (f), r).
\end{aligned}
\end{equation}
Set $\mathbf{e}  ^{(\ell)} _{j} (\pi , \Pi , \varrho ) =  g_{\delta _j, \delta _j } (\pi , \Pi  ,\varrho )   $ and $\mathbf{e}  ^{(\ell)}=(\mathbf{e}  ^{(\ell)} _{1},...,\mathbf{e}  ^{(\ell)} _{\mathbf{n}})$. Then  for $\varrho =\Pi (f)$
\begin{equation*}  \label{eq:key0}   \begin{aligned} &
 \{  H_2 ^{(\ell )}   , \chi \} ^{\pi cl}=     \im  \sum _{|\mu +\nu |=\ell +1}\mathbf{e}  ^{(\ell)}  ( \pi , \Pi ,\varrho  )\cdot (\mu -
\nu) z^{\mu} {\overline{z}}^{\nu} b_{\mu \nu }  ( \pi , \Pi ,\varrho  ) \\& +  \sum _{|\mu +\nu |=\ell  }  z^{\mu} {\overline{z}}^{\nu} \langle  (   \mathbf{e}  ^{(\ell)}  (   \pi , \Pi ,\varrho   )\cdot (\mu -
\nu)-  \mathcal{H} _{\pi})J^{-1}B _{\mu \nu }( \pi , \Pi ,\varrho   ),P_c(\pi )f
\rangle  +   \mathcal{{K}} \end{aligned}
\end{equation*}
with $ \mathcal{{K}} =   \mathcal{{K}} \big ( \pi , \Pi ,\varrho ,b( \pi , \Pi ,\varrho),B( \pi , \Pi ,\varrho)  \big)$
 a special case of a function like
\begin{equation}\label{eq:tildeKrho} \begin{aligned} & \mathcal{{K}} :=
 \sum _{|\mu +\nu |=\ell +1}  {k}_{\mu\nu} \big ( \pi , \Pi ,\varrho  ,b( \pi , \Pi ,\varrho),B( \pi , \Pi ,\varrho) \big ) z^{\mu} \overline{z}^{\nu} \\&+\sum _{ {i=0}}^{1}\sum _{ {j=1}}^{4}
 \sum _{|\mu +\nu |=\ell } z^{\mu}
\overline{z}^{\nu}
 \langle   J ^{-1}  \Diamond ^i_{j}  K_{ij \mu   \nu
}\big (\pi , \Pi , \varrho ,  b( \pi , \Pi ,\varrho),B( \pi , \Pi ,\varrho)   \big )
  , f \rangle ,\end{aligned}
\end{equation}
with $  \mathcal{{K}}( p_0 , p_0 ,0 ,b,B  )\equiv 0$
and dependent linearly on $(b ,B )$     with  $b_{\mu \nu}\in \C$, $B_{\mu \nu}\in \Sigma _{k'}$ and $k_{\mu \nu}\in \C$, $K_{ij\mu \nu}\in \Sigma _{k' }$. In fact our  $\mathcal{K}$ is
 \begin{align} &  \mathcal{{K}} =
   \sum _{|\mu +\nu |=\ell  }  z^{\mu} {\overline{z}}^{\nu} \langle  [1-(P _c(p_0 ) P_c (\pi )P _c (p_0 ) )^{-1} ] B _{\mu \nu }( \pi , \Pi ,\varrho   ),J^{-1}\mathcal{H} _{\pi} P_c(\pi )f
\rangle   \nonumber \\& + \sum _{|\mu +\nu |=\ell  }  z^{\mu} {\overline{z}}^{\nu}  \mathbf{e}  ^{(\ell)}  (   \pi , \Pi ,\varrho   )\cdot (\mu -
\nu)\langle  J^{-1}   B _{\mu \nu }( \pi , \Pi ,\varrho   ),  P_d(\pi )f
\rangle  +  \label{eq:tildeKrho1} \\&  \sum _{\substack{|\mu +\nu |=2\\ (\mu , \nu )\neq (\delta _j , \delta _j ) \ \forall \ j}  }  a_{\mu \nu}^{(\ell   )}( \pi , \Pi ,\varrho    )
\big [\sum _{|\mu '+\nu '|=\ell +1}
\{ z^\mu
{\overline{z}}^\nu ,  z^{\mu'}
{\overline{z}}^{ \nu '} \} ^{(\pi ) }  b_{\mu ' \nu '} ( \pi , \Pi ,\varrho  )  \nonumber\\& -\im
\sum _{|\mu '+\nu '|=\ell  }
\{  z^\mu
{\overline{z}}^\nu ,z ^{\mu'}
{\overline{z}}^{ \nu '} \} ^{(\pi ) } \langle  J^{-1}  B_{\mu ' \nu '}( \pi , \Pi ,\varrho  )
  , f \rangle
 \big ]  \ . \nonumber
\end{align}
By the argument in Theorem 6.4 \cite{Cu0}, the monomials of degree $\ell+1$ in the expansion
 of    $ H_2 ^{(\ell )}\circ
\phi $ of degree $\le 1$ in $f$ are of the form
\begin{equation}\label{eq:key1}
 \{  H_2^{(\ell )}  , \chi \}  ^{\pi cl} + \mathcal{G}+ \mathcal{{K}}
\end{equation}
with $ \mathcal{{K}}$ like in \eqref{eq:tildeKrho} and with
\begin{equation*}\label{eq:Krho}
\mathcal{G}(\pi , \Pi , \varrho ) :=\sum _{|\mu +\nu |=\ell +1} g_{\mu\nu} (\pi , \Pi , \varrho   ) z^{\mu} \overline{z}^{\nu} +\im \sum _{|\mu +\nu |=\ell} z^{\mu} \overline{z}^{\nu}
 \langle  J    G _{\mu   \nu
}(\pi , \Pi , \varrho  )
  , f \rangle
\end{equation*}
the monomials of degree $\ell+1$ in the expansion
 of    $ H  ^{(\ell )}$. Beyond \eqref{eq:tildeKrho1}, the $\mathcal{K}$
 in \eqref{eq:key1} contains contributions from the pullback of ${
\psi} (\pi, \Pi ,\Pi ( f))$, ${\textbf{R} _2} $  and $\resto ^{1,2}_{k,m}(\pi , \Pi ,\Pi (f), f)$. See \cite{Cu0}  for the details.

 By gathering in a polynomial $Z _{\ell +1 }$ the monomials in normal form in  \eqref{eq:key1}, we obtain that  \eqref{eq:key1}$=Z _{\ell +1 }$ by
 picking $(b,B)$ such that
\begin{align}&
\mathbf{e} ^{(\ell )}
\cdot (\mu - \nu) b_{\mu \nu}(\pi , \Pi , \varrho   ) + a_{\mu \nu}(\pi , \Pi , \varrho   ) + {k}_{\mu \nu} =0  \nonumber \\&  B_{\mu   \nu } (\pi , \Pi , \varrho   ) = \label{eq:NLhom11}\\& (P_c(\pi ) P_c(p_0 ) ) ^{-1}  R_{\mathcal{H}_{\pi }}( \mathbf{e}  ^{(\ell )} \cdot (\mu - \nu) ) (P_c(p_0) P_c(\pi ) ) ^{-1}\left [ G_{\mu \nu }(\pi , \Pi , \varrho   ) +   {{K}}  _{\mu \nu} \right ] \nonumber
\end{align}
with the first (resp. second) line for all pairs $(\mu ,  \nu)$  with $|\mu |+|  \nu|=\ell +1$ (resp. $|\mu |+|  \nu|=\ell  $) not in normal form. We observed already in Remark \ref{lem:normfor} that these indexes do not depend on $\pi$.
\eqref{eq:NLhom11} is a linear system in $(b (\pi , \Pi , \varrho   ),B(\pi , \Pi , \varrho   ))$
and  is a small perturbation of the case with $(k,K)=0$ when $(\pi , \Pi , \varrho   )$
is close to $(p_0 , p_0 , 0   )$.
The discussion for $\ell =1$ is similar, except now  the  $\mathcal{K}$ in \eqref{eq:key1}  is
like in \eqref{eq:tildeKrho} but with nonlinear dependence on the  $(b (\pi , \Pi , \varrho   ),B(\pi , \Pi , \varrho   ))$ and defined for $(b,B)$ small.
However, since $\mathcal{G}(\pi , \Pi , \varrho )=0$  for $(\pi , \Pi , \varrho )=(p_0 , p_0 , 0   )$,
\eqref{eq:NLhom11}   can   be solved by implicit function theorem.  For the details see \cite{Cu0}.

\qed

We will write  $H^{(\pi )}= K_0\circ  \mathfrak{F} ^{(\ell )} $ for
$\ell = 2N_1+1$  and we will set
\begin{equation}
 \label{eq:partham}H^{(\pi )}=\psi(\pi, \Pi, \Pi(f)) + H'_2+Z_0+Z_1+\resto
\end{equation}
with $H'_2$ like \eqref{eq:ExpH2}, $Z_0$   like \eqref{e.12c}, $Z_1$   like \eqref{e.12a}
and
\begin{equation}
 \label{eq:partham2}\resto =\sum _{j=2}^4\mathbf{R}_j +\resto ^{1,2}_{k,m}(\pi , \Pi ,\Pi (f), f).
\end{equation}

\section{Differentiations}
\label{sec:differentiations}

For a function $F(\pi , u)$ and  $\{\cdot,\cdot \}^{(\pi)}$
   the Poisson bracket w.r.t.\ $\Omega^{(\pi)}$ (see \eqref{eq:poiss})
\begin{equation} \label{eq:differentiations} \begin{aligned} &    \frac d{dt}   F(\pi ,
u) =   \partial _\pi   F(\pi , u) \dot \pi  +
  \{   F(\pi , u)  , \textbf{E}(u)  \} \\& =   \partial _\pi   F(\pi , u) \dot \pi  +
  \{   F(\pi , u (\pi , \tau ,\Pi , z,f) )  , \textbf{E}(u (\pi , \tau ,\Pi , z,f))  \}  ^{(\pi )} ,
\end{aligned}  \end{equation}

In the sequel  we  denote by $(\tau ' , \Pi ' , r')$ initial and  by $(\tau   , \Pi   , r )$
final coordinates.

\subsection{Differentiation in $\pi $}
\label{subsec:diffpi}

  The  coordinates depend on  the parameter $\pi$.  Specifically
we have functions $\tau ' (u)$,   $\Pi ' (u)$, and $r'(\pi , u)$. In particular, by Lemma \ref{lem:gradient R}
we have   $r'(\pi , u)\in C ^j ( D _{\R^4} (p_0, a) \times D _{H^l} (\Psi _{p_0}, a)   , H^{l -j}  )$  for
$j=0,1$ and $l=0,1$, with $a>0$ small and $D _{X} (u_0, a)$  the open disk of center $u_0$ and
radius $a$ in $X$.  The coordinate changes have not affected the variable $\Pi$ so that
$\Pi (u)=\Pi ' (u)$.

 In the remainder of the paper we will drop  indexes  $(k,m)$ and  we will write
\begin{equation} \label{eq:newsym} \begin{aligned} & \resto ^{i,j}= \resto ^{i,j}_{k,m} \text{ and }  \textbf{S} ^{i,j}=\textbf{S} ^{i,j}_{k,m} \text{ assuming $k$ and $m$   large at will,}
   \end{aligned}
\end{equation}
  abusing notation. This can be done. Indeed
   while  pull  back  and other operations
decrease  the values of   $(a,b)$ of
      $   \resto ^{i,j}_{a,b}$ and  $  \textbf{S} ^{i,j}_{a,b}$, by taking the initial  $(a,b)$'s large enough and by the control on their decrease by   Lemmas \ref{lem:ODEdomains}--\ref{lem:transf} , we can always assume that after all the various operations, the symbols are of the type in \eqref{eq:newsym}.

\begin{lemma}\label{lem:reg1} The final functions   $(\tau      , r )$
have the following dependence   on  $(\pi , u)$.  Set  $\varrho  =\Pi (r )$.  Then $\exists$   $a>0$
such that for $D _{\R^4} (p_0, a)$  (resp.   $D _{H^1} (\Phi _{p_0}, a)$)  the open disk  in $\R^4$ (resp.
  $ H^1 $)   of center $p_0$   (resp.
  $ \Phi _{p_0} $)   and  radius $a$, and for    ${\mathcal V}=\cup _{\tau \in \R^4}e^{J\tau \cdot \Diamond} D _{H^1} (\Phi _{p_0}, a)=\cup _{\tau \in \R^4} D _{H^1} (e^{J\tau \cdot \Diamond}\Phi _{p_0}, a) $, we have:

\begin{itemize}
\item[(1)]     $\tau (\pi , u), \varrho (\pi , u) \in C ^1 (   D _{\R^4} (p_0, a) \times {\mathcal V}   , \R ^4  )$;
\item[(2)]     $r (\pi , u)
\in C ^i (   D _{\R^4} (p_0, a) \times {\mathcal V} , H^{1-i})$  for $i=0,1$.
\item[(3)]     $\partial _\pi \tau (\pi , u)_{|u=u(\pi , \tau   , \Pi   , r )}=\resto ^{0,2} $ and $\partial _\pi\varrho (\pi , u)_{|u=u(\pi , \tau   , \Pi   , r )}=\resto ^{0,2}$;
\item[(4)]       $\partial _\pi  r (\pi , u)_{|u=u(\pi , \tau   , \Pi   , r )}=  J\resto ^{0,2}\cdot \Diamond r+\textbf{S} ^{0,1}   $.

\end{itemize}
  \end{lemma}
\proof   Set $\varrho '=\Pi (r')$.  We know that

\begin{equation}\label{eq:reg2} \begin{aligned}  &
 \tau = \tau ' +\resto ^{0,2} (\pi , \Pi' , \varrho  ' , r') \quad , \quad \Pi  = \Pi '  \\&
\varrho =\varrho ' +  \resto ^{0,2} (\pi , \Pi ', \varrho  ' , r')  \quad ,\\&
 r= e^{J\resto ^{0,2} (\pi , \Pi ', \varrho '  , r')\cdot \Diamond } r' + \textbf{S} ^{1,1}(\pi , \Pi ', \varrho '  , r')).
\end{aligned}\end{equation}
Notice that  varying $\Pi '$ near $p_0$ in $\R ^4$ and $r'$ near 0 in $N^\perp _g({\mathcal H}_{p_0}) \cap H^1$ and for $\tau '\in \R ^3\times S^1$,  if we define $p'=p'(\pi ,\Pi , r ) $  using \eqref{eq:variables},  then
$u=u(p', \tau ',r') $, defined using \eqref{eq:coordinate}, spans a set  like $ {\mathcal V}$.

Substituting  $\tau ' =\tau ' (u)$,   $\Pi ' =\Pi ' (u)$,   $\varrho   ' =\varrho   ' (\pi ,u)$,  and $r'=r'(\pi , u)$ in
 \eqref{eq:reg2},
we obtain   functionals  $\tau  =\tau  (\pi ,u )$,   $\Pi   =\Pi ' (u)$,   $\varrho   =\varrho    (\pi ,u)$,  and $r =r (\pi , u)$.  By Lemma \ref{lem:gradient R}  we can apply the
  Chain Rule, which yields   Claims (1) and (2).

By   \eqref{eq:gradient R2}  we   have
$ \partial _\pi  r '(u,\pi ) \sim   (\partial _\pi P(\pi ))   r'$    and  so $  \partial _\pi  r '(u,\pi ) =\textbf{S} ^{0,1}$
by \eqref{eq:tranR}.  This yields also   $ \partial _\pi  \varrho  '(u,\pi ) = \resto ^{0,2} $.
By   the Chain Rule   we get

\begin{equation}\label{eq:reg3} \begin{aligned}  &
  \partial _\pi [\resto ^{0,2} (\pi , \Pi' , \varrho  ' , r')] =  \underbrace{\partial _\pi \resto ^{0,2}  }_{\resto ^{0,2} }+  \underbrace{\partial _{\varrho  '}  \resto ^{0,2} \cdot   \partial _\pi\varrho  '}_{\resto ^{0,2} \cdot \resto ^{0,2}} +  \partial _{r  '}  \resto ^{0,2} \cdot   \partial _\pi r  '  =\resto ^{0,2},
\end{aligned}\end{equation}
since    $\partial _{r  '}  \resto ^{0,2}$ has values in $B (\Sigma _{-k},
\R )$  and $\|   \partial _{r  '}  \resto ^{0,2} \|  _{B (\Sigma _{-k},
\R )}\le C \|   r \|  _{\Sigma _{-k}}$, for some large $k$. This yields Claim (3). Similarly

\begin{equation}\label{eq:reg4} \begin{aligned}  &
  \partial _\pi [\textbf{S} ^{1,1} (\pi , \Pi' , \varrho  ' , r')] =  \underbrace{\partial _\pi \textbf{S} ^{1,1}  }_{\textbf{S} ^{0,1} }+  \underbrace{\partial _{\varrho  '}  \textbf{S} ^{1,1} \cdot   \partial _\pi\varrho  '}_{\textbf{S} ^{0,1} \cdot \resto ^{0,2}} +  \partial _{r  '}  \textbf{S} ^{1,1} \cdot   \partial _\pi r  '  =\textbf{S} ^{0,1},
\end{aligned}\end{equation}
because $\partial _{r  '}   \textbf{S} ^{1,1} $ has values in $B (\Sigma _{-k},
\Sigma _{k} )$  for a large $k$  and    $\|   \partial _{r  '}  \textbf{S} ^{1,1}  \|  _{B (\Sigma _{-k},
\Sigma _{k}  )}\le C(  |\pi - \Pi | + |\varrho | + \|   r \|  _{\Sigma _{-k}})$. This yields Claim (4).

\qed

Rather than   $r\in N^\perp _g (\mathcal{H}^*_{p_0})$ we use  $(z, f)\in \C ^{\textbf{n}}\times L^2_c(p_0)$, related to $r$ by  \eqref{eq:decomp2}.

\begin{lemma}\label{lem:reg5}
The   functions   $(z      , f )$
have the following dependence   on  $(\pi , u)$.

\begin{itemize}
\item[(1)]   We have $z (\pi , u) \in C ^1 (   D _{\R^4} (p_0, a) \times D _{H^1} (\Psi _{p_0}, a)   , \R ^4  )$.
\item[(2)]   We have $f (\pi , u)
\in C ^i (   D _{\R^4} (p_0, a) \times D _{H^1} (\Psi _{p_0}, a)   , H^{1-i})$  for $i=0,1$.
\item[(3)]    We have $\partial _\pi z (\pi , u) _{|u=u(\pi , \tau   , \Pi   , r )}   =\resto ^{0,1} $.
\item[(4)]    Consider the function  $\partial _\pi  f (\pi , u)$
after expressing $u$  in terms of  $(\pi , \tau   , \Pi   , r )$. Then       $\partial _\pi  f  =J\resto ^{0,2}\cdot \Diamond f+\textbf{S} ^{0,1}   $.

\end{itemize}
\end{lemma}
\proof   (1) follows from $z_j=\Omega ( P (\pi )r, \overline{\xi} _j ^{(\pi )} ) $,
Lemmas \ref{lem:basis} and \ref{lem:reg1}, \eqref{eq:projreg}--\eqref{eq:projP} and
\eqref{eq:decomp2}. Then $P_c(\pi ) f$ satisfies the conclusions of (2).
 By \eqref{eq:Pd}, Lemma \ref{lem:basis} and by \eqref{eq:projreg}  $P_{d }(p)\in C^{\infty}
(\mathcal{P}, B (\mathcal{S}',\mathcal{S}))$. Then by $f= (1+P_{d }(p_0)-P_{d }(\pi) )^{-1}P_{c }(\pi)f$ we get (2).
Differentiating    \eqref{eq:decomp2} in $\pi $ we have
\begin{equation*}    \begin{aligned} &
 (\partial _\pi P (\pi ))   r + P (\pi ) \partial _\pi r=  (\partial _\pi z ) \xi    ^{(\pi )} + z  \partial _\pi  \xi    ^{(\pi )}
+\text{c.c.}+ P_c(\pi ) \partial _\pi  f  +  (\partial _\pi P _c(\pi )) f
 \end{aligned}
\end{equation*}
with c.c.   the complex conjugate of the first two terms in the rhs. So, by Claim (4) in Lemma
\ref{lem:reg1},       we have the following, which yields   Claims (3)--(4)
\begin{equation*}    \begin{aligned}     \xi    ^{(\pi )}\partial _\pi z +    \overline{\xi}    ^{(\pi )}  \partial _\pi \overline{z} +   P_c(\pi ) \partial _\pi  f  &=   P (\pi )  J\resto ^{0,2}\cdot \Diamond r  +\textbf{S} ^{0,1} =     J\resto ^{0,2}\cdot \Diamond f   +\textbf{S} ^{0,1}.
 \end{aligned}
\end{equation*}

\subsection{Differentiation in $( \Pi , z,f)  $}
\label{subsec:diffothers}

\begin{lemma}\label{lem:diffothers1}  Consider the last system of coordinates.
Then we have
\begin{equation}  \label{eq:diffothers2} \begin{aligned} &
  \partial _{z_j}r= (P(\pi ) P(p_0) )^{-1}\xi _j^{(\pi )}\, , \, \partial _{\overline{z}_j}r= (P(\pi ) P(p_0) )^{-1}\overline{\xi} _j^{(\pi )}\, , \\& \partial _{f}r= (P(\pi ) P(p_0) )^{-1} P_c(\pi ) P_c(p_0)
   \end{aligned}
\end{equation}
 For $p=\Pi -\Pi (P(p) r)$, we have
 \begin{equation}  \label{eq:diffothers3} \begin{aligned} &
  \partial _{z_j}p= \resto ^{0,1}\, , \, \partial _{\overline{z}_j}p= \resto ^{0,1}\, ,    \, \partial _{\Pi _k}p_j=\delta _{jk} +\resto ^{0,2}
  \, , \ \nabla _{f}p_j=   \mathbf{S}^{0,1}- \Diamond _jf.
   \end{aligned}
\end{equation}
 \end{lemma}
\proof  \eqref{eq:diffothers2} follows immediately from  $P(\pi ) r=z\cdot \xi  ^{(\pi )} +\overline{z}\cdot \overline{\xi}  ^{(\pi )} + P_c(\pi ) f    $, that is from \eqref{eq:decomp2}.  \eqref{eq:diffothers3} follows easily from $p=\Pi -\Pi (P(p) r)$.
\qed

\begin{lemma}\label{lem:diffothers4}  Consider the first system of coordinates $(\tau ', \Pi , r')  $.
Then
\begin{equation}  \label{eq:diffothers5} \begin{aligned} &
  \nabla  _{z,  \overline{z}}\tau ' = \resto ^{0,1}   \, , \, \partial _{\Pi }\tau '=  \resto ^{0,2}\, , \ \nabla  _{f}\tau '= \mathbf{S}^{0,1} .
   \end{aligned}
\end{equation}
   $p'=\Pi -\Pi (P(p') r')$ satisfies analogous formulas as $p$ in
 \eqref{eq:diffothers3}. Finally, we have
 \begin{align}    &
  \partial    _{z,  \overline{z}}r'=  e^{J\resto ^{0,2}\cdot \Diamond }
  \left ( \partial    _{z,  \overline{z}}r  +J\partial    _{z,  \overline{z}}\resto ^{0,2}\cdot \Diamond    r \right )  + \mathbf{S}^{1,0}    \nonumber \\&
   \partial    _{\Pi}r'=   e^{J\resto ^{0,2}\cdot \Diamond }J\partial    _{\Pi}\resto ^{0,2}\cdot   \Diamond    r + \mathbf{S}^{1,1}
  \label{eq:diffothers6}
  \\& \partial _{f}r'=  e^{J\resto ^{0,2}\cdot \Diamond }
  P_c(p_0) + \gimel +   e^{J\resto ^{0,2}\cdot \Diamond }
   J\partial _{f}\resto ^{0,2}\cdot \Diamond   r
  +\mathbf{S}^{1,1}  \partial _{f}\resto ^{0,2}
   ,\nonumber \end{align}
with $\gimel \in C ^m (  \U  , B (\Sigma _{-k},\Sigma _{ k})) $,
with $\U $ a neighborhood of $(p_0, p_0, 0, 0)$ in $\mathcal{P}^{-k}_3=(\pi ,
\Pi , \varrho , r)$, where $m$ and  $k$ are fixed large numbers.
  \end{lemma}
\proof  \eqref{eq:diffothers5} follows from $\tau ' = \tau +\resto ^{0,2}.$  Notice that this follows from
Lemma \ref{lem:ODE} and from the fact that  the $\resto ^{0,2}$ are preserved by changes of variables
as of Lemma \ref{lem:ODEdomains}.
  The claim for $p'$ follows from  $p ' = p +\resto ^{1,2}.$

 \noindent Differentiating $r'=  e^{J\resto ^{0,2}\cdot \Diamond }  ( r
  +\mathbf{S}^{1,1})$ we obtain  \eqref{eq:diffothers6}.
  In particular we have
\begin{equation}  \label{eq:parr'}    \begin{aligned} & \partial _{f}r'=  e^{J\resto ^{0,2}\cdot \Diamond }
  \left ( \partial _{f}r+\partial _{f}\mathbf{S}^{1,1}  + J\partial _{f}\resto ^{0,2}\cdot \Diamond  ( r
  +\mathbf{S}^{1,1})  \right ) .
  \end{aligned}
\end{equation}
Notice that we have $ e^{J\resto ^{0,2}\cdot \Diamond }J \Diamond  \mathbf{S}^{1,1}  =  \mathbf{S}^{1,1} .$
  By  \eqref{eq:diffothers2} we have
  \begin{equation}    \begin{aligned} & \partial _{f}r = (P(\pi ) P(p_0) )^{-1} P (\pi ) P_c(p_0)  -(P(\pi ) P(p_0) )^{-1} (P_d (\pi ) -P _{N_g}(\pi ))  P_c(p_0)\\& =P_c(p_0) -(P(\pi ) P(p_0) )^{-1} (P_d (\pi ) -P _{N_g}(\pi ))  P_c(p_0)  =P_c(p_0)+ \gimel_1 .
  \end{aligned} \nonumber
\end{equation}
We can set $\gimel :=\gimel _1+e^{J\resto ^{0,2}\cdot \Diamond }
   \partial _{f}\mathbf{S}^{1,1}$.
	
	The first two formulas in \eqref{eq:diffothers6} are obtained proceeding similarly.
\qed

We have the following result after setting
\begin{equation}  \label{eq:V2} \begin{aligned}& V_2(t)=\gamma V (\cdot + \textbf{v} t+ D') .
\end{aligned}
\end{equation}

\begin{lemma}\label{lem:potaverage} For a   $\gimel \in C ^m (  \U  , B (\Sigma _{-k},\Sigma _{ k})) $,
  $\U\subset \mathcal{P}^{-k}_3 $  as in Lemma \ref{lem:diffothers4}:
\begin{align}& \nonumber
\nabla _f  \langle \gamma V (\cdot +  \mathbf{{v}} t )  u
,  u  \rangle  =
  2\langle  \partial _{D'}V_2 (t) (\Phi _{p'}+P(p')P(\pi )r'),
\Phi _{p'}+P(p')P(\pi )r'\rangle   \mathbf{S}^{0,1}  \\& \label{eq:potaverage1} -
 2\langle   V_2 (t) (\Phi _{p'}+P(p')P(\pi )r'), \mathbf{S}^{0,0}   \rangle   \Diamond f
\\&   \nonumber+ 2\langle  V_2 (t) (\Phi _{p'}+P(p')P(\pi )r') , P(p')P(\pi )
e^{J\resto ^{0,2}\cdot \Diamond} J  \Diamond r \rangle  \mathbf{S}^{0,1}  \\& \nonumber +2  ( P_c^*(p_0)e^{-J\resto ^{0,2}\cdot \Diamond} +\gimel )  V_2 (t) (\Phi _{p'}+P(p')P(\pi )r')
\end{align}
with the $e^{-J\resto ^{0,2}\cdot \Diamond} $ in the last line   the inverse of  $e^{ J\resto ^{0,2}\cdot \Diamond} $
in  \eqref{eq:parr'};
\begin{align}  \nonumber &2^{-1}
\nabla _{z,\overline{z}}  \langle   \gamma V (\cdot +  \mathbf{{v}} t ) u
,  u  \rangle  =   \langle  V_2 (t) (\Phi _{p'}+P(p')P(\pi )r') , e^{J\resto ^{0,2}\cdot \Diamond} J  \Diamond r \rangle   \resto ^{0,1}
\\&  \label{eq:potaverage2}+
  \langle   V_2 (t)(\Phi _{p'}+P(p')P(\pi )r'), \mathbf{S}^{0,0}  \rangle   \\&  +
  \langle  \partial _{D'}V_2 (t) (\Phi _{p'}+P(p')P(\pi )r'),
\Phi _{p'}+P(p')P(\pi )r'\rangle  \resto^{0,1}  \nonumber .
\end{align}

\end{lemma}
\proof   Recall that
\begin{equation} \label{eq:pot} \begin{aligned}       \langle \gamma V u
,  u  \rangle  = \langle \gamma V (\cdot + \textbf{v} t+ D') (\Phi _{p'}+P(p')P(\pi )r') , \Phi _{p'}+P(p')P(\pi )r'  \rangle  . \end{aligned}  \end{equation}
By the Chain Rule  and for  ${\textbf{g}}$ equal either to $f$ or  to $(z,\overline{z})$, we get
\begin{align}&\nonumber  2^{-1}
\partial _{\textbf{g}}  \langle   \gamma V (\cdot + \textbf{v} t )  u
,  u  \rangle  = \langle  V_2 (t) (\Phi _{p'}+P(p')P(\pi )r'),
    P(p')P(\pi )\partial _{\textbf{g}} r'\rangle  \\& \label{eq:potaverage3}  =   2^{-1}
  \langle  \partial _{D'}V_2 (t) (\Phi _{p'}+P(p')P(\pi )r'),
\Phi _{p'}+P(p')P(\pi )r'\rangle  \partial _{\textbf{g}} D' \\& \nonumber +
  \langle  V_2 (t) (\Phi _{p'}+P(p')P(\pi )r'),
 \partial _{p'}\Phi _{p'}+ \partial _{p'}P(p')P(\pi )r'\rangle \partial _{\textbf{g}} p'
.
\end{align}
For  $\textbf{g}=f$ and by    \eqref{eq:diffothers6}, \eqref{eq:potaverage3} becomes
 \begin{align}& \nonumber   2^{-1}
  \langle  \partial _{D'}V_2 (t) (\Phi _{p'}+P(p')P(\pi )r')
\Phi _{p'}+P(p')P(\pi )r'\rangle  \nabla _f D' \\&  \nonumber  +
  \langle  V_2 (t)(\Phi _{p'}+P(p')P(\pi )r'),
 \partial _{p'}\Phi _{p'}+ \partial _{p'}P(p')P(\pi )r'\rangle  \nabla _f p'
\\& + \label{eq:potaverage4}
 2^{-1}
  \langle  V_2 (t)  (\Phi _{p'}+P(p')P(\pi )r') , \mathbf{S}^{1,1} \rangle \nabla _f \resto ^{0,2} \\&\nonumber  + \langle  V_2 (t) (\Phi _{p'}+P(p')P(\pi )r') , e^{J\resto ^{0,2}\cdot \Diamond} J  \Diamond r \rangle \nabla _f \resto ^{0,2} \\& \nonumber + (P_c^*(p_0)e^{-J\resto ^{0,2}\cdot \Diamond} +\gimel )
		  (\Phi _{p'}+P(p')P(\pi )r') .
\end{align}
Then we use : $\nabla _f D'$ and $\nabla _f \resto ^{0,2} $
 are  $ \mathbf{S}^{0,1} $; we have  $\nabla _f p'_j= -\Diamond  _jf + \mathbf{S}^{0,1} $; we have $ \partial _{p'_j}\Phi _{p'}+ \partial _{p'_j}P(p')P(\pi )r'= \mathbf{S}^{0,0}  $. Then    \eqref{eq:potaverage4}
    yields   the r.h.s of \eqref{eq:potaverage1}.

To get  \eqref{eq:potaverage2} we set $\textbf{g}=(z,\overline{z})$.  When we use that
$\nabla _{z,\overline{z}} D'$  and    $\nabla _{z,\overline{z}} p' $  are $\resto^{0,1} $, we get terms
which can be absorbed by the first two lines of  \eqref{eq:potaverage2}.
Using  \eqref{eq:diffothers6}, the fact that $\partial  _{z,\overline{z}} r= \mathbf{S}^{0,0}   $
and that    $\nabla _{z,\overline{z}}\resto^{0,2} =\resto^{0,1} $,
we then set
\begin{equation}
\label{eq:potaverage5} \begin{aligned}&  \partial  _{z,\overline{z}} r' = e^{J\resto ^{0,2}\cdot \Diamond }
   J \resto ^{0,1}\cdot \Diamond    r   + \mathbf{S}^{0,0}    .
\end{aligned}
\end{equation}
We get  Lemma  \ref{lem:potaverage}
by the following, which goes in the last two lines of   \eqref{eq:potaverage2},
 \begin{equation*}
\label{eq:potaverage6} \begin{aligned}&   \langle   V_2 (t) (\Phi _{p'}+P(p')P(\pi )r'),
    P(p')P(\pi )\partial  _{z,\overline{z}} r'\rangle  =   \langle   V_2 (t) (\Phi _{p'}+P(p')P(\pi )r'),
     \mathbf{S}^{0,0} \rangle   \\& + \langle  V_2 (t) (\Phi _{p'}+P(p')P(\pi )r'),
    P(p')P(\pi ) e^{J\resto ^{0,2}\cdot \Diamond }
   J \Diamond    r \rangle   \resto ^{0,1}      .\qed
\end{aligned}
\end{equation*}

\subsection{Formulation of the system}
\label{subsec:system}

From now we choose $\pi =\pi (t) =\Pi (t)$, as anticipated after \eqref{eq:conjNLS1}.

\noindent Recall that the variables  $\Pi _j$ are unchanged by the changes of coordinates,  see \eqref{eq:ODE}.  Then,
  for $a\le 3$ and $V_2(t)=\gamma V (\cdot + \textbf{v} t+ D')$, using \eqref{eq:differentiations} with $\Pi _a=F(u)$,
\begin{equation} \label{eq:mom} \begin{aligned}     & \dot   \Pi _a    =
\{  \Pi _a   , \textbf{E}( u)  \} ^{(\pi )} = -
  {\partial  _{D_a}} H^{(\pi )}  -  2^{-1}  \partial _{ D_a}
\langle  V_2(t) u , u  \rangle =  \\&     - 2^{-1}   \langle     \partial _{ D_a'} V_2(t) (\Phi _{p'}+ P(p')P(\pi )r') , \Phi _{p'}+ P(p')P(\pi )r'  \rangle
 {\partial _{ D_a}  D_a '}   . \end{aligned}  \end{equation}
Since the r.h.s. of \eqref{eq:tranR} is independent of $\tau$,  we   have  $ \partial  _{  D_a  } D_a ' =1$.  By Lemma  \ref{lem:ExpH11}  we have ${\partial  _{D_a}} H^{(\pi )}=0.$
  $Q$ is an invariant of motion.   Summing up, we have
    \begin{equation*}   \begin{aligned}   &   \dot  Q  =
0 \ ,   \  \dot   \Pi _a =- 2^{-1}   \langle     \partial _{ x_a } V_2(t) (\Phi _{p'}+ P(p')P(\pi )r') , \Phi _{p'}+ P(p')P(\pi )r'  \rangle
  . \end{aligned}  \end{equation*}
	By   $\tau '=\tau' (u)$ and by       \eqref{eq:differentiations} we have $\dot \tau _j  ' =\{   \tau _j  '  , \textbf{E}   \}$.
	So in particular, by \eqref{eq:differentiations} with $ D_a  '=F(u)$ and by \eqref{eq:PoissIdent2}, we have
	\begin{equation*} \begin{aligned} &    \dot D_a  '  - v'_a =
	\{ D_a  '    ,   K_0\} + 2^{-1}{\gamma }    \{ D_a  '      ,\langle  V (\cdot + \textbf{v} t) u, u\rangle\}    .\end{aligned}  \end{equation*}
	In terms of the last   coordinate system, where  $D ' = D +\resto ^{0,2} $    by \eqref{eq:tranR},
 by \eqref{eq:poiss}
   \begin{align}    \nonumber   \dot D_a  '  - v'_a &=
    {\partial   _{\Pi _a} } H^{(\pi )}
  +2^{-1}
    \partial _{\Pi _a}\langle   V_2(t) (\Phi _{p'}+ P(p')P(\pi )r') ,   \Phi _{p'}+ P(p')P(\pi )r'  \rangle \\& +\{ \resto ^{0,2}      ,   K_0\} + 2^{-1}{\gamma }   \{ \resto ^{0,2}        ,\langle  V (\cdot + \textbf{v} t) u, u\rangle\}    .\label{eq:eqD}\end{align}
		Notice that in the current notation,  \eqref{eq:weakint}  reads as    \begin{equation}
\label{eq:eqD0} \begin{aligned} &
 \gamma  \sup _{\text{dist} _{S^2}(  \overrightarrow{{e}}, \frac{\mathbf{v}}{|\mathbf{v}|} ) \le \delta _0}
 \|    \langle   |\mathbf{v}|\overrightarrow{{e}} t+D'(0)  \rangle  ^{-2}  \| _{L^1_t(\R _+)}  <10\epsilon .
\end{aligned}
\end{equation}

	\noindent 	Since $\vartheta  ' = \vartheta  +\resto ^{0,2} $  by \eqref{eq:tranR}, by  the argument used for \eqref{eq:eqD}  and   \eqref{eq:thetad}--\eqref{eq:PoissIdent2}
  \begin{equation}\label{eq:eqtheta}  \begin{aligned} &   \dot \vartheta  '  -  \omega '- 2^{-1}{(v')^2}  =
	\{ \resto ^{0,2}      ,   K_0\} + 2^{-1}{\gamma }    \{ \resto ^{0,2}        ,\langle  V (\cdot + \textbf{v} t) u, u\rangle\}
	\\&
	-
    {\partial  _{Q}} H^{(\pi )}
  -
  2^{-1}  \partial _{Q} \langle     V_2(t) (\Phi _{p'}+ P(p')P(\pi )r') ,   \Phi _{p'}+ P(p')P(\pi )r'  \rangle   .
\end{aligned}  \end{equation}
We have
\begin{equation} \label{eq:eqz} \begin{aligned} &      \dot z_j =-\im  \partial _{\overline{z}_j}H^{(\pi )} +     \dot \Pi _b \partial _{\pi _b} z_j
  -\im
   2^{-1}  \partial _{\overline{z}_j} \langle
			V(\cdot + \textbf{v}t ) u , u  \rangle    .\end{aligned}  \end{equation}
And finally we have
\begin{equation} \label{eq:eqf} \begin{aligned} &       \dot f =    \dot \Pi _b \partial _{\pi _b} f  + (P _c(p_0 ) P_c (\pi )P _c (p_0 ) )^{-1}   J \nabla _{f}K',
 \\&  \nabla _{f}K'=\nabla _{f}   H^{(\pi )}
  +
  2^{-1} \nabla _{f}  \langle      V(\cdot + \textbf{v}t ) u , u  \rangle
		  .\end{aligned}  \end{equation}
In \eqref{eq:eqz}--\eqref{eq:eqf} we can substitute the formulas from Lemma \ref{lem:potaverage}.

\begin{theorem}\label{thm:mainbounds} Consider the constants   $0<\epsilon <\epsilon _0  $  of Theorem \ref{theorem-1.1}. Then there   is a fixed
$C >0$ such that    for $I= [0,\infty )$    we have:
\begin{align}
&   \|  f \| _{L^p_t(I,W^{ 1 ,q}_x)}\le
  C \epsilon \text{ for all admissible pairs $(p,q)$,}
  \label{Strichartzradiation}
\\& \| z ^\mu \| _{L^2_t(I)}\le
  C \epsilon \text{ for all multi indexes $\mu$
  with  $\textbf{e}\cdot \mu >\omega _0 $,} \label{L^2discrete}\\& \| z _j  \|
  _{W ^{1,\infty} _t  (I )}\le
  C \epsilon \text{ for all   $j\in \{ 1, \dots , \textbf{n}\}$ }
  \label{L^inftydiscrete} \\&  \|  \omega  ' -  \omega _0  \| _{L_t^\infty  (I )}\le
  C \epsilon  \,  , \quad    \| v  ' -  v _0\| _{L^\infty_t (I ) }\le
  C \epsilon
  \label{eq:orbstab1} .
\end{align}
Furthermore,  there exist $\omega _+$ and  $v _+$ such that
\begin{align}
 &    \lim _{t\to +\infty} \omega  '(t)= \omega _+   \,  , \quad   \lim _{t\to +\infty}  v  ' (t)=  v _+
  \label{eq:asstab1}\\&  \label{eq:asstab2}  \lim _{t\to +\infty}  \dot D  ' (t)=  v _+ \, , \quad   \lim _{t\to +\infty}  \dot \vartheta   ' (t)=  \omega _++  4^{-1}{v^2_+} \\&  \label{eq:asstab3}  \lim _{t\to +\infty}  z(t)=0 .
\end{align}
\end{theorem}
Theorem  \ref{thm:mainbounds} will be obtained as a consequence of the following Proposition.

\begin{proposition}\label{prop:mainbounds} Consider the constants   $0<\epsilon <\epsilon _0  $  of Theorem \ref{theorem-1.1}.     There exist  a  constant $c_0>0$  such that
for any  $C_0>c_0$ there is an  $ \epsilon _0  >0 $ such that   if   the inequalities  \eqref{Strichartzradiation}--\eqref{L^inftydiscrete}
hold  for $I=[0,T]$ for some $T>0$
and for $C=C_0$, and if furthermore  for $t\in I$
\begin{align} \label{eq:unds1}  &
	 \| \dot D ' -v'  \| _{L^1 (0,t )}<  C \epsilon    \langle t \rangle \, , \\&  \label{eq:unds2}   \|  p'  -p_0 \| _{L^\infty  (I )}<   C \epsilon,
\end{align}
then in fact for $I=[0,T]$  the inequalities  \eqref{Strichartzradiation}--\eqref{L^inftydiscrete} hold  for   $C=C_0/2$  and the  inequalities  \eqref{eq:unds1}--\eqref{eq:unds2}
 hold  for   $C=c$  with $c$ a specific fixed constant.
\end{proposition}

\section{Proof that Prop. \ref{prop:mainbounds} implies Theor. \ref{thm:mainbounds}}
\label{sec:propthm}

By a continuation argument
Prop. \ref{prop:mainbounds}  gives  \eqref{Strichartzradiation}--\eqref{L^inftydiscrete} and \eqref{eq:unds2} 
    on $[0,\infty )$.    So   to complete to proof
of  Theor. \ref{thm:mainbounds},
it   remains   to prove formulas
\eqref{eq:orbstab1}--\eqref{eq:asstab3}.

\begin{lemma}\label{lem:proptheor1}
     Prop. \ref{prop:mainbounds}   implies     \eqref{eq:asstab3}, that is $\displaystyle \lim _{t\nearrow\infty} z  (t)=0$.
\end{lemma}
\proof   The limit \eqref{eq:asstab3} follows  from \eqref{L^2discrete}--\eqref{L^inftydiscrete} and  by $\dot z \in L^\infty [0,\infty )$, which follows
by  equation  \eqref{eq:eqz}.  Indeed for
$  (M-1)\textbf{e} _j >\omega _0$ the following limit exists,
\begin{equation*}
 \lim _{t\nearrow\infty} z ^{2M}_j(t)=z^{2M}_j(0)+ 2M\int _0^\infty \dot z_j (t) z ^{2M-1}_j(t) dt,
\end{equation*} since $\dot z_j    z ^{2M-1}_j\in  L^1 [0,\infty )$ with
$\| \dot z_j  z ^{2M-1}_j\| _{L^1} \le \| \dot z_j   \| _{L^\infty}\|  z_j \| _{L^\infty}\|   z ^{ M-1}_j\| _{L^2}^2.$  Since  $z ^{2M}_j(t)\in L^1 [0,\infty )$ by  \eqref{L^2discrete}, the above limit  is necessarily equal to 0.
\qed

\begin{lemma}\label{lem:scattf}   There is  a fixed $C$ and  $f_+  \in H^1$  and a function $\varsigma :[0, \infty   )\to \R ^4$
such that for the variable $f$ in
\eqref{Strichartzradiation} we have
\begin{equation}\label{eq:scattering1}
 \begin{aligned} &    \lim _{t\nearrow \infty } \|   f(t) - e^{J \varsigma (t) \cdot \Diamond }
  e^{ -J  t   \Delta    } f_+ \| _{H^1 }=0 .
\end{aligned}
\end{equation}

\end{lemma}

\proof  We     defer   the proof  until Sect. \ref{sec:scattf}.    \qed

 We will consider now  a number of technical lemmas.

\begin{lemma}\label{lem:movpot}    Consider the hypotheses of   Proposition \ref{prop:mainbounds},
  a fixed $S^{0,0}_{2k,0}$ where $k>3$  and a fixed $q \in \mathcal{S}(\R ^3)$. Then   for $\epsilon _0$ small enough
there exists a fixed constant   $c$ dependent on   $S^{0,0}_{2k,0}$  and  $q $ such that

\begin{equation} \label{eq:movpot1} \begin{aligned} &
   \|    \gamma q (\cdot +  \mathbf{{v}} t+ D')S^{0,0}_{2k,0} \| _{L^1((0,T),L^p_x)}  \le c \epsilon  \text{  for all $p\ge 1$}.
		  	 \end{aligned}  \end{equation}
\end{lemma}
\proof    We have  by   $k>3$ and Sobolev embedding,
\begin{equation*} \label{eq:potinter1} \begin{aligned} &    \|
  q(\cdot + \textbf{v} t+ D')     S^{0,0}_{2k,0}  \| _{L^p_x}    \le   C_{q,k}    \|      S^{0,0}_{2k,0}  \| _{\Sigma  _{2k}} \langle D '(t)+t\textbf{v}\rangle ^{-k  }    .
	 	 \end{aligned}  \end{equation*}
Then   for a fixed $C= C_{q,k,S}   $
\begin{equation*} \label{eq:inteqD} \begin{aligned} &   \|    \gamma q (\cdot + \textbf{v} t+ D')S^{0,0}_{2k,0} \| _{L^1((0,T),L^p_x)}  \le C \gamma     \|  \langle D '(s)+s\textbf{v}\rangle ^{- k }    \|   _{L^1( 0,T) }
				,\end{aligned}  \end{equation*}
				\begin{equation} \label{eq:inteqD1} \begin{aligned} &  \|  \langle D '(s)+s\textbf{v}\rangle ^{- k }    \|   _{L^1( 0,T) } = \|  \langle D '(0)+s \textbf{v}   +   I(s)\rangle ^{- k }    \|   _{L^1( 0,T) } \text{ where} \\&  I(s):=sv_0+   \int _0^s   (\dot D ' (\tau )
				- v_0) d\tau  \quad , \quad   |\dot I(s)|     \le  3C_0\epsilon  ,
				 \end{aligned}  \end{equation}
				where $|\dot I(s)|     \le  3C_0\epsilon$  follows by  \eqref{eq:unds1}--\eqref{eq:unds2} and  by $|v_0|\lesssim \epsilon  $, see \eqref{eq:sizeindata1}.				
				Then \eqref{eq:movpot1} follows by Lemma \ref{lem:drift} below. \qed

\begin{lemma}
  \label{lem:drift}   Assume \eqref{eq:eqD0} and \eqref{eq:inteqD1}. Then for $\epsilon _0>0$ in \eqref{eq:sizeindata} sufficiently small,  we have for a fixed $c$
	
\begin{equation}\label{eq:drift1}
\begin{aligned} &\gamma \|  \langle D '(s)+s \mathbf{{v}} \rangle ^{- k }    \|   _{L^1( 0,T) }<c \epsilon  .   \end{aligned}
\end{equation}

\end{lemma}
				\proof   If  $|\langle D '(0)+s \mathbf{{v}} \rangle | \ge 6C_0\epsilon s$ for all $s>0$, then since $|  I(s)|     \le  3C_0\epsilon s$ by
   \eqref{eq:inteqD1}
     we get  $  \langle D '(s)+s\mathbf{{v}}\rangle     \sim  \langle  D '(0)+s \textbf{v} \rangle       $
for all $s>0$.
 Then \eqref{eq:drift1} follows from
				\eqref{eq:eqD0}.

		\noindent	 Suppose for an $s_0>0$ that  $|  D '(0)+s_0 \mathbf{v}|< 6C_0\epsilon s_0$ . Squaring this inequality and for
   $C_1= (6C_0)
				^2  |\mathbf{v}|^{-2} $ we get
			\begin{equation*}\label{eq:drift2}
\begin{aligned} &  |\mathbf{v}|^2  ( 1-C_1 \epsilon ^2)  s_0^2+2D '(0) \cdot  \textbf{v} s_0 +|D '(0) |^2<0.
  \end{aligned}
\end{equation*}	
			This implies    $ ( D '(0) \cdot  \mathbf{v} )^2  >    |D '(0) | ^2    \,  |\textbf{v}| ^2   ( 1-C_1 \epsilon ^2) $  for the discriminant
and
\begin{equation*}
\begin{aligned} &   D '(0) \cdot  \mathbf{v}   < -     |D '(0) |     \,  |\mathbf{v}|      \sqrt{ 1-C_1 \epsilon ^2 } .   \end{aligned}
\end{equation*}
This implies $D '(0)\neq 0$ and $\text{dist}_{S^2}(-\frac{D '(0)}{|D '(0) |},  \frac{\textbf{v} }{|\mathbf{v} |}) =O(\epsilon ^2).$   From    \eqref{eq:eqD0} we get
		\begin{equation*} \begin{aligned} &    \frac{\gamma }{|\mathbf{v}|}
		\|    \langle D '_1(0)    -\frac{D '(0)}{|D '(0) |}  s   \rangle   ^{-   k } \|   _{L^1( \R _+) }   =\frac{\gamma }{|\mathbf{v}|}
		\|   \langle | D '_1(0) |   -   s     \rangle   ^{-   k } \|   _{L^1( \R _+) }    <10 \epsilon
		    .
			 \end{aligned}  \end{equation*}	
			For $\epsilon _0>0$ in \eqref{eq:sizeindata} small,   we get $  \gamma  |\mathbf{v}|^{-1} < \kappa \epsilon$    for  $\kappa =20 / \| \langle t \rangle ^{-k}\| _{L^1( \R  ) }$.
			We have
					\begin{equation*}   \begin{aligned} &  \gamma  \|  \langle D '(s)+s\mathbf{v}\rangle ^{- k }    \|   _{L^1( 0,T) } \le    \gamma  |\mathbf{v}|^{-1}
				 \|  \langle  D _1'(0)+s +   I_1 ( {s}/ |\mathbf{v} |   )   \rangle ^{- k }    \|   _{L^1( 0,|\mathbf{v} |T) },	
			 \end{aligned}  \end{equation*}
			where    $ \frac d {ds} [I_1 ( {s}/ |\mathbf{v} |   )]  \le 3C_0\epsilon   / |\textbf{v} |$. We complete  the proof of   \eqref{eq:drift1}
 by
			\begin{equation*}   \begin{aligned} &
	 \gamma  \|  \langle D '(s)+s\mathbf{v}\rangle ^{-  k }    \|   _{L^1( 0,T) }  \le     \gamma  |\textbf{v}|^{-1}
		\| \langle   D _1'(0)/ |\mathbf{v} |+s +   I_1 ( {s}/ |\textbf{v} |   )  \rangle ^{- k }
			     \|   _{L^1( 0,|\mathbf{v} |T) }\\&
		\le    2\gamma  |\mathbf{v}|^{-1}    \| \langle t\rangle ^{- k }\| _{L^1(\R )}  < 40\epsilon \text{ for $3C_0\epsilon _0  / |\mathbf{v} |<1/2$.}  \qquad \qquad  \qquad  \qquad \qquad \qed \end{aligned}  \end{equation*}

\begin{lemma}\label{lem:derpsi}   Consider the function $\psi (\pi, \Pi ,\Pi (P(\pi )f))$, see \eqref{eq:psi}.
At $\Pi =\pi $
\begin{equation} \label{eq:derpsi}
 \begin{aligned}  &  \text{ we have }    \partial _\Pi  \psi = 0  \text{ and }   \partial _\pi  \psi = \resto  ^{0,0} (\pi , \Pi , \Pi (f), f)  .
\end{aligned}
\end{equation}
 \end{lemma}
 \proof At $\pi =\Pi$ we have
\begin{equation*}  \begin{aligned}  &    \partial _{\Pi _j} \psi  (\pi , \Pi ,\Pi  (P(\pi )f)) = \partial _{\Pi  _j} d (\Pi -  \Pi  (P(\pi )f) ) \\& + \partial _{\Pi  _j} \lambda  (\Pi -  \Pi  (P(\pi )f) )  \cdot ( \pi -   \Pi  (P(\pi )f) )  =0
\end{aligned}\end{equation*}
because of $\partial _{p_j}d(p)=-p\cdot \partial _{p_j} \lambda (p)$. The other formula in
\eqref{eq:derpsi} is also elementary.
\qed

 \begin{lemma}\label{lem:dotPi}   Under  the hypotheses of   Proposition \ref{prop:mainbounds} and  for $\epsilon _0>0$ in \eqref{eq:sizeindata} small enough we have $\| \dot \Pi  _j\| _{L^1(I) } \le c \epsilon $ for  a fixed $c$ for all  $j$.
\end{lemma}
\proof   We have $\dot \Pi _4 =\dot Q=0$.
For $q =\gamma  \partial _{ x_a}
   V $,  the  rhs of  \eqref{eq:mom} for $a\le 3$  is like
\begin{equation} \label{eq:derPipotinter1} \begin{aligned} &
     \langle    q (\cdot + \textbf{v} t+ D') S ^{0,0} ,    S ^{0,0}   \rangle    +
		\langle    q (\cdot + \textbf{v} t+ D') S ^{0,0} ,  e^{J \resto ^{0,2}\cdot \Diamond}  f  \rangle   \\& + \langle    q (\cdot + \textbf{v} t+ D')  e^{J \resto ^{0,2}\cdot \Diamond }  f  ,  e^{J \resto ^{0,2}\cdot \Diamond}  f  \rangle      .
				\end{aligned}  \end{equation}
				So  by Lemma \ref{lem:movpot} and \eqref{Strichartzradiation},  for a fixed $c_0$ (this yields the desired result)	 \begin{equation} \label{eq:derPipotinter2} \begin{aligned} &
			\|   \text{\eqref{eq:derPipotinter1}}\| _{L^1(0,T)}  \le \| \gamma   q (\cdot + \textbf{v} t+ D') S ^{0,0}\| _{L^1(0,T)L^1_x}
			 \|    S ^{0,0}\| _{L^\infty (0,T)L^\infty  _x} \\& +  \|  \gamma  q (\cdot + \textbf{v} t+ D') S ^{0,0}\| _{L^1(0,T)L^2_x}
			\|   f\| _{L^\infty (0,T)L^2_x} +C_q\|   f\| _{L^2(0,T)L^6_x}^2  \\& \le c _0 \epsilon + c_0C_0^2  \epsilon ^2  <2c_0 \epsilon \ .
		\qquad 	\qquad \qquad \qquad \qquad \qquad \qquad \qquad 	\qed\end{aligned}  \end{equation}

Lemmas \ref{lem:pfunds1} and  \ref{lem:vel} below are part of the proof of   Proposition \ref{prop:mainbounds}
but are also  used to prove  Theorem  \ref{thm:mainbounds}.
\begin{lemma}\label{lem:pfunds1}    Consider the hypotheses of   Proposition \ref{prop:mainbounds}.
 Then there exists a fixed constant   $c$ s.t.,  for $\epsilon _0>0$ in \eqref{eq:sizeindata}
 small enough,   \eqref{eq:unds1} holds for $C= c$.
 \end{lemma}
\proof
Consider \eqref{eq:eqD}.  By the notation of    Lem.  \ref{lem:ExpH11}  at
 $\Pi =\pi  (t)$   by Lem.   \ref{lem:derpsi}
\begin{equation*}  \begin{aligned}  &    \partial _{\Pi  }  H^{(\pi )} =\partial _{\Pi }  \psi +\partial _{\Pi }  H' _2 +\partial _{\Pi }\textbf{R}  =\partial _{\Pi }  H' _2 +\partial _{\Pi }\textbf{R} ,    \text{ in the  notation of \eqref{eq:ExpH11}}.
\end{aligned}\end{equation*}
   \eqref{Strichartzradiation}--\eqref{L^2discrete}   in $[0,T]$ and  Lemma \ref{lem:ExpH11}
give $|\partial _{\Pi _j}  H^{(\pi )}| \le \epsilon  $  for  $\epsilon _0 $ small.
   Then  $\|\partial _{\Pi _j}  H^{(\pi )}\|_{L^1(0,t)}\le t\epsilon $, which for $t\ge 1$ is the main term in \eqref{eq:unds1}.
 Consider
\begin{equation} \label{eq:derPipotinter} \begin{aligned} &    \gamma
    \partial _{\Pi _a}\langle    V (\cdot + \textbf{v} t+ D') (\Phi _{p'}+P(p')r') ,   \Phi _{p'}+ P(p')P(\pi )r'  \rangle
			\\&=\gamma   \langle    \partial _{D' _b}V (\cdot + \textbf{v} t+ D') (\Phi _{p'}+ P(p')P(\pi )r') ,   \Phi _{p'}+P(p')r'  \rangle     \partial _{\Pi _a} D' _b \\& + \gamma   \langle    V (\cdot + \textbf{v} t+ D') (\Phi _{p'}+ P(p')P(\pi )r') ,  \partial _{\Pi _a}  (\Phi _{p'}+ P(p')P(\pi )r') \rangle .
				\end{aligned}  \end{equation}
   \eqref{eq:tranR} yields $D'=D+\resto ^{0,2}$. So   $ \partial _{\Pi _a} D' _b= \resto ^{0,2} =O(\epsilon ^2 )$  by \eqref{Strichartzradiation}--\eqref{L^inftydiscrete}.  Since  the other factors  of the first term in the rhs are
$O(1)$,  this yields the desired bound $c\epsilon t$   for any given  $c>0$ on the $L^1(0,t)$ norm of  the first term in the rhs.

    \noindent  The  last line in  \eqref{eq:derPipotinter}  is the sum of   terms like \eqref{eq:derPipotinter1}, with $L^1(0,T)$ norm   $\le c\epsilon$,
\begin{equation*}   \begin{aligned} & \text{ and }   \gamma   \langle    V (\cdot + \textbf{v} t+ D') (\Phi _{p'}+ P(p')P(\pi )r') ,   P(p')P(\pi )\partial _{\Pi _a} r' \rangle  .
				\end{aligned}  \end{equation*}
Substituting \eqref{eq:diffothers6}, up to terms already bounded we get
\begin{equation*}  \begin{aligned} &    \gamma   \langle    V (\cdot + \textbf{v} t+ D') (\Phi _{p'}+ P(p')P(\pi )r') ,  \Diamond  f \rangle  .
				\end{aligned}  \end{equation*}
The  $L^1(0,T)$ of this term is bounded, for $\epsilon _0>0$ small enough,  by
\begin{equation*}  \begin{aligned} &   \gamma    \|    V (\cdot + \textbf{v} t+ D')  S^{0,0} \|  _{L^1_tL^2_x}  \|  f  \|  _{L^\infty _tH^1_x}    +\gamma   c_V \|  f  \|  _{L^2 _tL^6_x}  \|  f  \|  _{L^2 _tW^{1,6}_x}   \le \epsilon .
				\end{aligned}  \end{equation*}
 To complete the proof of Lemma \ref{lem:pfunds1}
 we need to   bound the second line of \eqref{eq:eqD}.
Using the fact that,  by Lemma \ref{lem:ODEdomains}, the functions
  $ \resto ^{0,2}      $ are preserved by
  the changes of coordinates in Lemma \ref{th:main}, by \eqref{eq:partham} we have
 \begin{equation*}  \begin{aligned}   &
 \{ \resto ^{0,2}      ,   K_0\} = \{ \resto ^{0,2}      ,   H ^{(\pi )}\}  ^{(\pi )} =\{ \resto ^{0,2}      ,   \psi  + H'_2+Z_0+Z_1+\resto\}  ^{(\pi )}
 . \end{aligned}  \end{equation*} Using the fact that
 \begin{equation} \label{eq:Pipot0} \begin{aligned}   &
 \text{$  \nabla _{\tau , \Pi ,z, \overline{z} }\resto ^{0,2}=\resto ^{0,1}$ and    $\nabla _{f}\resto ^{0,2}=\textbf{S} ^{0,1}$, }\end{aligned}  \end{equation}
 from \eqref{eq:poiss}
   it is elementary to get
 \begin{equation} \label{eq:Pipot1} \begin{aligned}   &
 \text{$ \{ \resto ^{0,2}      ,   K_0\}  = O(\epsilon ( |z| +\| f \| _{L^{2,-s}} ))$ for any preassigned $s>1$
 . }\end{aligned}  \end{equation}
 For example we have

 \begin{equation*} \begin{aligned} &    \langle  \nabla _f\resto ^{0,2}, (P _c(p_0 ) P_c (\pi )P _c (p_0 ) )^{-1}J \nabla _fE_P(f) \rangle  =  \langle  \mathbf{S}  ^{0,1},    \beta (|f|^2)f \rangle \\&    \sim \|  \langle x \rangle ^{-s}  \beta (|f|^2)f\| _{L^1_x}  \lesssim \| f \| _{L^{2,-s}} \|  \beta (|f|^2) \|  _{L^2_x} \le C_0 ^2 \epsilon ^2 \| f \| _{L^{2,-s}} .
\end{aligned}\end{equation*}
\eqref{eq:Pipot1} implies  $| \{ \resto ^{0,2}      ,   K_0\} | < \epsilon $.  Set $\langle  V  \rangle :=\langle  V (\cdot + \textbf{v} t) u, u\rangle $.
Then     by  \eqref{eq:poiss}
\begin{equation*}  \begin{aligned} &  \{ \resto ^{0,2}      ,   \langle  V  \rangle \}  = \partial _{\tau  _j}\resto ^{0,2} \partial _{\Pi _j}\langle  V  \rangle- \partial _{\Pi _j}\resto ^{0,2} \partial _{\tau  _j}\langle  V  \rangle -\im \partial _{z_j}\resto ^{0,2}\partial _{\overline{z}_j}\langle  V  \rangle
\\&+\im \partial _{\overline{z}_j}\resto ^{0,2}\partial _{z_j}\langle  V  \rangle +  \langle  \nabla _f\resto ^{0,2}, (P _c(p_0 ) P_c (\pi )P _c (p_0 ) )^{-1}J \nabla _f\langle  V  \rangle \rangle.
\end{aligned}\end{equation*}
 By Lemma \ref{lem:potaverage} we have
\begin{equation*}  \begin{aligned} &  \| \nabla _{z, \overline{z} }\langle  V  \rangle \|  _{L^\infty  _t} +\|
\nabla _f\langle  V  \rangle \| _{   L^\infty  _t L^2_x }  \le C  \|  z  \|  _{L^\infty  _t} +\|
f  \| _{L^\infty  _tH^1_x }  \le C C_0 \epsilon .
\end{aligned}\end{equation*}
By the arguments used to bound \eqref{eq:derPipotinter},  we have  $\| \nabla _{\Pi, \tau  }\langle  V  \rangle \|  _{L^\infty  (0,\infty )} \lesssim \epsilon .$    This and \eqref{eq:Pipot0} imply the following, which  completes the proof of Lemma \ref{lem:pfunds1},
\begin{equation} \label{eq:Pipot2} \begin{aligned}   &
 \text{$ \{ \resto ^{0,2}        ,\langle  V (\cdot + \textbf{v} t) u, u\rangle\}   = O(\epsilon ( |z| +\| f \| _{L^{2,-s}} ))$
 . }\end{aligned}  \end{equation}						
\qed

 \begin{lemma}\label{lem:vel}   Under  the hypotheses of   Prop. \ref{prop:mainbounds} and  for for $\epsilon _0>0$ in \eqref{eq:sizeindata}
 small enough,  \eqref{eq:unds2} holds for   $ C=c< C_0/2$ with $c$ a fixed constant. Furthermore,  \eqref{eq:orbstab1} holds for $ C=c< C_0/2$ with $c$ a fixed constant.
\end{lemma}
\proof  By Lemma \ref{lem:transf},   \eqref{eq:decomp2},  \eqref{Strichartzradiation} and  by \eqref{L^inftydiscrete}, for a fixed $c_1$ we have
\begin{equation} \label{eq:vel1} \begin{aligned} &   |Q(\Phi _{p'})  -Q|=|Q(\Phi _{p'})  -Q(\Phi _{p_0})| = Q(P(p')
r')\le c_1 C_0^2 \epsilon ^2\le  \epsilon
				.\end{aligned}  \end{equation}
Similarly, for $a\le 3$  we have
\begin{equation*} \label{eq:vel2} \begin{aligned} &   |\Pi _a(\Phi _{p'})  -\Pi _a |= |\Pi (P(p')
r')|\le c_1 C_0^2 \epsilon ^2\le  \epsilon
				.\end{aligned}  \end{equation*}
Furthermore
\begin{equation} \label{eq:vel3} \begin{aligned} &   | \Pi _a(\Phi _{p_0}) -\Pi _a | \le  \|\dot
\Pi _a \| _{L^1(0,t)} \le c      \epsilon
				 \end{aligned}  \end{equation}
by Lemma  \ref{lem:dotPi}.
Hence  for a fixed $c$ we have
\begin{equation} \label{eq:vel4} \begin{aligned} &   | \Pi _a(\Phi _{p_0}) -\Pi _a(\Phi _{p'}) | \le  c      \epsilon   .
				 \end{aligned}  \end{equation}
From  \eqref{eq:vel1} and  \eqref{eq:vel4} we get $|p'-p_0|\le c     \epsilon  $ for a fixed $c$.  This proves the  first sentence of   Lemma \ref{lem:vel}. The second sentence
follows  from the first by  formulas \eqref{eq:LagrMult} and the fact that the Lagrange multipliers
$\lambda (p)$ depend smoothly on $p$.
\qed

\begin{lemma}\label{lem:divcenter}
  Assume hypotheses  and conclusions of  Prop. \ref{prop:mainbounds}.
 Then

 \begin{equation} \label{eq:divcenter1}   | D'(t)+t {\mathbf{v}}|\ge t2^{-1}{|\mathbf{v}|} - | D'(0) |
 \end{equation}
\end{lemma}
\proof  We have
\begin{equation*}  \label{eq:divcenter2}\begin{aligned} &    D'(t)+t\textbf{v} = D'(0) +  t \textbf{v}+ tv_0
+\int _0^t (v'(s)-v_0) ds +  \int _0^t (D'(s) -v'(s) ) ds.
\end{aligned}
\end{equation*}
  Lemmas \ref{lem:pfunds1}-- \ref{lem:vel}  and    $|v_0|\lesssim \epsilon  $ yield (giving  \eqref{eq:divcenter1}
for   $\epsilon _0>0$ in \eqref{eq:sizeindata}
 small) $ |D'(t)-t\textbf{v} |\ge
 t (|\textbf{v}  |  -   3c   \epsilon )-  | D'(0) | . $ \qed

Lemmas  \ref{lem:partialmod}  and \ref{lem:limv}  together   imply the limits \eqref{eq:asstab2}.
\begin{lemma}\label{lem:partialmod}
  Assume hypotheses  and conclusions of  Prop. \ref{prop:mainbounds}.
 Then   we have

 \begin{equation} \label{eq:parasstab2}  \lim _{t\to +\infty} ( \dot D  '     -v ' )=0 \, , \quad   \lim _{t\to +\infty}    (\dot \vartheta   '  - \omega ' -  4^{-1}{(v')^2 }    )=0.
 \end{equation}

\end{lemma}
\proof
By  \eqref{eq:eqD} we have
\begin{equation}\label{eq:eqD1}  \begin{aligned} &    |\dot D   '  - v '|\le
   | {\partial   _{\Pi } } H^{(\pi )}| + |\{ \resto ^{0,2}      ,   K_0\} |+ 2^{-1}{\gamma }    |\{ \resto ^{0,2}        ,\langle  V (\cdot + \textbf{v} t) u, u\rangle\} |
       \\&
  +2^{-1}{\gamma}
    |\partial _{\Pi  }\langle    V (\cdot + \textbf{v} t+ D') (\Phi _{p'}+P(p')r') ,   \Phi _{p'}+P(p')r'  \rangle |
		 .\end{aligned}  \end{equation}
We now prove separately
   \begin{align} &    \label{eq:partialmod1}
  \lim _{t\nearrow \infty}  {\partial   _{\Pi } } H^{(\pi )} =0 \\& \label{eq:partialmod3}
  \lim _{t\nearrow \infty}
     \partial _{\Pi  }\langle    V (\cdot + \textbf{v} t+ D') (\Phi _{p'}+P(p')r') ,   \Phi _{p'}+P(p')r'  \rangle  =0
			 \\& \label{eq:partialmod21}  \lim _{t\nearrow \infty} \{ \resto ^{0,2}      ,   K_0\}  =0 \\& \label{eq:partialmod22}  \lim _{t\nearrow \infty} \{ \resto ^{0,2}        ,\langle  V (\cdot+\textbf{v} t) u, u\rangle\}   =0.\end{align}
		We have by \eqref{eq:partham}
\begin{equation*}  \begin{aligned}  &
|\partial _{\Pi _a}  H ^{(\pi )} | \le
|\partial _{\Pi _a} ( H' _2  + Z_0)|  +  |\partial _{\Pi _a} Z_1|  +|\partial _{\Pi _a} \resto | .
\end{aligned}\end{equation*}
  $|\partial _{\Pi _a} ( H' _2  + Z_0)|  \le C |z|^2$ by \eqref{L^inftydiscrete} and by \eqref{eq:asstab3}   this goes to 0 as $t\nearrow \infty .$ Also
\begin{equation*}  \begin{aligned}  &
|\partial _{\Pi _a} Z_1|\le C \sum _{\mathbf{e} \cdot \mu >\omega _0 }| z ^\mu  |  \|  f \| _{L^2_x }.
\end{aligned}\end{equation*}
By \eqref{Strichartzradiation}--\eqref{L^inftydiscrete}  $\partial _{\Pi _a} Z_1\in L^1(0,\infty )$. It is   easy to see that $\partial _{\Pi _a} Z_1$  has bounded time derivative and so converges to 0 for $t\nearrow \infty .$
Similarly, $\partial _{\Pi _a} \resto \in  L^1(0,\infty )$,    has bounded time derivative and so  goes to 0 for $t\nearrow \infty .$   So we get \eqref{eq:partialmod1}.

Turning to   \eqref{eq:partialmod3}      notice that  the function in  \eqref{eq:partialmod3} is of   form \eqref{eq:derPipotinter1}.
 We have
			\begin{equation*} \label{eq:partialmod4} \begin{aligned} &
			 |   \text{\eqref{eq:derPipotinter1}} |    \le \|   q (\cdot + \textbf{v} t+ D') S ^{0,0}\| _{ L^1_x}
			 \|    S ^{0,0}\| _{ L^\infty  _x} +  \|   q (\cdot + \textbf{v} t+ D') S ^{0,0}\| _{ L^2_x}
			\|   f\| _{ L^2_x} \\& +| \Upsilon (t) |\text{  with }
			  \Upsilon  (t):=
	  \langle    q (\cdot + \textbf{v} t+ D')  e^{J \resto ^{0,2}\cdot \Diamond }  f  ,  e^{J \resto ^{0,2}\cdot \Diamond}  f  \rangle
				\end{aligned}  \end{equation*} 	
The   first line   $ \le C_S (1+ \|   f\| _{ L^2_x} ) \langle  D'(t)+ \textbf{v} t
\rangle ^{-2}  $ which by Lemma \ref{lem:divcenter} goe to 0 as $t \nearrow \infty$
  To finish   with  \eqref{eq:partialmod3}   we show
 $
  \lim _{t\nearrow \infty}
 \Upsilon (t)  =0  .$
First of all we have
	\begin{equation*} \label{eq:partialmod6} \begin{aligned} &
		 \| \Upsilon  \| _{L^1(0,\infty )} \le
			C_q\|   f\| _{L^2((0,\infty ) ,L^6_x)}^2 <\infty.
				\end{aligned}  \end{equation*} 	
 Furthermore,  it is simple to check that $\Upsilon  (t)$  is differentiable with   $ \dot \Upsilon   \in  L^\infty (0,\infty )$.
Then, as in Lemma \ref{lem:proptheor1}  we conclude that  $
  \lim _{t\nearrow \infty}
 \Upsilon (t)  =0  .$

We turn now to     \eqref{eq:partialmod21} and    \eqref{eq:partialmod22}. By \eqref{eq:Pipot1} and \eqref{eq:Pipot2}
it is enough to check the limit as $t\nearrow \infty $ of
$ O(  |z| +\| f \| _{L^{2,-s}} )$. By Lemma
\ref{lem:proptheor1} we know that $z(t)$ converges to 0.
With standard arguments we can see that for  any $u_0\in L^2$    and for any function  $ \tau (t)$
		we have  for  $s>3/2$
		\begin{equation} \label{eq:w-conv}
\begin{aligned} &  \lim _{t\nearrow  \infty}  e^{J \tau (t) \cdot \Diamond }e^{-tJ\Delta }u_0=0   \text{  in $L^{2,-s}(\R ^3)$}.
\end{aligned} \end{equation}
Then \eqref{eq:scattering1} implies  $\displaystyle \lim _{t\nearrow \infty} \| f (t)\| _{L^{2,-s}}=0 .$ This yields
\eqref{eq:partialmod21}--\eqref{eq:partialmod22}.

\qed

  Then we get the following result which, with Lemma \ref{lem:partialmod}, yields \eqref{eq:orbstab1}--\eqref{eq:asstab2} thus completing the proof
of  Theorem \ref{thm:mainbounds}.

\begin{lemma}\label{lem:limv}    There exist $\omega _+$ and $v_+$    such that the limits  \eqref{eq:asstab1}
  are
true.
\end{lemma}
\proof   By \eqref{eq:variables}, for  $r'=  e^{J\resto ^{0,2}\cdot \Diamond }  ( r
  +\mathbf{S} )$,
\begin{equation} \label{eq:limv1}
\begin{aligned} &
  \Pi _j(t)= {p'_j(t)} +  \Pi _j(  r'(t))   - \Pi _j( (P (p' (t))-P (p_0  )) r'(t)) \\& +\langle r'(t),
 \Diamond _j (P (p' (t))-P (p_0  )) r'(t)\rangle   ,
\end{aligned} \end{equation}
 $
  |\resto ^{0,2}|   + \|\mathbf{S}  \| _{H ^{1, s}} ^{2}  \le C  (|z|+\| f \| _{L^{2,-s}})^2 ,$
for a preassigned large $s>1$.
We   have
\begin{equation}\label{eq:limv11}
\begin{aligned} & \Pi _j (   r' ) =\Pi  _j(  r  )   + \Pi _j(  \mathbf{S} )     +\langle r ,
 \Diamond _j \mathbf{S} \rangle
\end{aligned} \end{equation}
with  by \eqref{eq:decomp2}
\begin{equation}\label{eq:limv10}
   r  = (P(\Pi ) P(p_0))^{-1}  ( z \cdot \xi
     ^{(\Pi )} + \overline{z} \cdot \overline{\xi  }  ^{(\Pi )}  + P_c(\Pi )f ) = f+\textbf{ S}^{0,1} (\Pi , z,f). \end{equation}
		 By  \eqref{eq:scattering1} and  \eqref{eq:w-conv} which  imply  $\displaystyle \lim _{t\nearrow \infty} \| f (t)\| _{L^{2,-s}}=0  $ and by \eqref{eq:asstab3}
	 we have  $\displaystyle  \lim _{t\nearrow  \infty} \textbf{ S}^{0,1} (\Pi , z,f) =\lim _{t\nearrow  \infty} \textbf{ S}   =0$  in $H^1$ and  by Lemma \ref{lem:scattf}  and \eqref {eq:limv11}--\eqref {eq:limv10}
\begin{equation}\label{eq:limv12}
   \lim _{t\nearrow  \infty} \Pi  _j(  r'  )  =\lim _{t\nearrow  \infty} \Pi  _j(  r  )  = \lim _{t\nearrow  \infty} \Pi  _j(  f  )  =\Pi  _j(  f _+ )   . \end{equation}
Furthermore  by   \eqref{eq:limv11}--\eqref{eq:limv10} we get
\begin{equation}\label{eq:limv13}
\begin{aligned} &   |\Pi  (   r'(t) )  - \Pi  (  f _+ )  |\le |\Pi  (  f(t) )  - \Pi  (  f _+ )  | + C  (|z|+\| f \| _{L^{2,-s}})^2
\le C \epsilon   .
\end{aligned} \end{equation}
By   $\|\dot \Pi \| _{L^1(\R)}\le c \epsilon $, see  Lemma \ref{lem:dotPi}, we get   that $\displaystyle \lim _{t\nearrow \infty}\Pi (t)=\Pi  ^+$ exists and $|\Pi (t)-\Pi  ^+|\le c \epsilon .$
Then  taking the limit in  \eqref{eq:limv1} we get
\begin{equation} \label{eq:limv3}
\begin{aligned} &
  \Pi _j ^+= \lim _{t\nearrow  \infty} {p'_j(t)} +  \Pi _j(  f_+)   -  \lim _{t\nearrow  \infty} \Pi _j( (P (p' (t))-P (p_0  )) r'(t)) \\& + \lim _{t\nearrow  \infty} \langle r'(t),
 \Diamond _j (P (p' (t))-P (p_0  )) r'(t)\rangle  =   \lim _{t\nearrow  \infty} {p'_j(t)} +  \Pi _j(  f_+).
\end{aligned} \end{equation}
So   $\displaystyle  \lim _{t\nearrow  \infty} {p' (t)} = p^{\prime +}$  exists  and hence also  the limits in  \eqref{eq:asstab1}.

 \qed

\section{Proof of Proposition \ref{prop:mainbounds} }
\label{sec:pfprop}
In Lemmas \ref{lem:pfunds1} and  \ref{lem:vel} we have already proved that if   \eqref{Strichartzradiation}-- \eqref{L^inftydiscrete} and   \eqref{eq:unds1}--\eqref{eq:unds2} hold for any  given $C=C_0>c_0$ in $[0,T]$, then,
for $\epsilon _0  >0$ sufficiently small,
   \eqref{eq:unds1}--\eqref{eq:unds2}  hold for   $C=c$  with  $c$ fixed  which can be taken $c<c_0/2$.  Hence we can say that
   \eqref{eq:unds1}--\eqref{eq:unds2} hold for  $C=C_0/2$.   Now we need  to prove that also
\eqref{Strichartzradiation}-- \eqref{L^inftydiscrete} hold   for  $C=C_0/2$.
The following argument was initiated in \cite{BP2,SW3}.
We first   get estimates on  $f $  from the estimates on $z$.
 When   $\gamma V(x+\textbf{v}t)$    is absent in \eqref{eq:NLSP} and we are in the
situation of \cite{Cu3}, the main ingredients of the first step are  the estimates on the groups $e^{\im t \mathcal{K}_\omega }$ in \cite{Cu6,CPV},   see also    \cite{schlag}, and an integrating factor
argument as in \cite{beceanu}.  The additional presence of   $\gamma V(x +\textbf{v}t)$
is treated here following the theory in \cite{RSS1}  on  {\it charge transfer models}, see also \cite{cai}.

The second step consists in observing that the part of the Hamiltonian in \eqref{eq:partham} which effectively  couples
discrete and continuous modes is $Z_1$. In particular this yields the main source term in the equation of the
continuous mode $f$ and the main coupling term in the equations of the discrete modes.

The   last step of   the  proof of  Proposition \ref{prop:mainbounds} is where
the   approximation of our equations by    linear equations breaks down. In fact
the linear approximate equation of the discrete modes does not yield decay. The key ingredient
yielding decay is instead the loss of energy of the discrete modes due to the nonlinear   coupling  with the continuous mode $f$,
by means of the Fermi Golden Rule.

The three steps of the proof  of  Proposition \ref{prop:mainbounds} are split in
  Sections \ref{subsec:cont mode}, \ref{subsec:coupling} and \ref{subsec:FGR}.
The linear estimates needed   in  Sect. \ref{subsec:cont mode} are proved in Sect. \ref{sec:dispersion}.

It is important     to observe that the elaborate procedure we have just outlined, is necessary
only on large intervals of time  $[0,T]$, say for $T>\epsilon _0^{-1}$, since for in intervals with  $T\le \epsilon _0^{-1}$, the proof of Proposition \ref{prop:mainbounds} is much simpler and does not need  the Fermi Golden Rule,
see Remark \ref{rem:tspan1}  and the beginning of Sect. \ref{subsec:FGR}.

\subsection{Analysis of the continuous mode}
\label{subsec:cont mode}

Our aim is now to show the following lemma.
\begin{lemma}\label{lem:conditional4.2} Assume
the hypotheses of Prop. \ref{prop:mainbounds}.    Then  there is a fixed $c$
such that for all admissible pairs $(p,q)$
\begin{equation}
  \|  f \| _{L^p_t([0,T],W^{ 1 ,q}_x)}\le
  c  \epsilon  + c   \sum _{\textbf{e}\cdot \mu >\omega _0 }| z ^\mu  |  _{L^2_t( 0,T  )}
  \label{4.5}
\end{equation}
where we sum only on multiindexes such that $\textbf{e}\cdot \mu -
\textbf{e} _j <\omega _0$ for any  $j$ such that
for the $j$--th component of $\mu $ we have   $\mu _j\neq 0$.
\end{lemma}

\begin{remark}
\label{rem:tspan1}   For  $T\le \epsilon _0^{-1}$,   \eqref{L^inftydiscrete}
in  $[0,T]$,  \eqref{4.5} and  the fact that $|\mu |\ge 2$  in  \eqref{4.5}  imply   $| z ^\mu  |  _{L^2_t( 0,T  )}+ \|  f \| _{L^p_t([0,T],W^{ 1 ,q}_x)}\le
  2c  \epsilon   $  yielding  \eqref{Strichartzradiation}--\eqref{L^2discrete} for $C=C_0/2$  since we are free to choose $C_0>4c$.
\end{remark}

\proof We rewrite \eqref{eq:eqf}
\begin{equation} \label{eq:eqf1} \begin{aligned} &            \dot f =  J \nabla _{f}K' +   \dot \Pi _b \partial _{\pi _b} f + \mathbf{P}(p,p_0)J \nabla _{f}K' \\&  \mathbf{P}(p,p_0):=
   (P _c(p_0 ) P_c (\pi )P _c (p_0 ) )^{-1}   - P _c (p_0 )    .\end{aligned}  \end{equation}
	 We have  for any $k\in \N$
	\begin{equation} \label{eq:eqP1} \begin{aligned} & \mathbf{P}(p,p_0)=
   (P _c(p_0 )  +  P _c(p_0 ) (    P _d(p_0 ) -P_d (\pi )   )P _c (p_0 ) )^{-1} - P _c (p_0 )  \\& \|     \mathbf{P}(p,p_0)P _c(p_0 )  \| _{\Sigma _{-k}\to  \Sigma _{ k}} \le C_k |p_0 - \pi | \le  C_k  c \epsilon ,
	 \end{aligned}  \end{equation}
	where the  last inequality is  obtained by  \eqref{eq:vel3}   and $Q(\Phi _{p_0})= Q(\Phi _{\pi})$.

 \noindent By Lemma  \ref{lem:reg5}   we have for a fixed $c$
 \begin{equation} \label{eq:eqf11} \begin{aligned} &             \dot \Pi _b \partial _{\pi _b} f  =\dot \Pi _b   J\resto ^{0,1}\cdot \Diamond f+\dot \Pi _b\textbf{S} ^{0,1} \\&
   \|   \dot \Pi _b\textbf{S} ^{0,1} \|  _{L^1((0,T), H^1) }+ \|   \dot \Pi _b\resto ^{0,1} \|  _{L^1 ( 0,T) } \le c C_0  \epsilon ^2  .\end{aligned}  \end{equation}
Recall now

 \begin{equation*} \label{eq:eqf12} \begin{aligned} &   \nabla _{f}K'=\nabla _{f} H^{(\pi )}
  + \frac{\gamma }{2} \nabla _{f}
     \langle      V(\cdot +\textbf{v}t) u , u \rangle
		  .\end{aligned}  \end{equation*}
We have
\begin{equation} \label{eq:nablaH} \begin{aligned} &   \nabla _{f} H^{(\pi )}
   =  \nabla _{f} \psi + \nabla _{f}H'_2+ \nabla _{f} Z_0 + \nabla _{f} Z_1 +  \nabla _{f} \resto
		  .\end{aligned}  \end{equation}
We have
\begin{equation} \label{eq:nablaH1} \begin{aligned} &
     \nabla _{f} (\psi +  H'_2 +Z_0 )=-\im  P_c^*(p_0)\mathcal{H}_{\pi
 } P_c(\pi ) f   + {A}\cdot    P_c^*(p_0) \Diamond f,
 \\&  {A}:= \nabla _{\Pi (f)} \psi (\Pi (f))
     +\sum _{\substack{ |\mu +\nu |\ge 2\\
\textbf{e}  \cdot (\mu -\nu )=0}}
 \nabla _{\Pi (f)}  a_{\mu \nu} ( \pi , \Pi , \Pi  (f) )  z^\mu
\overline{z}^\nu
		  .\end{aligned}  \end{equation}
			Notice that by \eqref{L^inftydiscrete} we have $A\in L^\infty( (0,T)  \R ^4)$ with,  for $\mathcal{A}=A$,
			
			\begin{equation} \label{eq:nablaH2} \begin{aligned} &
     \|  \mathcal{A} \| _{L^\infty( (0,T), \R ^4)}\le c C_0^2 \epsilon ^2
		  .\end{aligned}  \end{equation}

\noindent Furthermore,  by standard arguments that we skip we have
\begin{equation} \label{eq:nablaH3} \begin{aligned} &
     \nabla _{f} \resto =  \mathcal{{A}}_{\resto}\cdot    P_c^*(p_0) \Diamond f + R_1+ R_2 \text{ with for any preassigned $s$}
\\& \| R_1 \| _{L^1_t([0,T],H^{ 1 }_x)}+\|
   R_2 \| _{L^{2
 }_t([0,T],H^{ 1
,
 s}_x)}\le C( s,C_0
 ) \epsilon^2    \end{aligned}  \end{equation}
 and with $\mathcal{{A}}=\mathcal{{A}}_{\resto}$ satisfying \eqref{eq:nablaH2}.
 We  have
  \begin{align}   \label{eq:nablaZ1} &
     \nabla _{f} Z_1=  \im  J   \sum _{| \textbf{e}   \cdot(\mu-\nu)|>\omega   _0  }
z^\mu \overline{z}^\nu    G_{\mu \nu}( t ,\Pi (f)   )  +{B}\cdot    P_c^*(p_0) \Diamond f \ , \ B_j:= \\&  \im  \sum _{|\textbf{e}   \cdot(\mu-\nu)|>\omega   _0  }  z^\mu \overline{z}^\nu \langle  J  \nabla _{\Pi _j(f)}G_{\mu \nu}( t ,\Pi (f)   ),f\rangle
		  \text{ , $G_{\mu \nu}( t, \Pi (f)   ):=G_{\mu \nu}( \Pi ,\Pi ,\Pi (f)   )$.}\nonumber \end{align}
			    \eqref{eq:nablaH2}  is true for  $\mathcal{A}=B$.
  Summing up and by  Lemma  \ref{lem:disppot}  below, we get
			  \begin{align}   \label{eq:eqf20}  &            \dot f =-\im    P_c (p_0)\mathcal{H}_{\pi
 } P_c(\pi ) f  +  J \gamma P_c(p_0)   V (\cdot + \textbf{v} t+ D' +\resto ^{0,2} )
	f \\&+   JA'(t)\cdot    P_c (p_0)  \Diamond f  -\im      \sum _{|\textbf{e} \cdot(\mu-\nu)|>\omega   _0  }
z^\mu \overline{z}^\nu    G_{\mu \nu}( t ,\Pi (f)   ) +R_1'+ R_2   \text{ with}  \nonumber  \\& \label{eq:nablaH44}
     \|  A'(t)\| _{L^\infty( (0,T), \R ^4) + L^1( (0,T), \R ^4)}+    \| R_1' \| _{L^1 ([0,T],H^{ 1 } )} \le c \epsilon  ,
		 \end{align}
 		for a fixed $c$.
		We now write
			\begin{equation} \label{eq:eqf21} \begin{aligned} &         P_c(\pi ) = P_c (p_0) + (P_d (p_0) - P_d (\pi )  ) \,   ,\\&    \mathcal{H}_{\pi }  =  \mathcal{H}_{p_0}+\im J (\lambda (\pi ) - \lambda (p_0) ) \cdot \Diamond +\im J (\nabla ^2 E_P(\Phi _{\pi }  ) -\nabla ^2 E_P(\Phi _{p_0}  ) ) .
		  \end{aligned}  \end{equation}
		Then by   $Q=Q(\Phi _{p_0 } )$ and \eqref{eq:vel3}, we have
  \begin{align}
  \nonumber &
         \|  (\nabla ^2 E_P(\Phi _{\pi }  ) -\nabla ^2  E_P(\Phi _{p_0}  ) )   f \| _{  L^2((0,T), W^{1, \frac 65}) }+ \\&
            \nonumber  \|  (P_d (p_0) - P_d (\pi ))  f \| _{  L^2((0,T), W^{1, \frac 65}) }
            \lesssim   |\pi - p_0| _{L^\infty (0,T)}   \|     f \| _{  L^2((0,T), L^{6} )}   \lesssim C_0 \epsilon  ^2 , \\ &   \|   \lambda (\pi ) - \lambda (p_0) \|   _{L^\infty( (0,T), \R ^4) }  \le c \epsilon .    \nonumber
		  \end{align}
		Then   we get to
			\begin{equation} \label{eq:eqf22} \begin{aligned} &          \im   \dot f =  \mathcal{H}_{p_0
 }   f  +\im J\gamma  P_c(p_0)   V (\cdot + \textbf{v} t+ D' +\resto ^{0,2} )
	f +\im J  A'(t)\cdot    P_c (p_0) \Diamond f\\&+     \sum _{| \textbf{e}   \cdot(\mu-\nu)|>\omega   _0  }
z^\mu \overline{z}^\nu  G_{\mu \nu}( t ,\Pi (f)   ) +E_1' + E_2'  ,    \end{aligned}  \end{equation}
		  $E_2 '$  (resp. $E_1 '$) satisfying  the same estimates of $R_2$  (resp. $R'_1 $).  We now consider
		\begin{equation} \label{eq:eqf23} \begin{aligned} &        \widetilde{ h}\in L^2_c(\mathcal{H}_{\omega _0
 } ) \to  e^{-J \frac{v_0\cdot x }2}\widetilde{h}=:f\in L^2_c(\mathcal{H}_{p _0
 } )   ,  \end{aligned}  \end{equation} which is an
		isomorphism.
Using \eqref{eq:decomp22}--\eqref{eq:decomp23} and
\eqref{eq:charge1} we have

\begin{equation*}  \begin{aligned} &       e^{ J \frac{v_0\cdot x }2}    P_c(p_0)    e^{ -J \frac{v_0\cdot x }2}  =
P_c(\omega _0) \ ,  \quad
e^{ J \frac{v_0\cdot x }2}   P_c(p_0)   V
	 e^{ -J \frac{v_0\cdot x }2} \widetilde{h}     = P_c(\omega _0)   V
  \widetilde{h}   \ ,  \\&   Q (f) = Q (h ) \ , \quad    \Pi _a (f )= \Pi _a(h ) +   \frac{v_{0a}}2 Q (h
   )  \text{ for $a\le 3$} \ ,\\&  e^{ J \frac{v_0\cdot x }2}  \Diamond  _a f =
	e^{ J \frac{v_0\cdot x }2} J \partial  _{x_a} f = \Diamond  _a \widetilde{h} +  \frac{v_{0a}}2\widetilde{h}=\Diamond  _a \widetilde{h} +  \frac{v_{0a}}2\Diamond  _4\widetilde{h}    \text{ for $a\le 3$}
	,\end{aligned}  \end{equation*}
 with $v_{0a}$ the components of $v_{0 }$.
Then    we get to
			\begin{equation} \label{eq:eqh22} \begin{aligned} &        \im     \dot {\widetilde{h}}  =\mathcal{H}_{\omega _0
 }   \widetilde{h}  +\im J\gamma  P_c(\omega _0)   V (\cdot + \textbf{v} t+ D' +\resto ^{0,2} )
	\widetilde{h} +
\im J\textbf{A}(t)\cdot    P_c ( {\omega _0
 }  ) \Diamond  \widetilde{h}
  \\& +
      \sum _{|\textbf{e}    \cdot(\mu-\nu)|>\omega   _0  }
z^\mu \overline{z}^\nu  e^{J \frac{v_0\cdot x}2 } G_{\mu \nu}(  t,\Pi (f)   ) +R_1 ^{\prime\prime }+ R_2 ^{\prime\prime }    \text{  where:}   \end{aligned}  \end{equation}
     $\textbf{A}_a=A'_a$  for $a\le 3$  and   $\textbf{A}_4=A'_4 +\frac 12 A'\cdot v_0$, so
$\textbf{A}$	 satisfies the estimates of $A'$   in \eqref{eq:nablaH44};  $R_1^{\prime\prime }$
		and $ R_2^{\prime\prime }$ satisfy the estimates of  $R_1^{\prime  }$    in \eqref{eq:nablaH44}
		and $ R_2 $    in \eqref{eq:nablaH3}. Set
\begin{equation*}
  h=M^{-1}\widetilde{h} \text{ with $M$ as in \eqref{eq:Homega1}.}
\end{equation*}
Then,  for $\mathcal{K}_{\omega
 }  =  M^{-1} \mathcal{H}_{\omega
 }    M$, see \eqref{eq:Homega1},   we have
\begin{equation} \label{eq:eqh222} \begin{aligned} &        \im     \dot { {h}}  =\mathcal{K}_{\omega _0
 }    {h}  + \sigma _{3}\gamma  P_c(\mathcal{K}_{\omega _0
 }  )   V (\cdot + \textbf{v} t+ D' +\resto ^{0,2} )
	 {h} +
\sigma _3\textbf{A}_4(t)     P_c (\mathcal{K}_{\omega _0
 })  {h}\\& - \sum _{a=1}^{3}\im \textbf{A}_a(t)     P_c (\mathcal{K}_{\omega _0
 }) \partial _{x_a} {h}+
      \sum _{|\textbf{e}    \cdot(\mu-\nu)|>\omega   _0  }
z^\mu \overline{z}^\nu \mathbf{G}_{\mu \nu}(  t,\Pi (f)   )   +R_1 ^{\prime\prime }+ R_2 ^{\prime\prime }     \, ,   \end{aligned}  \end{equation}
where:
\begin{equation} \label{eq:eqh223} \begin{aligned} &     \mathbf{G}_{\mu \nu}(  t,\Pi (f)   ):= M^{-1} e^{J \frac{v_0\cdot x}2 } G_{\mu \nu}(  t,\Pi (f)   ) \, , \quad f=e^{-\frac 12 J v_0 \cdot x}   Mh \ ,  \end{aligned}  \end{equation}
with $G_{\mu \nu}(  t,\Pi (f)   )$ the coefficients of $Z_1$, see \eqref{eq:partham} and Lemmas
\ref{lem:ExpH11} and \ref{th:main}.
We recall also that  the coefficients of  \eqref{eq:eqh222} satisfy
 \begin{equation} \label{eq:nablaH444} \begin{aligned} &
     \| \textbf{ A}  \| _{L^\infty(  0,T)   + L^1( 0,T) }+    \| R_1^{\prime\prime } \| _{L^1 ([0,T],H^{ 1 } )} +    \| R_2^{\prime\prime } \| _{L^2 ([0,T],H^{ 1 ,s} )} \le c \epsilon  .
		  \end{aligned}  \end{equation}  Then  from Theorem \ref{thm:strich} below we get
the following, which implies  \eqref{4.5}:
			\begin{equation}
  \| h \| _{L^p_t([0,T],W^{ 1 ,q}_x)}\le
  C  \epsilon  + C   \sum _{\textbf{e}   \cdot \mu >\omega _0 }| z ^\mu  | ^2_{L^2_t( 0,T  )} . \qquad    \qquad \qquad \qquad\qquad \qed
  \nonumber
\end{equation}

\begin{lemma}
  \label{lem:disppot}    We have the following formulas for a fixed $c$:
\begin{align} &
\nabla _f  \langle  V(\cdot +  {\mathbf{v}} t) u
,  u  \rangle  =    P_c(p_0)   V (\cdot +\mathbf{v} t+ D' +\resto ^{0,2} )
	f  + \mathcal{C}(t) \cdot \Diamond f + R_3\label{eq:disppot1} \\& \label{eq:nablaH5}
     \|  \mathcal{C}\| _{L^1( [0,T], \R ^4)}\le c   \epsilon    \quad , \quad     \| R_3 \| _{L^1 ([0,T],H^{ 1 }_x)} \le c \epsilon .
		  \end{align}
\end{lemma}
	 \proof
	     We denote the coefficient on the second line of    \eqref{eq:potaverage1} by
	\begin{equation*}   \begin{aligned} &
     \mathcal{C}:= -\gamma
 \langle   V (\cdot+ \mathbf{v} t+ D') (\Phi _{p'}+P(p')P(\pi )r'), \mathbf{S}^{0,0}   \rangle .
		  \end{aligned}  \end{equation*}
 For symmetry reasons, $ \mathcal{C}$ is an $\R ^4$ valued vector.    $ \|  \mathcal{C}\| _{L^1( [0,T], \R ^4)}\le c   \epsilon$   by the argument in Lemma
\ref{lem:dotPi}.   This settles    the  term $\mathcal{C}(t) \cdot \Diamond f $   in  \eqref{eq:disppot1}.

The   term on the rhs of the  first line  of   \eqref{eq:potaverage1} is of the form $ \varphi \textbf{S}^{0,1}$
with $ \varphi$ of the form \eqref{eq:derPipotinter1}. By estimate  \eqref{eq:derPipotinter2} this term can be
absorbed in $R_3 $. The same is true for the term in the third line of   \eqref{eq:potaverage1}.

\noindent  We now focus on the fourth line of    \eqref{eq:potaverage1}.
By
  Lemma \ref{lem:movpot} we have  for a fixed $c$
	\begin{equation} \label{eq:nablaH6} \begin{aligned} &
  \gamma    \|    (P_c(p_0) e^{-J \resto ^{0,2}\cdot \Diamond }+\gimel )   V (\cdot + \mathbf{v} t+ D')  \Phi _{p'} \| _{L^1( (0,T),  H^{1,s} )}\\&  \le  c' \gamma    \|     V (\cdot + \mathbf{v} t+ D')  \Phi _{p'} \| _{L^1( (0,T),  H^{1,s} )}   \le c   \epsilon  .
		  \end{aligned}  \end{equation}

\noindent  Next, we focus on the term
\begin{equation} \label{eq:nablaH7} \begin{aligned} & \gamma
    \gimel   V (\cdot + \textbf{v} t+ D')  P(p')P(\pi )r'  \\& = \gamma     \gimel   V (\cdot + \textbf{v} t+ D')  P(p')P(\pi )
	(	e^{J\resto ^{0,2}\cdot \Diamond }  r+   \textbf{S}^{1,1})
		  \end{aligned}  \end{equation}
The term with   $\textbf{S}^{1,1}$ is like in  \eqref{eq:nablaH6}.
Since  $
    [P(\pi ) ,
 	e^{J\resto ^{0,2}\cdot \Diamond }]  r =\textbf{S} ^{0,3}
		  $, see in Lemma 4.1 \cite{Cu0},
 and using  the expansion \eqref{eq:decomp2}, it is elementary to see that   \eqref{eq:nablaH7}  is the sum
of   terms  like in  \eqref{eq:nablaH6} and of
$$
    \gimel  \gamma  V (\cdot + \textbf{v} t+ D')  P(p') e^{J\resto ^{0,2} \cdot \Diamond }P(\pi )
	f.$$
Then
\begin{equation*} \label{eq:nablaH91} \begin{aligned} & \gamma
  \|   \gimel   V (\cdot + \textbf{v} t+ D')  P(p') e^{J\resto ^{0,2} \cdot \Diamond }P(\pi )
	f  \|  _{L^1 _tH^1_x}     \\& \lesssim\gamma   \|       V (\cdot + \textbf{v} t+ D')  P(p') e^{J\resto ^{0,2} \cdot \Diamond }P(\pi )
	f  \|  _{L^1 _tL^{2,-s}_x} \\&    \lesssim   \gamma  \|     \langle \cdot \rangle ^{-s}  V (\cdot + \textbf{v} t+ D')   \|  _{L^1 _tL^{2 }_x}
	\|
	f  \|  _{L^\infty  _tL^{2 }_x} \le cC_0\epsilon ^2.
		  \end{aligned}  \end{equation*}
Similarly, for $(*)$ a sum
of   terms  like in  \eqref{eq:nablaH6}
 \begin{equation} \label{eq:nablaH10} \begin{aligned} &   P_c(p_0)
e^{-J\resto ^{0,2}\cdot \Diamond}
 V (\cdot +  \mathbf{{v}} t+ D')  P(p')P(\pi )r'\\&= P_c(p_0)
e^{-J\resto ^{0,2}\cdot \Diamond}
 V (\cdot + \mathbf{{v}} t+ D')  P(p')P(\pi )
	(	e^{J\resto ^{0,2}\cdot \Diamond }  r+   \textbf{S}^{1,1})
\\&=
    P_c(p_0)   V (\cdot + \mathbf{{v}} t+ D' +\resto ^{0,2} )
	f  +(*).
		  \end{aligned}  \end{equation}
\qed

\begin{lemma}
  \label{lem:disppot1}    $\exists$ fixed $c>0$ s.t.
     $\|  \partial _{z_j}  \langle  V(\cdot + \mathbf{{v}} t) u
,  u  \rangle \| _{L^1( (0,T), \C  )}\le c    \epsilon $
  for  $i=1,2$.
\end{lemma}
\proof    We have  by \eqref{eq:potaverage2}
\begin{equation*}   \begin{aligned}
     &  {2}^{-1} \|
\partial _{z_j} \langle   \gamma V (\cdot + \mathbf{{v}} t ) u
,  u  \rangle  \|_{L^1_t}\le  \|
  \langle   V_2 (t)(\Phi _{p'}+P(p')P(\pi )r'), \mathbf{S}^{0,0}  \rangle    \|_{L^1_t} \\& + \| \langle  V_2 (t) (\Phi _{p'}+P(p')P(\pi )r') , e^{J\resto ^{0,2}\cdot \Diamond} J  \Diamond r \rangle   \|_{L^1_t}  \| \resto ^{0,1}  \|_{L^\infty _t}
\\&  +
 \| \langle  \partial _{D'}V_2 (t) (\Phi _{p'}+P(p')P(\pi )r'),
\Phi _{p'}+P(p')P(\pi )r'\rangle  \|_{L^1_t}  \| \resto ^{0,1}  \|_{L^\infty _t} .\end{aligned}  \end{equation*}
  $\| \resto ^{0,1}  \|_{L^\infty _t}\lesssim C_0 \epsilon$;    other factors are   $\le c\epsilon$
for a fixed $c$, see proof Lemma \ref{lem:dotPi}.
\qed

\subsection{Effective    discrete and continuous coupling}
\label{subsec:coupling}

For $\textbf{D}$ defined by the formula below,
we rewrite  \eqref{eq:eqh22} in the form
	   \begin{align} &  \nonumber       \im     \dot h  =\mathcal{K}_{\omega _0
 }   h  +\sigma _3\gamma  P_c(\mathcal{K}_{\omega _0
 } )   V (\cdot +\overbrace{ \textbf{v} t+ D' +\resto ^{0,2}}^{  \textbf{D} :=\quad \quad } )
	h + \sigma _3\textbf{A}_4(t)   P_c (\mathcal{K}_{\omega _0
 } )  h\\& - \label{eq:eqh23}    \sum _{a=1}^{3} \im  \textbf{A}_a(t)     P_c (\mathcal{K}_{\omega _0
 } ) \partial _{x_a} h+ \sum _{| \textbf{e}   \cdot(\mu-\nu)|>\omega   _0  }
z^\mu \overline{z}^\nu    \textbf{G}_{\mu \nu}  +\textbf{T}\\&   \nonumber
\textbf{T}:=  \sum
z^\mu \overline{z}^\nu    [ \textbf{G}_{\mu \nu}(t,  \Pi (f)   ) - \textbf{G}_{\mu \nu}]
+R_1 ^{\prime\prime }+ R_2 ^{\prime\prime }    \, , \quad   \textbf{G}_{\mu \nu}:=\textbf{G}_{\mu \nu}( t, 0   )  \end{align}
For  $T> \epsilon _0^{-1}$,  we   eliminate  the second term in  the second line of \eqref{eq:eqh23} setting
\begin{equation}
  \label{eq:g variable}
g=h+ Y \,  , \quad Y:=\sum _{| \textbf{e}  \cdot(\mu-\nu)|>\omega _0} z^\mu
\overline{z}^\nu
   R ^{+}_{\mathcal{K}_{\omega _0
 }  }  (\textbf{e}  \cdot(\mu-\nu) )
     \textbf{G}_{\mu \nu}  .
\end{equation}
\begin{lemma}\label{lemma:bound g} Assume the hypotheses of Prop. \ref{prop:mainbounds}  and let  $T> \epsilon _0^{-1}$. Then   for fixed $s>1$ there  exist
a fixed $c$
such that if   $\epsilon_0$ is sufficiently small, for any preassigned and large $L>1$  we have $\| g
\| _{L^2((0,T ), L^{2,-s}_x)}\le ( c  +\frac {C_0}L) \epsilon  $.
\end{lemma}
We postpone the proof to Sect. \ref{lem:g}.

\subsection{The discrete modes and the  Fermi Golden Rule}
\label{subsec:FGR}

Now we need   show that in $[0,T]$ inequalities \eqref{Strichartzradiation}--\eqref{L^inftydiscrete} hold for $C=C_0/2$ if   they hold for $C=C_0$.  First  we consider the elementary case when $T\le  \epsilon _0^{-1}$.
We have already seen in Remark \ref{rem:tspan1} that  \eqref{L^2discrete}
for $C=C_0/2$ is straightforward  when $T\le  \epsilon _0^{-1}$.
Elementary is also to get  \eqref{L^inftydiscrete} for $C=C_0/2$.
Indeed, we write
\begin{equation}
  \label{eq:tspan2 } \begin{aligned} &   \im \dot z_j =\partial _{\overline{z}_j}(  {H} _2 '+  Z_0+  Z_1+   \resto  +     \langle  V(\cdot + \textbf{v} t) u
,  u  \rangle ) + \im \dot \Pi _b \partial _{\pi _b} z_j.\end{aligned}
\end{equation}
Then   for $0\le t\le T\le \epsilon ^{-1}$ we have
\begin{equation*}
  \label{eq:tspan3 } \begin{aligned} &   |z _j(t)|^2=|z _j(0)|^2+\\& 2\int _0^t \Im   \{ (\partial _{\overline{z}_j} Z_1+\partial _{\overline{z}_j}  \resto  +  \partial _{z_j}  \langle  V(\cdot + \textbf{v} t) u
,  u  \rangle + \im \dot \Pi _b \partial _{\pi _b} z_j)  \overline{z}_j \}   ds \le  |z _j(0)|^2+ \\&   C  (C_0)(   T \epsilon ^3 + \epsilon ^2+ \|  \partial _{z_j}  \langle  V(\cdot + \textbf{v} t) u
,  u  \rangle \| _{L^1 (0,T) } )  \epsilon  \le |z _j(0)|^2+ C  (C_0)   \epsilon ^2,  \end{aligned}
\end{equation*}
 using      the fact that
{Lemma}
  \ref{lem:disppot1}   proves
$
  \|  \partial _{z_j}  \langle  V(\cdot + \textbf{v} t) u
,  u  \rangle \| _{L^1(  0,T    )} =O( \epsilon ^2)$ and Lemmas \ref{lem:dotPi} and \ref{lem:reg5} (3).   For $C_0$ large, this yields  \eqref{L^inftydiscrete}
for $C=C_0/2$.

We   assume now $T>  \epsilon _0^{-1}$. On large intervals  $[0,T]$   the  coupling between discrete and continuous modes is  crucial   and we need
the Fermi Golden Rule.

We have, summing over appropriate indexes $(\mu ,\nu )$ s.t. $|\textbf{e}  \cdot (\mu   -\nu )| > \omega
_0$,
\begin{equation*}   \begin{aligned} &  \partial _{\overline{z}_j}Z_1= \im    \sum   \nu _j {z ^\mu
 \overline{ {z }}^ { {\nu} } }{\overline{z}_j}^{-1}  \langle f ,J
  {G}
_{\mu \nu }(t,\Pi (f) )\rangle      \\& =  \im      \sum   \nu _j {z ^\mu
 \overline{ {z }}^ { {\nu} } }\ {\overline{z}_j}^{-1}  \langle e^{-\frac 12 J v_0 \cdot x}   Mh  ,J
   {G}
_{\mu \nu }(t,\Pi (f) )\rangle     \\&  =        \sum   \nu _j {z ^\mu
 \overline{ {z }}^ { {\nu} } }\ {\overline{z}_j}^{-1}  \langle    h  ,M^T\im J \overline{M}\, \overline{M}^{-1}
  e^{ \frac 12 J v_0 \cdot x}  {G}
_{\mu \nu }(t,\Pi (f) )\rangle  \\&  =     2    \sum   \nu _j {z ^\mu
 \overline{ {z }}^ { {\nu} } }\ {\overline{z}_j}^{-1} \langle    h  ,\sigma _3   \overline{\textbf{G}}
_{\nu \mu }(t,\Pi (f) )\rangle ,
\end{aligned}
\end{equation*}
where we used \eqref{eq:Homega2},\eqref{eq:eqh223} and the symmetry \eqref{eq:symm} which is satisfied by ${\textbf{G}}
_{\mu \nu }$.

\noindent Set $R_{\mu \nu }^+=R_{ \mathcal{K}_{\omega _0} }^+ (\textbf{e}\cdot (\mu   -\nu
) ).$
Using \eqref{eq:partham} and \eqref{eq:g variable} we have
\begin{equation}\label{eq:FGR0} \begin{aligned} & \im \dot z _j
=
\partial _{\overline{z}_j}(H_2'+Z_0) +2   \sum  _{
 |\textbf{e}  \cdot (\mu   -\nu )| > \omega
_0 }  \nu _j\frac{z ^\mu
 \overline{ {z }}^ { {\nu} } }{\overline{z}_j}  \langle g ,\sigma _3   \overline{\textbf{G}}
_{\nu \mu }(t,\Pi (f) )\rangle       +
\partial _{  \overline{z} _j}  \resto  \\&  +  \im    \dot \Pi _b \partial _{\pi _b} z_j
  +
   \frac{\gamma }{2}   \partial _{\overline{z}_j} \langle
			V(\cdot + \textbf{v}t ) u , u  \rangle
 \\ &  -2 \sum  _{ \substack{| \textbf{e}   \cdot
(\alpha    -\beta )|> \omega _0
\\
 |\textbf{e}    \cdot (\mu -\nu )|> \omega
_0 }}  \nu _j\frac{z ^{\mu +\alpha }  \overline{{z }}^ { {\nu}
+\beta}}{\overline{z}_j} \langle   R_{   \alpha \beta}^+ \textbf{G}
_{\alpha \beta } , \sigma _3   \overline{\textbf{G}}
_{\nu \mu }(t,\Pi (f) )\rangle     .
\end{aligned}  \end{equation}
Like in   \cite{Cu2},  there is a  new set of variables
   $\zeta =z+O(z^2)$ s.t.     for a fixed $C$
   \begin{equation}  \label{equation:FGR3} \begin{aligned}   & \| \zeta  -
 z  \| _{L^2_t}
\le CC_0\epsilon ^2\, , \quad  \| \zeta  -
 z \| _{L^\infty _t} \le C \epsilon ^3 \quad \text{ and}
\end{aligned}
\end{equation}
\begin{equation} \label{equation:FGR4} \begin{aligned} &
   \im \dot \zeta
 _j=
\partial _{\overline{\zeta}_j}H_2 (\zeta , h ) +
\partial _{\overline{\zeta}_j}Z_0 (\zeta , h )+  \mathcal{D}_j
   +  \im    \dot \Pi _b \partial _{\pi _b} z_j
  +
   \frac{\gamma }{2}   \partial _{\overline{z}_j} \langle
			V(\cdot + \textbf{v}t ) u , u  \rangle   \\&  -2 \sum  _{ \substack{ \textbf{e}
\cdot  \alpha =\textbf{e}\cdot   \nu   > \omega _0
   \\ \textbf{e}
\cdot  \alpha -\textbf{e} _k   < \omega _0   \, \forall \, k \, \text{
s.t. } \alpha _k\neq 0\\ \textbf{e} \cdot  \nu -\textbf{e} _k  <\omega _{0}  \,
\forall \, k \, \text{ s.t. } \nu _k\neq 0}} \nu _j \frac{\zeta ^{
\alpha } \overline{ \zeta}^ { \nu }}{\overline{\zeta}_j}  \langle   R_{   \alpha 0}^+  \textbf{G}
_{\alpha 0 } , \sigma _3\overline{\textbf{G}}
_{  \nu  0} \rangle  ,
\end{aligned}
\end{equation}
  where    for a
fixed constant $c_0$ and a preassigned large $L>1$ we have for $ {\mathcal{D} }_j'=\mathcal{D}_j$
\begin{equation}\label{eq:FGR7} \sum _{j=1}^{\textbf{n}}\| {\mathcal{D}}_j' \overline{\zeta} _j\|_{
L^1[0,T]}\le (1+C_0)( c _0 +\frac {C_0}L)  \epsilon ^{2}
 . \end{equation}
 By  Lemmas \ref{lem:dotPi} and  \ref{lem:disppot1}  we have also  for fixed $c_0$
\begin{equation*}  \sum _{j=1}^{\textbf{n}}\|( \im    \dot \Pi _b \partial _{\pi _b} z_j
  +
   \frac{\gamma }{2}   \partial _{\overline{z}_j} \langle
			V(\cdot +\textbf{v}t ) u , u  \rangle) \overline{\zeta} _j\|_{
L^1[0,T]}\le  c_0 \epsilon ^{2}
 , \end{equation*}
so that we can absorb the   $\mathcal{D}_j$ and the second line of  \eqref{equation:FGR4}  in a single
remainder term   $\mathcal{D}_j'$ which satisfies  \eqref{eq:FGR7} too.
   We have, see \cite{Cu2},
   \begin{align} \label{eq:FGR5}
 &\partial _t \sum _{j=1}^{\textbf{n}} \textbf{e} _j
 | \zeta _j|^2  =  2  \sum _{j=1}^{\textbf{n}} \textbf{e} _j\Im \left (
\mathcal{D}_j '\overline{\zeta} _j \right ) -\\&    -4 \sum _{
\substack{ \textbf{e} \cdot  \alpha =\textbf{e}  \cdot   \nu   >
\omega _0
   \\ \textbf{e}
\cdot  \alpha -\textbf{e}   _k   < \omega _0   \, \forall \, k \, \text{
s.t. } \alpha _k\neq 0\\ \textbf{e} \cdot  \nu -\textbf{e}  _k  <\omega _0
\, \forall \, k \, \text{ s.t. } \nu _k\neq 0}} \textbf{e}  \cdot \nu
\Im  \left ( \zeta ^{ \alpha } \overline{\zeta }^ { \nu  } \langle
R_{ \alpha 0}^+   \textbf{G}_{ \alpha 0} , \sigma _3   \overline{\textbf{ G} }  _{ \nu 0
}\rangle \right ) .\nonumber
\end{align}

\noindent In the second line of  \eqref{eq:FGR5} we have a sum
\begin{equation}   \label{eq:FGR8} \begin{aligned} &      \Gamma (\zeta ):=     4\sum _{\Lambda  >\omega
_0 } \Lambda      \Im \left    \langle R_{ \mathcal{K}_{\omega _0}}^+ (\Lambda  )
\textbf{G} (\Lambda , \zeta),
 \sigma _3 \overline{\textbf{G} (\Lambda , \zeta)        }\right \rangle , \\&   \text{ for }  \textbf{G} (\Lambda , \zeta)
:=\sum _{
\textbf{e}  \cdot \alpha =\Lambda   }\zeta ^{ \alpha }   \textbf{G}_{ \alpha 0} .
\end{aligned}
\end{equation}
  For
   $W(\omega ) =\lim_{t\to\infty}e^{-\im t \mathcal{K}_\omega }e^{\im t\sigma_3
(-\Delta+\omega  )}$, there exist  $ F   \in W^{k,p}(\R ^3, \C ^2)$
 for all $k\in \R$ and $p\ge 1$ with
 $ 2\textbf{G} (\Lambda , \zeta) =W(\omega _0)F  $, \cite{Cu6,CPV}. Then for
 $^t{ {F}}=( {F}_1, {F}_2)  $
\begin{equation} \label{eq:FGR81} \begin{aligned} &
 \Gamma (\zeta ) = \lim _{\varepsilon \searrow 0} \int _{\R ^3} \frac{\varepsilon}{(\xi ^2- (\Lambda-\omega _0))^2+\varepsilon ^2  }  |\widehat{ {F}}_1 (\xi )|^2 d\xi \ge 0,
\end{aligned}
\end{equation}
  by \cite{Cu2} and  Lemma 4.1  \cite{Cu7}.
Now we   assume:

\begin{itemize}
\item[(H11)]   for some fixed constant $\Gamma >0$  and for all $\zeta \in
\mathbb{C} ^{\mathbf{n}}$ we have \end{itemize}
\begin{equation} \label{eq:FGR} \begin{aligned} &  \Gamma (\zeta )
 \ge \Gamma  \sum _{ \substack{ \mathbf{e}\cdot  \alpha
> \omega _0
\\
   \mathbf{e}
\cdot  \alpha -\mathbf{e} _k   < \omega _0 \\ \forall \, k \, \text{
s.t. } \alpha _k\neq 0}}  | \zeta ^\alpha  | ^2 .
\end{aligned}
\end{equation}
Then integrating    and exploiting \eqref{equation:FGR3} we get for $t\in [0,T]$
\begin{equation}\label{eq:FGRfgr} \sum _j \textbf{e} _j  |z
_j(t)|^2 +4\Gamma \sum _{ \substack{ \text{$\alpha$ as in (H11)}}}  \|z ^\alpha \| _{L^2(0,t)}^2\le
c(1 +  C_0  +   {C_0^2} L^{-1} )\epsilon ^2.
\end{equation}
By \eqref{eq:FGRfgr} and by Lemma \ref{lem:conditional4.2}
we conclude that for $\epsilon _0>0$  sufficiently small and any $T>0$,
\eqref{Strichartzradiation}--\eqref{L^inftydiscrete} in  $I=[0,T]$ and with $C=C_0$  implies
\eqref{Strichartzradiation}--\eqref{L^inftydiscrete} in  $I=[0,T]$  with $C=c(1+\sqrt{C_0} +\frac {C_0}L)$ for $c$.
This yields Proposition  \ref{prop:mainbounds}.

\section{Linear dispersion}
\label{sec:dispersion}

In this section we will set $\mathcal{K}_0= \sigma _3( -\Delta +\omega _0)$,
$\mathcal{K}_1=\mathcal{K} _{\omega _0}=\mathcal{K}_0 +\mathcal{V}_1   $,
 $\mathcal{K}_2= \mathcal{K}_0 +\mathcal{V}_2   $ where $\mathcal{V}_2=\gamma
 \sigma _3 {V}$.
We denote by   $\textbf{D}(t)$ a given $C^1([0,T] , \R^3) $ function such that
 \begin{equation} \label{eq:divcenter10}  \begin{aligned} & \gamma  \| \langle \textbf{D} (t) \rangle ^{-2} \|_{L^1(0,T)} <  c\epsilon   \text{   for a fixed  $c$}\\&     | \textbf{D} (t)| \ge (t-t_0) \frac{|\textbf{v}|}2 -|\textbf{D} (t_0)|   \text{   for    $0\le t_0<t \le T$}  .   \end{aligned}
 \end{equation}
The   $\textbf{D} (t) := \textbf{v} t+ D' +\resto ^{0,2} $ in  \eqref{eq:eqh23}   satisfies  \eqref{eq:divcenter10},   by
the arguments of
Lemmas   \ref{lem:drift}  and   \ref{lem:divcenter}       and assuming \eqref{Strichartzradiation}--\eqref{L^inftydiscrete}
with $C=C_0$ in $I=[0,T]$.
We set
	 \begin{align} &   \text{$ \mathcal{V}_2^D(t,x) = \mathcal{V}_2(x+\textbf{D}(t))$,
   $P_c =P_c(\mathcal{K}_1)$,  $P_d=1-P_c,$ $\mathcal{K }(t)=\mathcal{K}_0 +\mathcal{V}_1 +\mathcal{V}_2^D(t)$,}\nonumber\\&
  \mathcal{U}(t,s) =e^{-\im (t-s)\mathcal{K}_0}W^{-1}(t)W(s) \ ,    \text{ for }   W(t) = e^{-\int _0^t (v(s) \cdot \nabla +
  \im  \varphi (s) \sigma _3) ds }.  \label{eq:semigroup}
\end{align}

 In this section we seek the following  combination of the dispersive
 estimates for the
 \textit{Charge Transfer Models} in \cite{RSS1,cai}    with the integrating factor in \cite{beceanu}.
   \begin{theorem}\label{thm:strich}
   Consider  for $P_cF(t)=F(t)$ and  $P_cu(0)=u_0$
	the equation
   \begin{equation} \label{eq:strich1}\im \dot u -
P_c\mathcal{K}(t)P_cu- \im  P_c v(t) \cdot \nabla _x u  +  \varphi (t) P_c \sigma _3u=F
 \end{equation}
   for $(v(t),\varphi (t))\in C^1 ([0,T]  , \R ^3\times \R  )$.
   For $ \mathbf{{v}}$ the vector   in Theor.\ref{theorem-1.1}, set
   \begin{equation} \label{eq:strich2} \begin{aligned} &  c(T):=\|  (\varphi  (t), v(t),\mathbf{v}- \dot {\textbf{D}} (t)  )  \| _{L^\infty _t[0,T] + L^1 _t[0,T]   }   . \end{aligned}
 \end{equation}
    Then there exist fixed $\sigma >3/2 $,    $c_0>0$  and a  $C>0$  such that, if $c(T)<c_0  $  then    for any admissible pair   $ (p,q)$, see \eqref{admissiblepair}, we have for $i=1$
\begin{equation}\label{eq:strich3} \|  u \|
_{L^p_t( [0,T],W^{i,q}_x)}\le
 C (\|  u_0 \| _{H^i }  +  \|  F \|
_{L^2_t( [0,T],H^{i, \sigma}_x)+ L^1_t( [0,T],H^{i }_x)} ) .
\end{equation}
\end{theorem}

	Before proving Theorem \ref{thm:strich}
  for the $P_d$ in \eqref{eq:semigroup}  we  consider
\begin{equation} \label{eq:smooth0}\begin{aligned}  &\im \dot u -  \mathcal{K}_0
u- \im   v(t) \cdot \nabla _x u  +  \varphi (t) \sigma _3u =\mathcal{V}_2^D u+   {G} u -\im \delta P_d u
 \, , \quad u(t_0) =u_0,\end{aligned}
 \end{equation}
 where      $
 G(t):=   \mathcal{V}_1   - P_d\mathcal{K}(t)P_c-\mathcal{K}(t)P_d.$
For the
$u(t)$ in Theorem   \ref{thm:strich}
 we have $P_du(t)=0$,
 so
	$\im \delta P_d u  $ with fixed $\delta >0$ does not changes the problem.

\begin{proposition}
  \label{prop:weights2}    Let $U(t,t_0)$ be the group associated to   \eqref{eq:smooth0}.
	 Then for $\sigma >3/2$ there exists a fixed
	 $C>0$
  such that for all $0\le t_0 < t \le T$
\begin{equation}  \label{eq:weights2}
\begin{aligned}   &\| \langle x - x_0 \rangle ^{-\sigma} U(t,t_0)  \langle x - x_1 \rangle ^{-\sigma} \| _{2\to 2}\le C \langle t -t_0 \rangle ^{- \frac{3}{2}}  \quad \text{   $\forall \ ( x_0, x_1)\in \R^6$,}
\end{aligned}
\end{equation}
\begin{equation} \label{eq:weights21}
\begin{aligned}
\int_0^T ||\<x-x(t)\>^{-\sigma}U(t,t_0) u_0 ||_{L^2_x}^2  dt \leq C ||u_0||_{L^2_x} ^2   \text{   $\forall \ x(t) \in C^0([0,T], \R^3)$.}
\end{aligned}
\end{equation}
\end{proposition}
We have the following elementary lemma.
 \begin{lemma}
  \label{lem:weights1}   For any fixed $\sigma >3/2$ there is a $C>0$
  such that
\begin{align}  \nonumber &\| \langle x - x_0 \rangle ^{-\sigma} e^{ -{\im \mathcal{K}_0 (t-t_0)}}  \langle x - x_1 \rangle ^{-\sigma} \| _{2\to 2}\le C \langle t - t_0 \rangle ^{- \frac{3}{2}}  \text{ $\forall$ $(t,x_0, x_1)\in \R^7$.}
\end{align}
\end{lemma}
\proof For $|s |:= |t-t_0|\le 1$ follows by $\|  e^{ -{\im \mathcal{K}_0 s}}    \| _{2\to 2}\le 1$. For $|s |\ge 1$ by
\begin{align}  \nonumber &\| \langle x - x_0 \rangle ^{-\sigma} e^{ -{\im \mathcal{K}_0 s }}  \langle x - x_1 \rangle ^{-\sigma} \| _{2\to 2}\le  \|  \langle x   \rangle ^{-\sigma}  \| _{L^2} ^2  \|  e^{  {\im s\Delta  }}   \| _{1\to \infty }
 \le C |s |   ^{- \frac{3}{2}} .
\end{align}\qed

Before proving  Proposition \ref{prop:weights2}
we show that it implies Theorem \ref{thm:strich}.

  \proof [Proof of Theorem \ref{thm:strich}]
 It is   enough to
  prove  \eqref{eq:strich3} for $i=0$.
Write \eqref{eq:strich1} as

\begin{equation*} \label{eq:strich4}\begin{aligned}  \im \dot u -  \mathcal{K}_0
u- \im   v(t) \cdot \nabla _x u  +  \varphi (t) \sigma _3u=F(t)+\mathcal{V}_2^D(t) u+ \widehat{G}(t)   ,
 \end{aligned}
 \end{equation*}
\begin{equation}\label{eq:strich5} \widehat{G}(t)  :={G}(t)+ Y(t) \, , \quad Y(t):=-\im  P_d v(t) \cdot \nabla _x   +\varphi (t) P_d\sigma _3  .
\end{equation}
Then, we have
\begin{equation*}
\|  u \|_{L^p_t( [0,T],L^{ q}_x)}   \lesssim    \|  u_0 \| _{L^2 }  + \|  F \|_{L^2_t( [0,T],L^{ \frac{6}{5}}_x)+ L^1_t( [0,T],L^{2 }_x)}   + \| ( \mathcal{V}_2^D + \widehat{G} ) u \|_{L^2_t(  [0,T],L^{ \frac{6}{5}}_x)   }  .
\end{equation*}
Furthermore,   for   $M>0 $  fixed   large
\begin{equation*}\label{eq:strich6}
\begin{aligned}  \|  (\mathcal{V}_2^D + \widehat{G}) u \|
_{L^2_t L^{ \frac{6}{5}}_x}   \lesssim          \|  \langle  x +\textbf{D}\rangle ^{ -M   }   u\|
_{ L^2_tL^{ 2 }_x  } +   \|  \langle  x \rangle ^{ -M }   u\|
_{ L^2_tL^{ 2 }_x  }.
\end{aligned}
 \end{equation*}
Thus, it suffices to show   for any     $y(t)\in  C^0 ([0,T]  , \R ^3)$ the following inequality
\begin{equation}\label{eq:smooth1} \| \langle x-y(t)\rangle ^{-A}   u  \|
_{L^2_t( [0,T],L^2_x)}\le C \left (
  \|  u  (0) \| _{L^2 }  +  \|  F \|
_{L^2_t( [0,T],L^{ 2, \sigma}_x)+ L^1_t( [0,T],L^{2 }_x)}\right ) ,
\end{equation}
where $A>A_0>0$,   with $A_0>0$    and  $C>0$    some fixed constants.

To show \eqref{eq:smooth1}, we use Duhamel formula and expand $u$.
We  have
\begin{equation} \label{eq:resol} \begin{aligned} &
    u(t) =  \mathcal{U} (t,0) u(0) -\im \int _0^t   \mathcal{U} (t, s) F(s) ds    -\im \int _0^t    \mathcal{U}   (t, s) (\mathcal{V}_2^D + \widehat{G})u(s) ds
.\end{aligned}
\end{equation}
  Then the weighted norms of the first two terms in the rhs of  \eqref{eq:resol} can be bounded
	by the rhs of  \eqref{eq:smooth1}. For example, we have  for fixed $C$
\begin{equation*} \label{eq:smooth3}\begin{aligned}  \| \langle x-y(t)\rangle ^{-A}  \int _0^t   \mathcal{U}   (t, s)F(s) ds  \|
_{L^2_t( [0,T],L^2_x)} & \lesssim   \|    \int _0^t   \mathcal{U}   (t, s)F(s) ds  \|
_{L^2_t( [0,T],L^6_x)} \\& \lesssim      \|  F \|
_{L^2_t( [0,T],L^{ \frac{6}{5}}_x)+ L^1_t( [0,T],L^{2 }_x)}
.\end{aligned}
 \end{equation*}
Therefore, it suffices to bound the following by the rhs of \eqref{eq:smooth1}:
\begin{equation}\label{eq:smooth2} \| \langle x-y(t)\rangle ^{-A}   \int _0^t   \mathcal{U}   (t, s) (\mathcal{V}_2^D + \widehat{G})u(s) ds  \|
_{L^2_t( [0,T],L^2_x)}.
\end{equation}
We expand $u(s)$ by using the group $U(t,t_0)$ associated to   \eqref{eq:smooth0}:
\begin{equation}\label{eq:duh00}
\begin{aligned}
u(s)=U(s,0)u_0 -	\im \int_0^s U(s,\tau) F(\tau)\,d\tau -\im \int_0^s U(s,\tau)Y(\tau)u(\tau)\,d\tau.
\end{aligned}
\end{equation}

\noindent  Substituting  $u(s)$ in  \eqref{eq:smooth2}   with  the rhs of  \eqref{eq:duh00} we
reduce to bound three terms. We first estimate the contribution from  the first term in rhs of \eqref{eq:duh00}:
\begin{align}
&\nonumber  \| \langle x-y(t)\rangle ^{-A}   \int _0^t    \mathcal{U}   (t, s)(\mathcal{V}_2^D + \widehat{G})U(s,0)u_0 ds  \|
_{L^2_t( [0,T],L^2_x)} \lesssim \\& \nonumber
\left \| \int_0^t \<t-s\>^{-\frac 3 2}     (\|\<x+\textbf{D} \>^A \mathcal{V}_2^D U(s,0)u_0 \|_{L^2_x}+\|\<x\>^A  \widehat{G} U(s,0)u_0 \|_{L^2_x} )ds\right  \|_{L^2_t( 0,T)}\\& \nonumber
\lesssim  ||\<x +\textbf{D}\>^{-A}U(t,0)u_0||_{L^2_t( [0,T],L^2_x)}   +||\<x\>^{-A}U(t,0)u_0||_{L^2_t( [0,T],L^2_x)} \lesssim ||u_0||_{L^2_x},  \nonumber
\end{align}
where we have used Lemma \ref{lem:weights1} in the second line, the inequalities
\begin{equation*}\label{eq:GY}
||\<x+\textbf{D}\>^A\mathcal{V}_2^D u||_{L^2_x}\lesssim ||\<x+\textbf{D}\>^{-A}u||_{L^2_x}  \text {  and } ||\<x\>^A\widehat{G}u||_{L^2_x}\lesssim ||\<x\>^{-A}u||_{L^2_x}
\end{equation*}
in the third line
 and \eqref{eq:weights21} in  Proposition \ref{prop:weights2} in the final line.

We next estimate the contribution to \eqref{eq:smooth2} of the second term
in the rhs of \eqref{eq:duh00}.
First, we bound it by $||F||_{L^1_t([0,T] ,L^2_x)}$.    We have
\begin{equation*}\label{eq:duh002}
\begin{aligned}
 & \| \langle x-y(t)\rangle ^{-A}   \int _0^t   \mathcal{U}   (t, s) \mathcal{V}_2^D(s)   \int_0^s U(s,\tau) F(\tau)\,d\tau ds  \|
_{L^2_t( [0,T],L^2_x)}\\&
\lesssim   \left \| \int_0^t \<t-s\>^{-3/2} \int_0^s\|\<x+\textbf{D}(s)\>^{-A}U(s,\tau)F(\tau)\|_{L^2_x}\,d\tau ds\right \|_{L^2_t( [0,T])}\\&
\lesssim   \int_0^T d\tau
\left  \| \int_\tau ^t \<t-s\>^{-3/2} \|\<x+\textbf{D} (s)\>^{-A}U(s,\tau)F(\tau)\|_{L^2_x} ds\right \|_{L^2_t( [0,T])}
\\&
\lesssim   \int_0^T   ||\<x+\textbf{D}(s)\>^{-A}U(t,\tau)F(\tau)||_{L^2_t ([0,T], L^{2 }_x)}\,d\tau  d\tau
\leq   C ||F||_{L^1_t ([0,T], L^{2 }_x)}.
\end{aligned}
\end{equation*}
The term with $\mathcal{V}_2^D$ replaced by $ \widehat{G}$ is obtained similarly with $\textbf{D}$ replaced by 0.

\noindent We  bound the same terms    by $||F||_{L^2_t([0,T],L^{2,\sigma}_x)}$ using the Young's inequality:
\begin{equation}\label{eq:duh003}
\begin{aligned}
&  \| \langle x-y(t)\rangle ^{-A}   \int _0^t    \mathcal{U}   (t, s)  (\mathcal{V}_2^D+\widehat{G})  \int_0^s U(s,\tau) F(\tau)\,d\tau ds  \|
_{L^2_t( [0,T],L^2_x)}\\&
\leq C \left \| \int_0^t \<t-s\>^{-3/2} \int_0^s\<s-\tau\>^{-3/2}\|F(\tau)\|_{L^{2,\sigma}_x}\,d\tau ds\right \|_{L^2_t( [0,T])}\\&
\leq C \left \| \int_0^t\<t-\tau\>^{-3/2}\|F(\tau)\|_{L^{2,\sigma}_x}\,d\tau  \right  \|_{L^2_t( [0,T])}
\leq C ||F||_{L^2_t([0,T],L^{2,\sigma}_x)}.
\end{aligned}
\end{equation}
We   bootstrap the final term  of   \eqref{eq:duh00}  by using the smallness of $c(T)$:
\begin{equation}\label{eq:duh004}
\begin{aligned}
&\| \langle x-y(t)\rangle ^{-A}   \int _0^t    \mathcal{U}   (t, s)(\mathcal{V}_2^D+\widehat{G}) \int_0^s U(s,\tau)Y(\tau)u(\tau)\,d\tau ds  \|
_{L^2_t( [0,T],L^2_x)}\\&
\leq c(T) \left \| \int_0^t \<t-s\>^{-3/2}\int_0^s\<s-\tau\>^{-3/2}||\<x\>^{-\sigma}u(\tau)||_{L^2_x}\,d\tau ds\right \| _{L^2_t([0,T])}\\&
\leq c(T)||\<x\>^{-\sigma }u||_{L^2_{t,x}([0,T])} \text{ where we used:}
\end{aligned}
\end{equation}
 in the 1st inequality   $\|\<x\>^AY(\tau)u(\tau)||_{L^2_x}\leq C c(T)||\<x\>^{-\sigma }u(\tau)||_{L^2_x}$,
by   definition of $Y(t)$ in  \eqref{eq:strich5};  in the 2nd inequality    Young's inequality twice as  in \eqref{eq:duh003}.
\qed

\subsection{Proof of {Proposition} \ref{prop:weights2}}
\label{subsec:smooth2}

The proof     follows the argument in   \cite{RSS1,cai}.  We will need the following lemma.

\begin{lemma}
  \label{lem:error1}   Let $M>5/2$ and $\alpha\in [0,1/2)$.
Then there exists  $C>0$ such that
   \begin{align}  &   \| \langle x  \rangle ^{-M}\left ( W^{-1}(t)W(s) -1\right ) e^{ -{\im \mathcal{K}_0 (t-s)}}  \langle x   \rangle ^{-M } \| _{2\to 2} \label{eq:error11}\\&  \le C \psi _\alpha (t  -s)   \| (v, \varphi )\| ^\alpha _{L^1 ([s,t]) +L^\infty ([s,t])}  \text{ with $W(t)= e^{-\int _0^t (v(s) \cdot \nabla +
  \im  \varphi (s) \sigma _3) ds }$}\nonumber \\& \text{and $\psi _\alpha (t   ) =
\langle t\rangle ^{- \frac{3}{2} +\alpha   }$ for $t\ge 1$ and  $\psi _\alpha (t   ) =
  t  ^{-  \alpha   }$  for $t\in (0,1)$.}\label{eq:psialpha}\end{align}

  \end{lemma}
\proof   We follow  Lemma 2.16  \cite{beceanu}. First of all, we split the operator in \eqref{eq:error11} and reduce to the bound
\begin{equation}\label{eq:error12}
\begin{aligned}  &  \|  \langle x  \rangle ^{-M}\left ( e^{-\im \sigma _3 \int _s^t
     \varphi (s')   ds'} -1\right )  e^{-\int _s^t  v(s') \cdot \nabla  ds '}  e^{ -{\im \mathcal{K}_0 (t-s)}}  \langle x   \rangle ^{-M } \| _{2\to 2} +  \\&  \|   \langle x  \rangle ^{-M}  e^{-\im \sigma _3 \int _s^t
     \varphi (s')   ds'} \left (  e^{-\int _s^t  v(s') \cdot \nabla  ds '} -1\right )     e^{ -{\im \mathcal{K}_0 (t-s)}}  \langle x   \rangle ^{-M }\| _{2\to 2} . \end{aligned}
\end{equation}
For $M>3/2 $ the first line is bounded  for a fixed $C$ by
\begin{equation*}
\begin{aligned}  &  C\langle t  -s\rangle ^{-  \frac{3}{2}     } \left | e^{-\im \sigma _3 \int _s^t
     \varphi (s')   ds'} -1\right |   \le  C 2^{1-\alpha} \langle t  -s\rangle ^{-  \frac{3}{2} + \alpha   } \|   \varphi  \| ^\alpha _{L^1 ([s,t]) +L^\infty ([s,t])}. \end{aligned}
\end{equation*}
We turn now to the second line of \eqref{eq:error12}. We can drop the factor  $e^{-\im \sigma _3 \int _s^t
     \varphi (s')   ds'}$  and replace $\mathcal{K}_0$ with $-\Delta$. Set also $\delta v=-\int _s^t  v(s')    ds '.$ Then
\begin{equation*}
\begin{aligned}  &   e^{\delta v \cdot \nabla   } e^{  {\im \Delta (t-s)}}  f(x)=
(t-s) ^{-\frac{3}{2}} \int _{\R^3}  e^{\im \frac{(x-y+\delta v)^2}{t-s}} f(y) dy=
e^{  {\im \Delta (t-s)}}  f(x) + \\&  (t-s) ^{-\frac{3}{2}} \int _{\R^3}  e^{\im \frac{(x-y)^2}{t-s}}   \left [ e^{2\im \frac{(x-y)\cdot\delta v }{t-s}} \left ( e^{\im \frac{ (\delta v ) ^2}{t-s}} -1 \right )  + \left (e^{2\im \frac{(x-y)\cdot\delta v }{t-s}} -1\right ) \right ] f(y) dy \end{aligned}
\end{equation*}
where we ignored irrelevant constants.
It remains to bound the terms in the last line. We focus only on the last. The other can be treated similarly. Set
\begin{equation*}
\begin{aligned}  &   e^{2\im \frac{(x-y)\cdot\delta v }{t-s}} -1   = e^{2\im \frac{x\cdot\delta v }{t-s}}  \left ( e^{-2\im \frac{ y\cdot\delta v }{t-s}} -1 \right ) +(e^{2\im \frac{x\cdot\delta v }{t-s}}-1) .\end{aligned}
\end{equation*}
Looking again  only at the last difference, we have for $M> 5/2$
\begin{equation*}
\begin{aligned}  &    \|   \langle x  \rangle ^{-M} (e^{2\im \frac{x\cdot\delta v }{t-s}}-1)  e^{  {\im \Delta (t-s)}} (\langle \cdot    \rangle ^{-M } f)\| _{L^2} \\& \le C (t-s) ^{ -\alpha }
 |\delta v|^{\alpha}   \|  \langle x  \rangle ^{-M+\alpha } e^{  {\im \Delta (t-s)}} \langle x  \rangle ^{-M } \| _{2\to 2} \| f\| _{L^2}\\& \le C  \psi _\alpha (t-s) \|  v  \| ^\alpha _{L^1 ([s,t]) +L^\infty ([s,t])}\| f\| _{L^2}. \end{aligned}
\end{equation*}
Other terms can be treated similarly.
\qed

\proof [Proof of {Proposition} \ref{prop:weights2}]
First of all, our first claim is  that  for any $T<\infty$ there is a $C=C_T$ such that \eqref{eq:weights2}--\eqref{eq:weights21}
are true.
We  rewrite \eqref{eq:smooth0} as (notice $u=P_cu$)
\begin{equation} \label{eq:smooth5}
  \im \dot u -  \mathcal{K}_0
u- \im   v(t) \cdot \nabla _x u  +  \varphi (t) \sigma _3u = \overbrace{(\mathcal{V} _1 -\im \delta P_d  - \mathcal{K}_1P_d)}^{ V_1V_2:=\quad \quad  \quad } u  +   P_c\mathcal{V} _2^D P_c u ,
 \end{equation}
where $V_1$ is a rapidly decreasing potential and where $V_2\in B (H ^{k,-s}, H ^{k, s})$ for any $(k,s)$.
The  factorization
  $V_1V_2$ is    as in \cite{beceanu} and is required by the additional presence,
   with respect to \cite{RSS1},
   of the terms $ \im   v  \cdot \nabla _x u  $ and $\varphi  \sigma _3u$  in \eqref{eq:smooth5}. We have

\begin{equation} \label{eq:duh1}
\begin{aligned} &U(t,t_0)        =
\mathcal{U}   (t, t_0)
-\im \int _{t_0}^t    \mathcal{U}   (t, s)V_1V_2  U(s,t_0)  ds
\\&   -\im \int _{t_0}^t   \mathcal{U}   (t, s)P_c\mathcal{V} _2^D(s)P_c  U(s,t_0)  ds
 =:I_1+I_2+I_3
.\end{aligned}
\end{equation}
By     Lemma \ref{lem:weights1}  and by the endpoint  Strichartz estimate,   there exists a fixed  $C$  independent
of $x_0,x_1$ and $x(t)$ such that
\begin{equation} \label{eq:duhfirst}
\begin{aligned}
&\|\<x-x_0\>^{-\sigma} \mathcal{U}   (t, t_0)\<x-x_1\>^{-\sigma}\|_{2\to 2}\leq C\<t-t_0\>^{-3/2}   , \\&
\int_{ 0}^T\|\<x-x(t)\>^{-\sigma}\mathcal{U}   (t, t_0) u_0\|_{L^2_x}^2 \,dt\leq C \| u_0 \|_{L^2}^2.
\end{aligned}
\end{equation}
Our first claim  follows by \eqref{eq:duh1}--\eqref{eq:duhfirst} and Gronwall's inequality.
Furthermore,    $C(T)=C_{T-t_0} $, that is it depends on  $T-t_0.$
We will show   that $C(T)$  can be taken independent of $T$.
  We   follow the bootstrap   argument in \cite{RSS1}  based on the observation that it is enough to  show that if   \eqref{eq:weights2}--\eqref{eq:weights21} hold  for  $C=C(T)$
for some large $ C(T)$,   they hold  also  for  $C=C(T)/2$. In the argument we can assume $T-t_0>10 A$  for $A$   some fixed  but arbitrarily  large number.
  From    \eqref{eq:duhfirst},  $I_1$  satisfies the desired bounds.
To bound the contributions of $I_2$ and $I_3$, we split $[t_0, t]=[t_0,t_0+A]\cup [t_0+A,t-A] \cup [t-A,t]$ and
set
 \begin{align}
&\nonumber I_{2,t_0}:=  \int _{t_0}^{t_0+A}    \mathcal{U}   (t, s)V_1V_2  U(s,t_0)  ds    ,  \quad  I_{2,A}:=  \int _{t_0+A}^{t-A}    \mathcal{U}   (t, s)V_1V_2  U(s,t_0)  ds    , \\&  \nonumber
I_{2,0}:=  \int _{t-A}^t    \mathcal{U}   (t, s)V_1V_2  U(s,t_0)  ds   ,  \
I_{3,t_0}:=  \int _{t_0}^{t_0+A}   \mathcal{U}   (t, s)P_c\mathcal{V} _2^D(s)P_c  U(s,t_0)  ds   , \\&  \nonumber
I_{3,A}:=  \int _{t_0+A}^{t-A}   \mathcal{U}   (t, s)P_c\mathcal{V} _2^D(s)P_c  U(s,t_0)  ds    ,  \\&
I_{3,0}:=  \int _{t-A}^t   \mathcal{U}   (t, s)P_c\mathcal{V} _2^D(s)P_c  U(s,t_0)  ds \  . \nonumber
\end{align}
We bound $I_{j,t_0}$ and $I_{j,A}$ $(j=2,3)$   as follows.
First, for $I_{3,t_0}$,
\begin{equation*}\label{eq: Ijminor1}
\begin{aligned}
&||\<x-x_0\>^{-\sigma}I_{3,t_0} \<x-x_1\>^{-\sigma}||_{2\to2}\\ &\lesssim \int_{t_0}^{t_0+A}\<t-s\>^{-3/2}||\<x+\textbf{D}(s)\>^{-\sigma}U(s,t_0)\<x-x_1\>^{-\sigma}||_{2\to2}\,ds\\& \le  C(A) \int_{t_0}^{t_0+A}\<t-s\>^{-3/2}\<s-t_0\>^{3/2}
\le  C_A\<t-t_0\>^{-3/2}.
\end{aligned}
\end{equation*}
  $C_A$ does not depend on $T$, so   this  bound is    of the desired type. We have also
\begin{equation*}\label{eq: Ijminor2}
\begin{aligned}
&   \int_0^T||\<x-x(t)\>^{-\sigma}I_{3,t_0} u_0||_{L^2_x}^2\,dt \\&\lesssim       \int_0^T
  \left (     \int_{t_0}^{t_0+A} ds \<t-s\>^{-\frac 32}||\<x+\textbf{D}(s)\>^{-\sigma}U(s,t_0)u_0||_{L^2_x} \,dt\right ) ^{2} \\&
\lesssim  ||\<x+\textbf{D}(s)\>^{-\sigma}U(s,t_0)u_0||_{L^2_t( (t_0,t_0+A),L^2_x)}^2 \le  C_A  ||u_0||_{L^2_x}^2.
\end{aligned}
\end{equation*}
$I_{2,t_0}$ can   bounded similarly  replacing $\textbf{D}$  with 0.
Next, for $I_{j,A}$  $j=2,3$     we have
\begin{equation*}\label{eq: Ijminor3}
\begin{aligned}
||\<x-x_0\>^{-\sigma}I_{j,A}\<x-x_1\>^{-\sigma}||_{2\to 2}&\lesssim
C(T)\int_{t_0+A}^{t-A}\<t-s\>^{-\frac 3 2}\<s-t_0\>^{-\frac 3 2} ds\\&
\le C A^{-\frac 1 2}C(T)\<t-t_0\>^{-\frac 3 2} \text{ for a fixed $C$}.
\end{aligned}
\end{equation*}
Here, since $A\gg 1$, this is a bound of the desired type. Similarly
\begin{equation*}\label{eq: Ijminor4}
\begin{aligned}
&\int_0^T||\<x-x_0\>^{-\sigma}I_{3,A}u_0||_{L^2_x}^2dt\\&\lesssim \int_0^T  \left ( \int_{t_0+A}^{t-A}\<t-s\>^{-\frac 32 }||\<x+\textbf{D}(s)\>^{-\sigma}U(s,t_0)u_0||_{L^2_x}  ds\right ) ^2 dt\\&
\le k(A) || \<x+\textbf{D}(t)\>^{-\sigma}U(t,t_0)u_0||_{L^2 ( ( 0,T),L^2_x)}^2  \le  k (A) C^2(T) \|u_0\|_{L^2_x}^2
 \end{aligned}
\end{equation*}
   by Young's inequality, with
\begin{align} \nonumber & k (A):=\sup _{0\le s\le T-A } \left  (   \int _{s+A}^T \<t-s\>^{-\frac 32 } dt     \right )
\sup _{A\le t\le T  } \left  (   \int^{t-A }_0 \<t-s\>^{-\frac 32 } ds   \right )\\& \lesssim C A^{-1}. \label{eq: Ijminor41}
\end{align}
Then   $k (A) C^2(T)$ is an arbitrarily small fraction of  $C^2(T)$. The same bound  holds for $I_{2,A}$
replacing $\textbf{D}$ with 0.

\noindent To estimate   $I_{2,0}$  we expand $V_2U(s,t_0)$    following \cite{beceanu}.
By  \eqref{eq:duh1} we get
\begin{equation} \label{eq:duh2}
\begin{aligned} (1+\im  {T}_0)   V_2U(\cdot,t_0)   (s)     & =\im
({T}_0 -  \widetilde{T}_0) V_2U(\cdot,t_0)   (s)
 +  V_2
\mathcal{U}   (s, t_0)
\\&   -\im  V_2\int _{t_0}^s    \mathcal{U}   (s, \tau ) P_c \mathcal{V} _2^D P_c  U(\tau ,t_0) d\tau \, ,
 \end{aligned}
\end{equation}
where  $ \widetilde{T}_0$, $  {T}_1$ and  $  {T}_1$ are for $\mathcal{M}_0= \mathcal{K}_0 $ and $\mathcal{M}_1= \mathcal{K}_1P_c-\im \delta P_d$ defined as follows:
 \begin{align}\nonumber
 & \widetilde{T}_0f (s)    := V_2 \int _{t_0}^s    \mathcal{U}   (s, \tau  )V_1f(\tau   ) d\tau \, ,  \quad  {T}_0f (s)      := V_2 \int _{t_0}^s   e^{-\im (s-\tau )\mathcal{M}_0 }  V_1f(\tau   ) d\tau   \, ,   \\&  {T}_1f (s)      := V_2 \int _{t_0}^s   e^{-  \im (s-\tau  )\mathcal{M}_1 }  V_1f(\tau   ) d\tau   \,
 .  \label{eq:duh3}
 \end{align}
 Then by    \cite{beceanu} or    by    formula (B.28)
\cite{NS}  the  equality  $\mathcal{M}_1 -\mathcal{M}_0  =V_1V_2 $ implies
\begin{equation}\label{lem:inverse}
    (1-\im  {T}_1)(1+\im  {T}_0)=1.
\end{equation}
Indeed we have the following, which implies \eqref{lem:inverse}:
 \begin{equation*}
\begin{aligned}   & T_1T_0f(t) =V_2\int  _{t_0}^t ds e^{-\im (t-s )\mathcal{M}_1 } V_1V_2\int  _{t_0}^s  e^{-\im ( s-\tau  )\mathcal{M}_0 }V_1 f(\tau ) d\tau \\&   =V_2\int  _{t_0}^t  d\tau  \int  _{\tau }^t ds   e^{-\im (t-s )\mathcal{M}_1 } (\mathcal{M}_1 -\mathcal{M}_0)   e^{-\im ( s-\tau  )\mathcal{M}_0 } V_1f (\tau )    \\&    = -\im V_2\int  _{t_0}^t  d\tau     (  e^{-\im (t-\tau )\mathcal{M}_0 } -   e^{-\im (  t-\tau  )\mathcal{M}_1 } )  V_1f (\tau ) =  \im T_1 f(t) -\im T_0f(t)
 \end{aligned}
\end{equation*}
 By \eqref{eq:duh2} and   \eqref{lem:inverse}, we have
\begin{equation} \label{eq:duh4}
\begin{aligned}    V_2U(s,t_0)      & =(1-\im  {T}_1) V_2
 \mathcal{U}   (\cdot, t_0)    (s) + (1-\im  {T}_1)\im
({T}_0 - \widetilde{T}_0) V_2U(\cdot,t_0)   (s)
\\&   -\im (1-\im  {T}_1)  V_2\left [\int _{t_0}^\cdot    \mathcal{U}   (\cdot, \tau ) P_c\mathcal{V} _2^D(\tau )P_c
 U(\tau,t_0) d\tau   \right ] (s).
 \end{aligned}
\end{equation}

We set
$ I_{2,0}=I_{2,1}+I_{2,2}+I_{2,3}$, where
\begin{equation*}
\begin{aligned} & I_{2,1} :=-\im \int _{t-A}^t    \mathcal{U}   (t, s)V_1(1-\im  {T}_1) V_2
\mathcal{U}   (\cdot, t_0)    (s)  ds\\&
I_{2,2} :=  \int _{t-A}^t   \mathcal{U}   (t, s)V_1(1-\im  {T}_1)
({T}_0 - \widetilde{T}_0) V_2U(\cdot,t_0)   (s)  ds\\&
I_{2,3} := - \int _{t-A}^t   \mathcal{U}   (t, s)V_1 (1-\im  {T}_1)  V_2\left [\int _{t_0}^\cdot    \mathcal{U}   (\cdot, \tau ) P_c\mathcal{V} _2^D(\tau )P_c
 U(\tau,t_0) d\tau   \right ] (s)  ds
 \end{aligned}
\end{equation*}

We estimate $I_{2,1}$.
We split $I_{2,1}=-\im \varsigma -\varpi$, where
\begin{equation*}\label{eq: I211}
\begin{aligned}
&\varsigma:=   \int _{t-A}^t    \mathcal{U}   (t, s)V_1 V_2
\mathcal{U}   (s, t_0)      ds\, , \
\omega:=  \int _{t-A}^t    \mathcal{U}   (t, s)V_1  {T}_1 V_2
\mathcal{U}   (\cdot, t_0)    (s)  ds.
\end{aligned}
\end{equation*}

We estimate $\varsigma$. We have:
\begin{equation}
\begin{aligned}
\|\<x-x_0\>^{-\sigma}\varsigma\<x-x_1\>   ^{-\sigma} \|_{2\to2}
&\lesssim \int_{t-A}^t\<t-s\>^{-\frac 3 2}\<s-t_0\>^{- \frac 3 2} ds
\lesssim \<t-t_0\>^{-\frac 3 2} ;
\end{aligned}\nonumber
\end{equation}
\begin{equation}
\begin{aligned}&
 \|\<x-x(t)\>^{-\sigma}\varsigma u_0 \|_{  L^2_t ((0,T), L^2_{x})}
\\&=     \|\<x- x (t)\>^{-\sigma} \int _{t-A}^t    \mathcal{U}   (t, s)V_1 V_2
\mathcal{U}   (s, t_0) u_0  ds\|_{  L^2  ((0,T), L^2_{x})} \lesssim \\&    \|  \int _\R    \<t-s\>^{-\frac3 2}  \| V_2
\mathcal{U}   (s, t_0) u_0\|  _{    L^2_{x} }   ds\|_{  L^2  ( 0,T) }
 \lesssim  \|   V_2
\mathcal{U}   (\cdot , t_0) u_0 \|_{  L^2  (0,T)  L^2_{x} } \lesssim
  ||u_0||_{L^2_x} .
\end{aligned}\nonumber
\end{equation}
We   estimate $\varpi$.  We will use an   analogue of Lemma \ref{lem:weights1}    for $\mathcal{K}_1P_c$ and  $\mathcal{K}_2$.
\begin{lemma}
  \label{lem:weights12}
  $\forall$ fixed $\sigma >3/2$ $\exists$ $C>0$
  such that for $\mathbf{K}=\mathcal{K}_1P_c,\mathcal{K}_2$ we have
  \begin{align} \label{eq:weights22}  &\| \langle x - x_0 \rangle ^{-\sigma} e^{ -{\im \mathbf{K} (t-t_0)}}  \langle x - x_1 \rangle ^{-\sigma} \| _{2\to 2}\le C \langle t - t_0 \rangle ^{- \frac{3}{2}},  \text{$\forall$ $(t,x_0, x_1)\in \R^7$}\\&
   \int _\R \| \langle x -x(t) \rangle ^{-\sigma} e^{ -{\im \mathbf{K} t}}    \| _{2\to 2}^2  dt\le C , \text{ $\forall$ $x(t)\in C^0([0,T], \R ^3)$.}\label{eq:weights33}
  \end{align}
\end{lemma}
 \proof  \eqref{eq:weights22} follows   by $ \|   e^{ -{\im \mathbf{K} t}}P_c \| _{2\to 2} \le C$ and $\|   e^{ -{\im \mathbf{K} t}}    \| _{1\to \infty}\le C | t| ^{- \frac{3}{2}}$.  For $\mathbf{K}=\mathcal{K}_2$  see \cite{JSS,Y1}.
 For  $\mathbf{K}=\mathcal{K}_1P_c$ see \cite{Cu6,CPV} or \cite{schlag}.

 \eqref{eq:weights33} is  obtained by endpoint  Strichartz estimate like the second line in \eqref{eq:duhfirst}
and is
due to   \cite{Y1} for  $\mathbf{H}=\mathcal{H}_2$
 and to   \cite{Cu6,CPV} for  $\mathbf{H}=\mathcal{H}_1P_c$.
\qed

\noindent Using \eqref{eq:weights22} we get for fixed constants
\begin{equation}
\begin{aligned}
&\|\<x-x_0\>^{-\sigma}\varpi\<x-x_1\> ^{-\sigma} \|_{2\to2} \lesssim \\&
\int_{t-A}^t\<t-s\>^{-\frac 32}\int_{t_0}^s\|\<x\>^{-\sigma } e^{-(s-\tau)(\im \mathcal{K}_1P_c+\delta P_d)}V_1V_2\mathcal{U}(\tau,t_0)\,d\tau \<x-x_1\>^{-\sigma}||_{2\to 2}\\&
\lesssim   \int_{t-A}^t\<t-s\>^{-\frac 32 }\int_{t_0}^s \| \<x\>^{-\sigma} e^{-(s-\tau)(\im\mathcal{H}_1P_c+\delta P_d)}\<x\>^{-\sigma}||_{2\to 2}\<\tau-t_0\>^{-\frac 32}d\tau ds\\&
\lesssim \int_{t-A}^t\<t-s\>^{-\frac 32}\int_{t_0}^s \<s-\tau\>^{-\frac 32}\<\tau-t_0\>^{-3/2}\,d\tau ds
\lesssim \<t-t_0\>^{-\frac 32} ,
\end{aligned}\nonumber
\end{equation}
\begin{equation}
\begin{aligned}
&\|\<x-x_1\>^{-\sigma}\varpi u_0\|_{  L^2  ((0,T), L^2_{x})}  \lesssim \\&
\|   \int_{t-A}^t\<t-s\>^{-\frac 3 2}    \int_{t_0}^s  \<s -\tau \>^{-\frac 3 2} \|   V_2\mathcal{U}(\tau,t_0) u_0
||_{   L^2_{x} }
\|_{  L^2  ( 0,T) }
\\& \lesssim   \| V_2\mathcal{U}(\tau,t_0) u_0  \|_{  L^2  ((0,T), L^2_{x})}
  \le C
  ||u_0||_{L^2_x} .
\end{aligned}\nonumber
\end{equation}
Since
$T_1$ doesn't affect     estimates   we will  estimate only terms  without $T_1$.

We estimate $I_{2,2}$. For a fixed $C$,   $c(T) +A^{-1/2}\ll 1 $ and the $\psi _{\alpha}$ in \eqref{eq:psialpha}
  \begin{equation}\label{eq:I22}
\begin{aligned}
&\|\<x-x_0\>^{-\sigma}\int _{t-A}^t   \mathcal{U}   (t, s)V_1
({T}_0 - \widetilde{T}_0) V_2U(\cdot,t_0)   (s)  ds\<x-x_1\>^{-\sigma}\|_{2\to2}\\&
\lesssim \int_{t-A}^t\<t-s\>^{-3/2}\int_{t_0}^s \| \<x\>^{-M}(W^{-1}(s)W(\tau)-1)e^{-i(s-\tau)\mathcal{H}_0}\<x\>^{-M}
\\&\qquad\qquad\times\<x\>^{M}V_2U(\tau,t_0)\<x-x_1\>^{-\sigma} \|_{2\to 2} \,d\tau ds \\&
\lesssim c(T)^{\alpha}C(T)\int_{t-A}^t\<t-s\>^{-3/2}\int_{s-A}^s
\psi _{\alpha} (s-\tau)  \<\tau-t_0\>^{-3/2}\,d\tau ds\\&
+C(T)\int_{t-A}^t\<t-s\>^{-3/2}\int_{t_0+A}^{s-A}\<s-\tau\>^{-3/2 }\<\tau-t_0\>^{-3/2}\,d\tau ds\\&
+C_A \int_{t-A}^t\<t-s\>^{-3/2}\int_{t_0}^{t_0+A}\<s-\tau\>^{-3/2}\<\tau-t_0\>^{-3/2}\,d\tau ds\\&
\le C \left (c(T)^{\alpha}C(T)+A^{-1/2}C(T)+C_A\right )\<t-t_0\>^{-3/2}.
\end{aligned}
\end{equation}
$I_{2,2}$   has a term with $T_1$ with   similar bound. With  $c(T)^{  \alpha}C(T)\ll C(T)$,
\begin{equation*}
\begin{aligned}
&\|   \<x-x(t)\>^{-\sigma}\int _{t-A}^t   \mathcal{U}   (t, s)V_1
({T}_0 - \widetilde{T}_0) V_2U(\cdot,t_0)   (s) u_0 ds\| _{  L^2  ((0,T), L^2_{x})}\\&
\lesssim   \|  \int_{t-A}^t\<t-s\>^{-\frac 32 } \int_{t_0} ^s  \| \<x\>^{-M}(W^{-1}(s)W(\tau)-1)e^{-i(s-\tau)\mathcal{K}_0}\<x\>^{-M}
\\&\qquad\qquad\times  \<x\>^{M}V_2U(\tau,t_0)u_0 \| _{    L^2_{x} }    ds \, d\tau    \|  _{  L^2  ( 0,T) }
\\&
\lesssim c(T)^{ \alpha} \|    \int_{t-A}^t\<t-s\>^{-\frac 32 }\int_{t_0}^s\psi _{\alpha} (s-\tau)\|\<x\>^{-\sigma}U(\tau,t_0)u_0\|_{L^2_x}^2\,d\tau ds  \|  _{  L^2  ( 0,T) }  \\&
\lesssim c(T)^{  \alpha}C(T) \| u_0\|_{L^2_x}^2.
\end{aligned}
\end{equation*}
The part of $I_{2,2}$ with $T_1$ is bounded similarly.
So   the estimate of $I_{2,2}$ is proved.

We now estimate $I_{2,3}$ restricting integration of $\tau$ in $I_{2,3}$ to   $[s-A,s]$, since
the estimates for   $[t_0,t_0+A]$ and $[t_0+A,s-A]$ are like those of $I_{2,t_0}$ and $I_{2,A}$.
We split $P_c\mathcal{V}_2^D P_c=\mathcal{V}_2^D-P_d\mathcal{V}_2^DP_c-\mathcal{V}_2^D P_d$.
We first bound the contributions of $ P_d\mathcal{V}_2^DP_c$ in $I_{2,3}$.
Using $P_d=\sum_j (\cdot,\tilde{e}_j)e_j$, where $e_j$, $\tilde{e}_j\in \mathcal{S}(\R^3)$,
\begin{align} \nonumber
&||\<x-x_0\>^{-\sigma}\int_{t-A}^t\mathcal{U}(s,t)V_1 V_2\int_{s-A}^s\mathcal{U}(s,\tau)P_d\mathcal{V}_2^D(\tau )P_cU(\tau,t_0) \<x-x_1\>^{-\sigma}u_0||_{L^2_x}\\& \label{eq:I23r1}
\lesssim \int_{t-A}^t\<t-s\>^{-\frac 32}\int_{s-A}^s\<s-\tau\>^{-\frac 32}\\&  \nonumber \qquad\times \sum_j |(\<x\>^{-\sigma}P_cU(\tau,t_0)\<x-x_1\>^{-\sigma}u_0,\<x\>^{\sigma}\mathcal{V}_2(\cdot +\textbf{D}(\tau ))\tilde{e}_j)|\,d\tau ds\\&
\lesssim   C(T) l(A)  \<t-t_0\>^{-\frac 32}\|u_0\|_{L^2_x}  \, , \, l(A):= \sum_j|| \<x\>^{\sigma} \mathcal{V}_2  (\cdot +\textbf{D})\tilde{e}_j||_{L^\infty_t ((A,T),L^2)}
. \nonumber
\end{align}
By \eqref{eq:divcenter10},  $l(A) $    is arbitrarily close to $0$ for   $A $  sufficiently large.
We have
\begin{align}
&   \label{eq:I23r10} \| \<x-x(t)\>^{-\sigma}\int_{t-A}^t\mathcal{U}(s,t)V_1\\& \nonumber\qquad\times V_2\int_{s-A}^s\mathcal{U}(s,\tau)P_d\mathcal{V}_2^D(\tau )P_cU(\tau,t_0)u_0 d\tau   ds    \| _{  L^2   (0,T)  L^2 _{x} }\\& \nonumber
\lesssim   \|  \int_{t-A}^t\<t-s\>^{-3/2}\int_{s-A}^s\<s-\tau\>^{-3/2}\times\\&  \sum_j|(\<x\>^{-\sigma}P_cU(\tau,t_0) u_0,\<x\>^{\sigma}\mathcal{V}_2(\cdot +\textbf{D}(\tau ) )\tilde{e}_j)|  \,d\tau ds     \| _{  L^2   (0,T) }
\lesssim C_T c(A) \|u_0\|_{L^2_x}, \nonumber
\end{align}
 $c(A)$ as in  \eqref{eq:I23r1}.
  $\mathcal{V}_2^DP_d$ can be treated similarly.
The remaining part of   $I_{2,3}$ is
\begin{equation}   \label{eq:duh43}
   \tilde{I}_{2,3}:=
  \int _{t-A}^t   \mathcal{U}   (t, s)V_1       V_2\int _{s-A}^s     \mathcal{U}   (s, \tau ) \mathcal{V} _2^D
 U(\tau,t_0) d\tau ds    .
\end{equation}
We have, for  $y(s, \tau) =\textbf{D}(\tau ) +\int _\tau ^s v(r) dr$,
\begin{equation}   \label{eq:duh44}
  V_2    \mathcal{U}   (s, \tau ) \mathcal{V} _2^D
  =    V_2    e^{-\im (s -\tau )\mathcal{K}_0}   e^{-\im \sigma _3 \int _\tau ^s
   \varphi (r)  dr } \mathcal{V} _2  (\cdot  + y(s, \tau) ) .
\end{equation}
 By \eqref{eq:divcenter10} and    since $c(T)$ in  \eqref{eq:strich2} is small,  we can assume
\begin{equation}   \label{eq:duh441}
\begin{aligned}       &   \int _{s-A}^s \langle y(s, \tau) \rangle ^{-\frac 32}  d\tau  =  \int _{s-A}^s \langle \textbf{D}(\tau ) +\int _\tau ^s v(r) dr \rangle ^{-\frac 32}  d\tau
\\& \le 2    \int _{s-A}^s \langle \textbf{D}(\tau )   \rangle ^{-\frac 32}  d\tau \le 2c \epsilon
\end{aligned}
\end{equation}
by $\langle \textbf{D}(\tau ) +\int _\tau ^s v(r) dr \rangle \sim \langle \textbf{D}(\tau )   \rangle$
which follows by $|\int _\tau ^s v(r) dr |\le c(T)A<1$.

\noindent We ignore  $e^{-\im \sigma _3 \int _\tau ^s
   \varphi (r)  dr }$ without harm.
As in \cite{RSS1}, for a smooth partition
\begin{equation}   \label{eq:duh45}
\begin{aligned}       &   1= F(|\im \nabla |\le N ) + F(|\im \nabla |\ge N )
\end{aligned}
\end{equation}
in low and high frequencies it     sufficies to estimate the following $J_{L}$ and $J_{H}$:
\begin{equation*}\label{eq:highlow}
\begin{aligned}
&J_{L}=\int _{t-A}^t   \mathcal{U}   (t, s)V_1       \int _{s-A}^s     V_2    e^{-\im (s -\tau )\mathcal{K}_0}  F(|\im \nabla |\le N ) \mathcal{V} _2  (\cdot + y(s, \tau) )
 U(\tau,t_0) d\tau ds \\&
 J_{H}=\int _{t-A}^t   \mathcal{U}   (t, s)V_1       \int _{s-A}^s     V_2    e^{-\im (s -\tau )\mathcal{K}_0}  F(|\im \nabla |\geq N ) \mathcal{V} _2  (\cdot + y(s, \tau) )
 U(\tau,t_0) d\tau ds
\end{aligned}
\end{equation*}

For $J_L$ and for $\mathcal{V} _2^\sigma  (x):= \mathcal{V} _2 (x)\langle  x \rangle ^{ \sigma }    $  we use the identity
\begin{align}    \label{eq:duh46}   &  \|   V_2    e^{-\im (s -\tau )\mathcal{K}_0}      F(|\im \nabla |\le N ) \mathcal{V} _2^\sigma   (\cdot +y(s, \tau)   )    u    \| _{L^2 }    =     \|   \int _{\R ^3}  k(x,\eta )   \widehat{u} (\eta )    d\eta   \| _{L^2 }
\\&  k(x,\eta ):=
V_2(x)   \int _{\R ^3} e^{-\im (s -\tau )\sigma _ 3 (\xi ^2 +\omega _0)
+\im  \xi \cdot  (x- y(s, \tau)) }      \chi (\frac \xi N ) \widehat{\mathcal{V}} _2^\sigma   (\xi -\eta    )    d\xi  e^{\im  \eta \cdot y(s, \tau)} .\nonumber
\end{align}
     In the   next computation we drop the irrelevant
the  $\omega _0$   in the phase.
We claim
\begin{equation}   \label{eq:duh461}
\begin{aligned}       &   |k(x,\eta )|  \le C_{N,M} \<x\>^{-M} \left \langle  y(s, \tau) \rangle ^{-M}    \langle  \eta \right  \rangle ^{-M}.
\end{aligned}
\end{equation}
Then  \eqref{eq:duh46}--\eqref{eq:duh461}   imply the following bound:
\begin{equation}   \label{eq:duh462}
\begin{aligned}       &  \|   V_2    e^{-\im (s -\tau )\mathcal{K}_0}      F(|\im \nabla |\le N )
 \mathcal{V} _2^\sigma   (\cdot +y(s, \tau)   )       \| _{2\to 2 }   \le C_{N,M}
 \left  \langle  y(s, \tau) \right \rangle ^{-M} .
\end{aligned}
\end{equation}
Therefore, we can bound $J_L$ as follows:
\begin{equation}
\begin{aligned} &
||\<x-x_0\>^{-\sigma}J_L \<x-x_1\>^{-\sigma} ||_{2\to 2}      \\&  \leq C(T) C_{N } \gamma \int _{t-A} ^t ds\int _{s-A} ^s  d\tau \<t-s\>^{-\frac 32}
\left \langle y(s, \tau) \right \rangle  ^{-\frac 32}   \<\tau -t_0\>^{-\frac 32} \lesssim
  C(T) C_{N } \gamma \\& \times  \<t -t_0\>^{-\frac 32}\int _{t-A} ^t ds\int _{s-A} ^s  d\tau \<t-s\>^{-\frac 32}
\left \langle y(s, \tau) \right \rangle  ^{-\frac 32}
\lesssim       C(T) C_{N } \epsilon
\<t-t_0\>^{-\frac 32} ;
\end{aligned}\nonumber
\end{equation}
\begin{equation}
\begin{aligned} &
 \| \<x-x(t)\>^{-\sigma} J_L u_0 \|_{  L^2  ((0,T), L^2_{x})}    \lesssim C_{N }  \gamma
  \|  \int _{t-A} ^t ds\int _{s-A} ^s d\tau \<t-s\>^{-\frac 32}  \\&    \times
\left \langle  \textbf{D}(\tau )  \right \rangle  ^{-\frac 32}
 \|   \langle  (x + y(s, \tau) )  \rangle  ^{-\sigma } U(\tau ,t_0) u_0 \|_{    L^2_{x} }     \|_{  L^2   (0,T) }
			\\&  \le  C_{N }    \sqrt{A}\gamma   \|  \left \langle  \textbf{D}(\tau )  \right \rangle  ^{-\frac 32}   \|_{  L^2(0,T) }
			 \|  U(t ,t_0) u_0   \|_{  L^\infty   ((0,T), L^2_{x})}\le   C(T) C_{N } \sqrt{A  \epsilon }  \|  u_0   \|_{  L^2_{x} }.
\end{aligned}\nonumber
\end{equation}
This yields the desired bound  picking $\epsilon _0>0$ sufficiently small.
  We get  \eqref{eq:duh461}     setting
$e ^{-\im  I_2 \xi \cdot     y(s, \tau)  }  =  (- \frac {y(s, \tau) \cdot \nabla _\xi }{|y(s, \tau)|^2} )^N e ^{-\im  I_2 \xi \cdot     y(s, \tau)  }  $ and   integrating  by parts
 in  \eqref{eq:duh46}.  Notice that  $\widehat{\mathcal{V}} _2^\sigma   (\xi -\eta    )  $ and its derivatives
are   $\le C(M ,N)   \langle  \eta  \rangle ^{-M}$ for $|\xi |\lesssim N$.

\noindent We now consider $J_H$.
Proceeding as in \cite{RSS1} we have
\begin{equation}   \label{eq:duh47}
\begin{aligned}       & \int _{s-A}^s \|   V_2    e^{-\im (s -\tau )\mathcal{K}_0}      F(|\im \nabla |\ge N ) \mathcal{V} _2^\sigma   (\cdot +y(s, \tau)   )      \| _{2 \to 2}   d\tau  \\& \lesssim   \int _{s-A}^s \|   [V_2,F(|\im \nabla |\ge N )]    e^{-\im (s -\tau )\mathcal{K}_0}       \mathcal{V} _2^\sigma   (\cdot +y(s, \tau)   )      \| _{2 \to 2}   d\tau \\& +  A^{ \frac{1}{2}}N ^{-\frac{1}{2}}
(   \int _{s-A}^s \|   {\langle \im \nabla \rangle } ^{-\frac{1}{2}} \nabla V_2   e^{-\im (s -\tau )\mathcal{K}_0}       \mathcal{V} _2^\sigma   (\cdot +y(s, \tau)   )      \| _{2 \to 2} ^2 d\tau )^{\frac 12}.
\end{aligned}\nonumber
\end{equation}
By
 Young's inequality $ \| [V_2,F(|\im \nabla |\ge N )]  \|  _{2 \to 2}\lesssim N^{-1}$,    \cite{cai}.  In the third line
\begin{equation*}
      {\langle \im \nabla \rangle } ^{-\frac{1}{2}} \nabla V_2 =  V_2 {\langle \im \nabla \rangle } ^{-\frac{1}{2}} \nabla+ [  {\langle \im \nabla \rangle } ^{-\frac{1}{2}} \nabla , V_2].
\end{equation*}
We have $\|   [  {\langle \im \nabla \rangle } ^{-\frac{1}{2}} \nabla , V_2] \| _{2\to 2} \lesssim 1$.
For $\sigma >1/2$ we have
\begin{equation}\label{eq:locsm}
    \int _{\R }  \|   \< x \> ^{-\sigma}   \langle \im \nabla \rangle   ^{-\frac{1}{2}}  \nabla  e^{-\im (s -\tau )\mathcal{K}_0}            \| _{2 \to 2} ^2 d\tau \le C _{\sigma},
\end{equation}
  by    the local smoothing effect, see Theorem 4.3 \cite{lp}.
 Therefore, we have
\begin{equation}
\begin{aligned}
&\|\<x-x_0\>^{-\sigma}J_H\<x-x_1\>^{-\sigma}\|_{2\to 2}\lesssim C(T)\<t-t_0\>^{-3/2}\int_{t-A}^t\<t-s\>^{-3/2}\\&\times\int _{s-A}^s \|   V_2    e^{-\im (s -\tau )\mathcal{K}_0}      F(|\im \nabla |\ge N ) \mathcal{V} _2^\sigma   (\cdot +y(s, \tau)   )      \| _{2 \to 2}   d\tau ds\\&
\lesssim C(T)\<t-t_0\>^{-3/2}\int_{t-A}^t\<t-s\>^{-1/2}\(N^{-1}A+{N} ^{-\frac 12} \sqrt{ {A}}\(1+A^{1/2}\)\)\\&
\lesssim C(T)AN^{-1/2}\<t-t_0\>^{-3/2} \le  C(T) L ^{-1} \<t-t_0\>^{-3/2}  \text{ for $N\gg 1$}
\end{aligned}\nonumber
\end{equation}
for any preassigned   $L\gg 1$.
  Finally, the estimate for $I_{23}$ is completed with
\begin{equation}
\begin{aligned}
& ||\<x-x(t)\>^{-\sigma}J_H u_0 \|_{  L^2  ((2A,T), L^2_{x})}  \lesssim \\&
\|  \int_{t-A }^t\<t-s\>^{-\frac 32}\int_{s-A  }^s \|   V_2    e^{-\im (s -\tau )\mathcal{K}_0}      F(|\im \nabla |\ge N ) \mathcal{V} _2^\sigma   (\cdot +y(s, \tau)   )      \| _{2 \to 2} \\&\times ||\<\cdot-y(s,\tau)\>^{-\sigma}U(\tau,t_0)u_0||_{L^2_x}   d\tau ds  \|_{  L^2  ( 0,T)  }  \lesssim
\|   \int_{t-A }^t\<t-s\>^{-\frac 32 }    \\&   \times \int_{s-A }^s \( {N} ^{-1} + {N} ^{-\frac 12}(\sqrt{ {A}} + {  A } ) \) ||\<\cdot+y(s,\tau)\>^{-\sigma}U(\tau,t_0)u_0||_{L^2_x} \,d\tau ds \|_{  L^2  ( A,T)  }\\&
\lesssim C(T) AN^{-\frac 12 }||u_0||_{L^2_x} .
\end{aligned}\nonumber
\end{equation}

\subsection{Estimate of $I_{3,0}$}
\label{subsec:duh51}

 \begin{lemma}
  \label{lem:commutator}  For  $\textbf{D}(t) \in C^1$   and  $\mathbf{v} \in \R ^3$
   the following operator
  \begin{equation} \label{eq:comm1}
\begin{aligned}   \widetilde{g}(t)u(t,x):=e^{\im \sigma _3(-  \frac{t}{4}\mathbf{v}^2  -\frac{\mathbf{v}\cdot x}{2} ) }
  u(t, x+\textbf{D}(t)) \text{ satisfies:}
\end{aligned}
\end{equation}   \begin{equation}\label{eq:comm2}\begin{aligned} &[\im \partial _t -\mathcal{K}_0,  \widetilde{g}(t)  ]
            =\im  \widetilde{g}(t)
(\dot {\textbf{D}}-\mathbf{v}  ) \cdot \nabla _x  \ , \ [\widetilde{g}(t)^{-1}, \partial_{x_j}]u  = - \im \sigma_3 \frac{\mathbf{v}_{ j}}{2}\widetilde{g}^{-1}(t)u \ , \\&
\widetilde{g}(t)^{-1}u=e^{ \im \sigma _3(  \frac{t}{4}\mathbf{v}^2   +\frac{\mathbf{v}\cdot (x-\textbf{D}(t))}{2} ) }
  u(t, x-\textbf{D}(t)) \ ,\\&
\ [\widetilde{g}(t)^{-1}, \im \partial _t -\mathcal{K}_0]u  =\im   (\dot {\textbf{D}}  -\mathbf{v} ) \cdot \nabla _x (\widetilde{g}^{-1}(t)u) \ .
\end{aligned}
\end{equation}
\end{lemma}
\proof  Follows by   elementary computations.\qed

Consider the  $\widetilde{g}(t)$ in \eqref{eq:comm1}.
Set  $(\hat v , \hat \varphi )=(v-(\mathbf{v}-\dot {\textbf{D}} ), \varphi - v \cdot \mathbf{v}/2)$  and     $g (t) =\widetilde{g}(t) e^{\im \sigma _3 \int _0^t  \hat \varphi  (s) ds} $. Then
 $\mathcal{V} _2^D(t)= g ^{-1}(t)\mathcal{V} _2 g (t)$.   \eqref{eq:smooth5}  can be rewritten as

\begin{equation*} \label{eq:reduh0}\begin{aligned} &  (\im \partial _t -  \mathcal{K}_0)
g^{-1} u- \im   \hat v (t) \cdot \nabla _x  g^{-1}u    = \mathcal{V} _2g^{-1}u
  + g^{-1} Q ,\\&   Q:=  \mathcal{V} _1 -\im \delta P_d  - \mathcal{K}_1P_d  +    P_d \mathcal{V} _2^D P_c -\mathcal{V} _2^DP_d  , \  \hat v := v-(\mathbf{v} -\dot {\textbf{D}}) .\end{aligned}
 \end{equation*}
      $ \hat v $ satisfies   estimates like  $v$.
  For simplicity
	we  drop the hat.  Let $ \mathcal{U}   (t, s )$ be   \eqref{eq:semigroup}
	with this $v$ and    $\varphi =0$.
We have  an expansion like  \eqref{eq:duh1}:
\begin{equation} \label{eq:reduh1}
\begin{aligned}  U(s,t_0)       = &  g(s) g^{-1} (t_0)
 \mathcal{U}   (s, t_0)
-\im  g(s)\int _{t_0}^s    \mathcal{U}   (s, \tau )W_1W_2  g^{-1} (\tau )   U(\tau ,t_0) d\tau
\\&   - \im  g(s)\int _{t_0}^s   \mathcal{U}   (s, \tau )  g^{-1} (\tau )   Q U(\tau ,t_0) d\tau
,\end{aligned}
\end{equation}
where   $\mathcal{V} _2=W_1W_2$ with both $W_j$  rapidly decreasing.
 As in \eqref{eq:duh2}--\eqref{eq:duh3}
\begin{equation} \label{eq:reduh2}
\begin{aligned} &(1+\im T_0)   W_2g^{-1} (\cdot )U(\cdot,t_0)   (s)      =\im
(T_0 - \widetilde{T}_0 ) W_2g^{-1} (\cdot )U(\cdot,t_0)   (s)
\\& +  W_2  g^{-1} (t_0)
\mathcal{U} ^{-1} (s) \mathcal{U}   (t_0)e^{-\im (s-t_0)\mathcal{K}_0}
\\&   -\im  W_2\int _{t_0}^s   e^{-\im (s-\tau )\mathcal{H}_0}
 \mathcal{U} ^{-1} (s) \mathcal{U}   (\tau )g^{-1} (\tau )Q  U(\tau ,t_0) d\tau  \  ,
 \end{aligned}
\end{equation}
where  $ \widetilde{T}_0$ and $  {T}_0$  are defined as in \eqref{eq:duh3}
but with $V_j$ replaced by $W_j$ and
\begin{equation}\label{eq:T1}
\begin{aligned}
 T_1f (s)      := W_2 \int _{t_0}^s   e^{-  \im (s-\tau  )\mathcal{K}_2}  W_1f(\tau   ) d\tau   .
\end{aligned}
\end{equation}
$(1 - \im T_1)(1+ \im T_0)=1$ by $\mathcal{K}_2   - \mathcal{K}_0=W_1W_2     $ and  \eqref{lem:inverse}.
  Then
\begin{align}\nonumber &   W_2g^{-1} (\tau )U(\tau ,t_0)       =(1-\im  {T}_1)\im
(\widetilde{T}_0 -{T}_0) W_2g^{-1} (\cdot )U(\cdot,t_0)   (\tau )
\\&  \label{eq:reduh4} + (1-\im  {T}_1) W_2 g^{-1} (t_0)
\mathcal{U} ^{-1} (\cdot ) \mathcal{U}   (t_0)e^{-\im (\cdot -t_0)\mathcal{K}_0}    (\tau )
\\&   -\im (\mathbf{K}-\im  {T}_1\mathbf{K})    (\tau ) \, , \,
\mathbf{K} (t ):=W_2
  \int _{t_0}^t    e^{-\im (t -r )\mathcal{K}_0}
 \mathcal{U} ^{-1} (\cdot ) \mathcal{U}   (r )g^{-1} (r )Q  U(r,t_0)
   dr    .\nonumber
 \end{align}
By the discussion for   \eqref{eq:I23r1}--\eqref{eq:I23r10}, we can replace $P_c\mathcal{V}_2^DP_c$ with $\mathcal{V}_2^D$ in $I_{3,0}$.
We denote by $\tilde{I}_{3,0}$ the resulting operator:
\begin{equation} \label{eq:reduh110}
\begin{aligned}       &   \tilde{I}_{3,0}= -\im \int _{t-A}^t    \mathcal{U}   (t, s) g  (s )\mathcal{V} _2 g^{-1} (s )  U(s,t_0) ds .
\end{aligned}
\end{equation}
We substitute    $   U(s,t_0)$   with the rhs of   \eqref{eq:reduh1}.
Set $\tilde{I}_{3,0}=I_{3,1}+I_{3,2}+I_{3,3}$, with
\begin{equation*}
\begin{aligned}
& I_{3,1} := -\im \int _{t-A}^t    \mathcal{U}   (t, s) g  (s )\mathcal{V} _2 g^{-1} (t_0)
 \mathcal{U}   (s, t_0)  ds\\&
 I_{3,2}:= -\int _{t-A}^t    \mathcal{U}   (t, s) g  (s )\mathcal{V} _2 \int _{t_0}^s    \mathcal{U}   (s, \tau )W_1W_2  g^{-1} (\tau )   U(\tau ,t_0) d\tau  ds\\&
 I_{3,3}:= - \int _{t-A}^t    \mathcal{U}   (t, s) g  (s )\mathcal{V} _2 \int _{t_0}^s   \mathcal{U}   (s, \tau )  g^{-1} (\tau )   Q U(\tau ,t_0) d\tau  ds.
\end{aligned}
\end{equation*}
The contribution of $I_{3,1}$ can be bounded by a constant independent of $T$.
   We focus now on $I_{3,3}$.  As before, the integral by $\tau$ can be restricted to $[s-A,s]$.
   Therefore, it suffices to estimate
\begin{equation} \label{eq:reduh111}
\begin{aligned}       &    \int _{t-A}^t    \mathcal{U}   (t, s) g  (s )\mathcal{V} _2
 \int _{s-A}^s    \mathcal{U}   (s, \tau )  g^{-1} (\tau )   Q U(\tau ,t_0) d\tau ds.
\end{aligned}
\end{equation}
For simplicity we will focus only on the contribution of  $\mathcal{V}_1$ to $Q$ and look at
\begin{equation} \label{eq:reduh112}
\begin{aligned}       &   \int _{t-A}^t    \mathcal{U}   (t, s) g  (s )\mathcal{V} _2
 \int _{s-A}^s    \mathcal{U}   (s, \tau )  g^{-1} (\tau )   \mathcal{V}_1 U(\tau ,t_0) d\tau ds.
\end{aligned}
\end{equation}
Then, for some $\varphi (s ,\tau )$ and for  $y(s, \tau) =-\textbf{D}(\tau ) +\int _\tau ^s v(r) dr$ we have
\begin{equation*}   \label{eq:reduh44}
\begin{aligned}       &  \mathcal{V} _2
      \mathcal{U}   (s, \tau )  g^{-1} (\tau )   \mathcal{V}_1
  =    \mathcal{V} _2      e^{-\im (s -\tau )\mathcal{K}_0}   e^{ \im \sigma _3 \varphi (s ,\tau )  }
	e^{\frac \im  2\sigma _3 \mathbf{v} \cdot x     }    \mathcal{V} _1  (x + y(s, \tau) ).
\end{aligned}
\end{equation*}
    Then for the operator \eqref{eq:reduh111}
  we have the desired bound,    which is independent of     $\varphi (s ,\tau )$,  by the argument used for  $J_L$ and $J_H$
starting from \eqref{eq:duh45}.

We finally  consider $   I_{3,2}$.
Once again it is enough to consider
\begin{equation} \label{eq:reduh211}
\begin{aligned}       \int _{t-A}^t    \mathcal{U}   (t, s) g  (s )\mathcal{V} _2 \int _{s-A}^s    \mathcal{U}   (s, \tau )W_1W_2  g^{-1} (\tau )   U(\tau ,t_0) d\tau  ds.
\end{aligned}
\end{equation}
We   substitute  $g^{-1} (\tau )   U(\tau ,t_0 )$ with the rhs of  \eqref{eq:reduh4}. Then
the  contribution of the first two lines of  \eqref{eq:reduh4}  is bounded by   $(C+L^{-1}{C(T)} ) $ for   fixed
$L\gg 1$  (using $c(T)$ small) and $C$. The last line  of  \eqref{eq:reduh4}   is bounded by $L^{-1}{C(T)} $ because $\mathbf{K}(t)$
can be treated by the argument
starting from \eqref{eq:duh45} and  used for $J_L$ and $J_H$ .

\section{Proof of Lemma \ref{lemma:bound g}}
\label{lem:g}
 Using the notation of Sect. \ref{sec:dispersion},  with in particular $\mathcal{K}_0 :=\sigma _3(-\Delta +\omega _0) $,  $\mathcal{K}_1 :=\mathcal{K} _{\omega _0} $ and $Y$ as in
 \eqref{eq:g variable}, we have for a preassigned $\delta >0$
\begin{equation}  \label{eq:eqg1}  \begin{aligned} &        \im     \dot g  =\mathcal{K}_{0}   g +     (\mathcal{V}  _1-\im \delta P_d )
	g +     \sigma _3\varphi  g   - \im    v  \cdot \nabla _x g
  +P_c   \mathcal{V}  _2^{D}
	h   +\sum _{b=1}^3\textbf{T}_b ,
 \\&  \textbf{T}_1:=
   -   \sigma _3\varphi  Y   + \im    v  \cdot \nabla _x Y
 \, , \,   \textbf{T}_2:=  \textbf{T} - \varphi  \sigma _3  P_d h - \im  v\cdot  \nabla _x  P_d h  ,\\&   \textbf{T} _3:=\sum _j \left [\partial _{z_j}Y (\im \dot z_j-\mathbf{e}_jz_j)+\partial _{\overline{z}_j}Y (\im \dot {\overline{z}}_j+\mathbf{e}_j\overline{z}_j)
 \right ]
  . \end{aligned}
\end{equation}
In particular we have $\varphi =\mathbf{A}_4$ and $v_a=\mathbf{A}_a$
for $a\le 3$ which  are such that
\begin{equation} \label{eq:eqg-1} \begin{aligned} &
     \|  (\varphi  , v) \| _{L^\infty( (0,T), \R ^4)+ L^1( (0,T), \R ^4) }\le C(C_0) \epsilon
		  .\end{aligned}  \end{equation}

 Set now  $ \mathcal{V}  _1-\im \delta P_d=V_1V_2$ with $V_1(x):= \langle x \rangle ^{-\sigma }$.
 Then we have
    \begin{align} \label{eq:eqg2}
 &   \|  (1+\im \widetilde{T}_0 )V_2g \| _{  L^2 ((0,T), L^{2   }_x)  }\le  \|   \mathcal{U}(t,0)(h(0)+Y
(0))\| _{ L^2 ((0,T), L^{2,-\sigma }_x)} +\\&  \| \int _0^t\mathcal{U}(t,s)P_c   \mathcal{V}  _2^{D}
	h  ds\| _{L^2 ((0,T), L^{2,-\sigma }_x)}
    +\sum _{b=1}^3
 \| \int _0^t\mathcal{U}(t,s)
\textbf{T}_b(s) ds\| _{   L^2 ((0,T), L^{2,-\sigma }_x)}\nonumber
    \end{align}
for $ \widetilde{T}_0  $ as in \eqref{eq:duh3}, but with the $ V_1 $
and
$ V_2$ introduced here   above \eqref{eq:eqg2}. We have
$\widetilde{T}_0= {T}_0 +(\widetilde{T}_0- {T}_0)$ with
  $\|\widetilde{T}_0- {T}_0\| _{ L^2 ((0,T), L^{2  }_x)\to L^2 ((0,T), L^{2  }_x)}\le C \epsilon ^{\alpha}$
by Lemma
\ref{lem:error1}    for some $\alpha \in (0,1).$ Furthermore,   \eqref{lem:inverse} continues to hold. Hence
 \begin{equation} \label{eq:eqg3}  \begin{aligned}
 &  \|   g \| _{  L^2 ((0,T), L^{2,-\sigma    }_x)}\le C  \|  (1+\im \widetilde{T}_0 )V_2g \| _{  L^2 ((0,T ) , L^{2   }_x)  }
    \end{aligned}
 \end{equation}
where $C=\|  1-\im  {T}_1\| _{L^2(\R^4,\C^2 ) \to L^2(\R^4, \C^2)}<\infty $.

By   \eqref{eq:eqg2}--\eqref{eq:eqg3} the proof of Lemma \ref{lemma:bound g} will follow by showing that the rhs of  \eqref{eq:eqg2} is bounded by $c \epsilon +C(C_0)\epsilon ^2 $.
We have
\begin{equation} \label{eq:eqg4}  \begin{aligned}
 &  \|    \mathcal{U}(t,0) h(0) \| _{  L^2 ((0,T), L^{2,-\sigma    } )}\lesssim  \|    e^{-\im t \mathcal{K}_0} h(0) \| _{  L^2 (\R , L^{6    } )} \lesssim  \|    h(0) \| _{    L^{2    }  }  \le     c \epsilon .
    \end{aligned}
 \end{equation}
We have
\begin{equation} \label{eq:eqg5}  \begin{aligned}
 &  R^{+}_{\mathcal{K}_1}(\lambda ) =R^{+}_{\mathcal{K}_0}(\lambda )-R^{+}_{\mathcal{K}_0}(\lambda ) \mathcal{V}_1R^{+}_{\mathcal{K}_1}(\lambda )  .
    \end{aligned}
 \end{equation}
 By the definition of $Y(t)$ in \eqref{eq:g variable}, for $R^{+}_{\mu \nu} =R^{+}_{\mathcal{K}_1}(\mathbf{e}\cdot (\mu -\nu )) $  as in Sect. \ref{subsec:FGR},
\begin{equation} \label{eq:eqg6}  \begin{aligned}
 &  \|    \mathcal{U}(t,0) Y(0) \| _{  L^2 ((0,T), L^{2,-\sigma    } )}\le \sum | z  (0)| ^2 \| \mathcal{U}(t,0)R^{+}_{\mu \nu}  \mathbf{G}_{\mu \nu} (0)   \| _{  L^2 ((0,T) , L^{2,-\sigma    }) } .
    \end{aligned}
 \end{equation}
We have $  |z(0)|^2\lesssim \epsilon ^2$ and by
\eqref{eq:eqg5} we have
  \begin{align}\nonumber
 &     \| \mathcal{U}(t,0) R^{+}_{\mu \nu}  \mathbf{G}_{\mu \nu} (0)   \| _{  L^2 (\R , L^{2,-\sigma    }) } \le    \| \mathcal{U}(t,0)  R^{+}_{\mathcal{K}_0}(\mathbf{e}\cdot (\mu -\nu )) \mathbf{G}_{\mu \nu} (0)   \| _{  L^2 (\R , L^{2,-\sigma    }) }
 \\&   + \| \mathcal{U}(t,0)  R^{+}_{\mathcal{K}_0}(\mathbf{e}\cdot (\mu -\nu )) \mathcal{V}_1R^{+}_{\mu \nu} \mathbf{G}_{\mu \nu} (0)   \| _{  L^2 (\R , L^{2,-\sigma    }) }  \le C \| \mathbf{G}_{\mu \nu} (0)   \| _{   L^{2, \sigma    }  }   \label{eq:eqg7}
    \end{align}
by   Lemma \ref{lem:lemg9} below and  $ \| R^{+}_{\mu \nu}P_c(\mathcal{K}_1) \| _{ L^{2, \sigma    } \to L^{2,-\sigma    }}< C,$    \cite{CPV}.
By \eqref{eq:eqg6}--\eqref{eq:eqg7}
\begin{equation} \label{eq:eqg8}  \begin{aligned}
 &  \|    \mathcal{U}(t,0) Y(0) \| _{  L^2 ((0,T), L^{2,-\sigma    } )}\le c \epsilon \text{ for a fixed $c$.}
    \end{aligned}
 \end{equation}
The next lemma  is a standard tool in the analysis of the coupling of discrete and continuous modes,   \cite{BP2,SW3}.  Here the novelty is the presence of  $e ^{  \int _s^tdt' v(t')\cdot \nabla }$.

\begin{lemma}\label{lem:lemg9}
Let $\overline{I}$ be a
compact subset of $(\omega _0,\infty)\cup (-\infty ,-\omega _0)$
and let $\sigma >9/2$. Assume
\eqref{eq:eqg-1},  more precisely assume     \begin{equation}\label{eq:lemg90}   \begin{aligned} &
     \|  v  \| _{L^\infty( (0,T), \R ^4)+ L^1( (0,T), \R ^4) }\le C(C_0) \epsilon
		  .\end{aligned}  \end{equation}
			Then there exists a $C=C(\overline{I})>0$
such that  for the  $\epsilon _0$
in Proposition \ref{prop:mainbounds} sufficiently small
we have for every $t\ge s \ge 0$, $\lambda \in \overline{I}$  and $u\in
\mathcal{S}(\R^3;\C^2)$
\begin{equation} \label{eq:lemg91}
\|  e ^{  -\int _s^tdt' v(t')\cdot \nabla }e^{-\im \mathcal{K}_0   (t-s)}R_{ \mathcal{K}_0 }^{+}( \lambda )
u \|_{L^{ 2, - \sigma}(\R^3)} \le C\langle  t -s\rangle ^{-\frac 32} \| u \|_{L^{ 2,   \sigma}(\R^3)} .
\end{equation}
\end{lemma}
\proof    We follow Sect.8 \cite{Cu7} which is like Proposition 2.2 \cite{SW3}.
  We will see that
inequality \eqref{eq:lemg90}   plays an important role.  It is not restrictive to assume $\overline{I}\subset J\subset \overline{J}  \subset (\omega _0,\infty)$ with $J$ open and $\overline{J}$ compact.  Let
$\chi  \in C^\infty (\R , [0,1])$ be a smooth cutoff equal to 1 in $\overline{I}$ and to 0 outside $\overline{J}$.  Let $\chi  _1=1-\chi $.
Then    we have
\begin{equation*} \label{eq:lemg92}\begin{aligned} &
\|     \mathcal{U} (t,s) R_{ \mathcal{K}_0 }^{+}( \lambda )\chi  _1 (\mathcal{K}_0)
u \|_{L^{ 2, - \sigma} }  \le  \\&   \|     \mathcal{U} (t,s)  \|_{L^{ 2,  \sigma} \to  L^{ 2, - \sigma} }
 \|    R_{ \mathcal{K}_0 }^+ ( \lambda )\chi  _1 (\mathcal{K}_0)  \|_{L^{ 2,  \sigma} \to  L^{ 2,  \sigma} }
\| u \|_{L^{ 2, - \sigma} }
\le C \langle   t-s  \rangle ^{-\frac 32} \| u \|_{L^{ 2,   \sigma} } \end{aligned}
\end{equation*}
 by  Lemma \ref{lem:weights1}.
 So   \eqref{eq:lemg91}  holds with  an additional
   $\chi  _1 (\mathcal{K}_0)$.
We need to prove
 \begin{equation} \label{eq:lemg94}\begin{aligned}&
\|    A  \langle  \cdot  \rangle ^{ \sigma } u \|_{L^{ 2 } }   \le C\langle  t-s \rangle ^{-\frac 32} \|  \langle x \rangle ^{ \sigma }u \|_{L^{ 2 } } \text{ where}   \\&    A:= \langle x \rangle ^{-\sigma }e ^{ - \int _s^tdt' v(t')\cdot \nabla } e^{-\im \mathcal{K}_0  (t-s)}R_{ \mathcal{K}_0 }^{+}( \lambda )\chi   (\mathcal{K}_0)
 \langle \cdot  \rangle ^{-\sigma }  .\end{aligned}
\end{equation}
For $\Gamma (t,s):=\int _s^t v(t') dt'$,  the operator $A$ has an integral kernel, given by
\begin{equation*} \label{eq:lemg95}\begin{aligned}
  A(x,y)=\frac{
  \langle x \rangle ^{-\sigma }  \langle y \rangle ^{-\sigma }}{(2\pi )^3} \int _{t-s}^\infty dt' \int _{\R ^3}  d\xi e^{-\im (\sigma _3(\xi ^2 +\omega )  -  \lambda  )t' +\im \xi \cdot (x- \Gamma (t,s)-y) }
  \chi (\xi ^2 +\omega _0) .\end{aligned}
\end{equation*}
  If $| x |+ |y| \ge   |   \Gamma (t,s)| $
  we integrate by parts using   $e^{  \im  \sigma _3 \xi ^2t'}=\sigma _3 (  2\im  |\xi |t' )^{-1}
\frac d {d|\xi |} e^{  \im  \sigma _3  \xi ^2t'}$
\begin{equation} \label{eq:lemg96}\begin{aligned}&
  |A(x,y)|\le  C '
  \langle x \rangle ^{-\sigma }  \langle y \rangle ^{-\sigma }
    (1+| x- \Gamma (t,s)|+ |y|)^3
  \int _{t-s} ^\infty dt'  \langle t' \rangle ^{-3}
  \\&
  \le C  \langle x \rangle ^{-\sigma +3 }  \langle y \rangle ^{-\sigma +3}  \langle t-s\rangle  ^{-2}  .\end{aligned}
\end{equation}
By  $\chi (\xi ^2 +\omega _0)\neq 0 \Rightarrow | \xi |  > C _J  $  for a fixed $C _J >0 $,
if  $| x |+ |y| <  |   \Gamma (t,s)| $
\begin{equation*} \label{eq:lemg97}\begin{aligned}&
  |\nabla _\xi (\mp  \xi ^2 t' +\xi \cdot (x- \Gamma (t,s)-y))| \ge 2|\xi | t' -| x- \Gamma (t,s)|- |y|\\&  \ge  2C _Jt' -2 |   \Gamma (t,s)|  \ge 2t'  (  C _J  - C(C_0)\epsilon  ) > t'    C _J
     \end{aligned}
\end{equation*}
for $2(C_0)\epsilon _0<C _J  $   and  $t'\ge 1$.  This implies that integrating by parts we get
\begin{equation*} \label{eq:lemg98}\begin{aligned}&
  |A(x,y)|
  \le C  \langle x \rangle ^{-\sigma   }  \langle y \rangle ^{-\sigma  }   \langle t-s\rangle  ^{-2} . \end{aligned}
\end{equation*}
So  $ |A(x,y)|\le $   rhs of
\eqref{eq:lemg96}     $\forall$ $(x,y)$.
For $\sigma >9/2$   we get   \eqref{eq:lemg91}.
\qed

\begin{lemma}
  \label{lem:Tj} For $\epsilon _0$ small enough and for all $b=1,2,3$ we have
  \begin{equation} \label{eq:Tj1}  \begin{aligned}
 &  \| \int _0^t\mathcal{U}(t,s)
\textbf{T}_b(s) ds\| _{   L^2 ((0,T), L^{2,-\sigma }_x)}\le   \epsilon . \end{aligned}
 \end{equation}
  \end{lemma}
  \proof   By \eqref{Strichartzradiation}  and \eqref{eq:eqg-1} we have
     \begin{align}
 &  \| \int _0^t\mathcal{U}(t,s)
\textbf{T}_2(s) ds\| _{   L^2 ((0,T), L^{2,-\sigma } )}\le C_1  \| \int _0^t\mathcal{U}(t,s)
\textbf{T}_2(s) ds\| _{   L^2 ((0,T), L^{6 } )}\nonumber \\& \le   C_2\|  \textbf{T}_2   \| _{   L^2 ((0,T), L^{2,-\sigma }_x)  +L^1 ((0,T), L^{2  }_x) }\le C '(C_0) \epsilon  ^2< \epsilon .  \label{eq:Tj2}\end{align}
To get the bound for $b=1$ we focus for definiteness only on the term with the gradient, since the other can be bounded similarly.  By \eqref{Strichartzradiation}--\eqref{L^inftydiscrete}  we have
 \begin{align}
 & \nonumber  \| \int _0^t z^{\mu}(s)\overline{z}^{\nu}(s)   \mathcal{U}(t,s)
v\cdot \nabla R ^{+}_{\mu \nu }\textbf{G} _{\mu \nu}(s)  ds\| _{   L^2 ((0,T), L^{2,-\sigma } )}\\&  \label{eq:Tj3}\le C '  \| \int _0^t |z^{\mu}(s)\overline{z}^{\nu}(s) |   \langle t-s\rangle ^{-\frac{3}{2}} |v(s) | \, \,
\| \textbf{G} _{\mu \nu}(s)\| _{H^{1,\sigma}}  ds\| _{   L^2  (0,T) } \\& \le C  \|  z^{\mu} \overline{z}^{\nu} v   \| _{   L^2  (0,T)+L^1  (0,T)  }   \|   \textbf{G} _{\mu \nu}(t)\| _{L^\infty ((0,T) ,H^{1,\sigma})}\le C(C_0)\epsilon ^2,
\nonumber \end{align}
by Lemma  \ref{lem:lemg9} and    $ \|   \textbf{G} _{\mu \nu}(t)\| _{L^\infty ((0,T) ,H^{1,\sigma})}\le C $ for a fixed $C$, consequence of  \eqref{eq:nablaZ1} and \eqref{eq:eqh222}, Claim 5 in   Lemma \ref{lem:ExpH11} and   Lemma  \ref{lem:dotPi}.

Case $b=3$ uses   Lemma  \ref{lem:lemg9} and  can be proved
as Lemma 10.7 \cite{boussaidcuccagna}.  We have
 \begin{align}  &\partial _{z_j}Y =\sum _{| \mathbf{e}  \cdot(\mu-\nu)|>\omega _0} \mu _j\frac{z^\mu
\overline{z}^\nu}{z_j}
   R ^{+}_{\mu \nu}
     \textbf{G}_{\mu \nu} \, \label{eq:g variable1}\\&  \im \dot z_j-\mathbf{e} _j  z_j  = \partial _{\overline{z}_j}(  \widetilde{H} _2+  Z_0+  Z_1+   \resto  +     \langle  V(\cdot + \textbf{v} t) u
,  u  \rangle ) + \im \dot \Pi _b \partial _{\pi _b} z_j  \text{ with} \nonumber\end{align}
\begin{equation*} \widetilde{H} _2= \sum _{\substack{ | \alpha  +\beta  |=2\\
\mathbf{e}  \cdot (\alpha -\beta )=0}}
 g_{\alpha \beta} (   \Pi , \Pi  (f) )  z^\alpha
\overline{z}^\beta     \text{   with $g_{\alpha \beta} (   \Pi , \Pi  (f) ) =O(| \Pi  (f)|)$,}
\end{equation*}
by (1) and (3) in Lemma \ref{lem:ExpH11} (recall now $\Pi =\pi$) where we use also $\pi =\pi (t) =\Pi (t)$ and $\mathbf{e}_j (\pi (t))\equiv \mathbf{e}_j $, see under \eqref{eq:conjNLS1}.
By (H9), \eqref{eq:ExpHcoeff1} and  \eqref{Strichartzradiation}
\begin{equation} \label{eq:Tj4}  \widetilde{H} _2= \sum _ {(i,j ) \text{ s.t. } \textbf{e}_i =\textbf{e}_j }
 g_{ij  }   z _i
\overline{z} _j     \text{   with $\| g_{ij}\| _{L^\infty (0,T)} \le C (C_0) \epsilon ^2  $.}
\end{equation}
To bound \eqref{eq:Tj1}   for $b=3$ we need to bound
\begin{equation}  \label{eq:Tj5} \begin{aligned} & \mu _j
\|  \int _0^t  e ^{   \int _t^sdt' v(t')\cdot \nabla }e^{-\im \mathcal{K}_0   (t-s)}\frac{z^\mu
\overline{z}^\nu}{z_j} \varpi
  R ^{+}_{\mu \nu}
     \textbf{G}_{\mu \nu} ds \|_{L^2 (0,T) L^{ 2, - \sigma}_x } \ll \epsilon \end{aligned}      \end{equation}
for $\varpi =\partial _{\overline{z}_j} I,\dot \Pi _b \partial _{\pi _b} z_j  $   terms  in the rhs of the second equation  of \eqref{eq:g variable1}.
Case $I=Z_0$ is detailed in  Lemma 10.7 \cite{boussaidcuccagna}. Cases $I=Z_1, \resto$ admit easier proofs,
  skipped here as in  \cite{boussaidcuccagna}.  Case  $I= \widetilde{H} _2$, untreated in  \cite{boussaidcuccagna}, can be treated analogously, but deserves some attention.  It is enough to consider
$\partial _{\overline{z}_j} I=a_{ij}     z_i  $.     By    Lemma  \ref{lem:lemg9}, inequality \eqref{eq:Tj5} is a consequence of
\begin{equation}\label{eq:Tj6}    \| a_{ij}\| _{L^\infty (0,T)}
\|  {{z} ^{ {\mu}  }  \overline{z} ^ { \nu }} {z _j} ^{-1}
        z_i  \| _{L^2(0,T)} \le  C (C_0)\epsilon  ^2   C _0  \epsilon  .\end{equation}
This last inequality follows from  $\| a_{ij}\| _{L^\infty (0,T)} \le C (C_0) \epsilon ^2  $, by  \eqref{Strichartzradiation},  and by
$$ \mu _j{{z} ^{ {\mu}  }  \overline{z} ^ { \nu }}{z _j} ^{-1}
        z_i  =\mu _j{z} ^{ {\mu}'  }  \overline{z} ^ { \nu ' }  \text{with }   \|  {z} ^{ {\mu}'  }  \overline{z} ^ { \nu ' }   \|   _{L^2(0,T)}
				\le  C _0  \epsilon ,$$
which   follows   by
					\eqref{L^2discrete} and $|\mathbf{e}  \cdot (\mu '  -\nu ' )|=|\mathbf{e}  \cdot (\mu   -\nu )|>\omega _0$, the latter consequence of
$\textbf{e}_i =\textbf{e}_j$ in   \eqref{eq:Tj4}.
					Case   $I=\langle   V(\cdot +t\textbf{v})  u
,  u  \rangle $  follows by  \begin{equation*}\label{eq:Tj7}    \| {{z} ^{ {\mu}  }  \overline{z} ^ { \nu }}{z _j} ^{-1}\| _{L^\infty (0,T)}
\|  \partial _{ \overline{z} _j}I   \| _{L^1(0,T)}^{\frac 12}  \|  \partial _{ \overline{z} _j}I  \| _{L^\infty (0,T)}^{\frac 12}  \le  C (C_0)   \epsilon  ^{2+\frac 12} ,   \end{equation*}
					  consequence of  \eqref{L^inftydiscrete},
					Lemma \ref{lem:disppot1}  and  $\|  \partial _{ \overline{z} _j}I   \| _{L^\infty (0,T)} \le C$ (which is easy).

For $ \varpi = \dot \Pi _b \partial _{\pi _b} z_j $, since  $\|  \partial _{\pi _b} z_j\| _{L^\infty (0,T)}\le C (C_0)   \epsilon$ by
Lemma \ref{lem:reg5},
\begin{equation*}     \| {{z} ^{ {\mu}  }  \overline{z} ^ { \nu }}{z _j} ^{-1} \| _{L^\infty (0,T)} \|  \partial _{\pi _b} z_j\| _{L^\infty (0,T)}
\|  \dot \Pi _b  \| _{L^1(0,T)}^{\frac 12}  \|  \dot \Pi _b  \| _{L^\infty (0,T)}^{\frac 12}  \le  C (C_0)   \epsilon  ^{2+\frac 12} ,   \end{equation*}
where      we use also Lemma \ref{lem:dotPi} and  $\|   \dot \Pi _b   \| _{L^\infty (0,T)} \le C$ (which is easy).

					\qed

To
complete the proof of Lemma \ref{lemma:bound g}   it  remains to  bound the term
$P_c\mathcal{V}  _2^{D}
	h  $ in the rhs of \eqref{eq:eqg1}.  We proceed as in Sect. \ref{sec:dispersion},
	where to bootstrap estimates we exploited  that  the interval $[0,T]$ was large.
		Notice that so far we have not exploited the fact that   in Lemma \ref{lemma:bound g}
we have $T\ge \epsilon _0^{-1}$.

\begin{lemma}
  \label{lem:VD} For any    $L>1$,
		if  $\epsilon _0>0$  is small enough
	we have
  \begin{equation} \label{eq:VD1}  \begin{aligned}
 & \| \int _0^t\mathcal{U}(t,s)P_c   \mathcal{V}  _2^{D}
	h  ds\| _{L^2 ((0,T), L^{2,-\sigma }_x)}\le   L^{-1}{C_0}   \epsilon . \end{aligned}
 \end{equation}
  \end{lemma}
  \proof We have for some fixed large $A\ll T$
  \begin{equation*} \label{eq:VD2}  \begin{aligned}
 & \| \int _0^t\mathcal{U}(t,s)P_c   \mathcal{V}  _2^{D}
	h  ds\| _{L^2 ((0,A), L^{2,-\sigma }_x)}\le   C_A  \epsilon . \end{aligned}
 \end{equation*}
  Next,  for $k(A)\sim A ^{-1}$
  defined  in \eqref{eq: Ijminor41}, we have
    \begin{align} \nonumber
 & \| \int _0^{t-A}\mathcal{U}(t,s)P_c   \mathcal{V}  _2^{D}
	h  ds\| _{L^2 ((A,T), L^{2,-\sigma }_x)} \lesssim  \left \| \int _0^{t-A}  \langle t-s \rangle ^{-\frac{3}{2}} \| h(s) \| _{  L^{6 }_x  }  ds\right \| _{L^2 (A,T) } \\&  \le  \sqrt{k(A)}  \| h \| _{ L^2 ((0,T) ,L^{6 }_x)} \le  \sqrt{k(A)} C_0\epsilon   .  \nonumber \end{align}

\noindent We write   $P_c   \mathcal{V}  _2^{D} h =   \mathcal{V}  _2^{D}
h- P_d   \mathcal{V}  _2^{D}
	h $. Then,  by \eqref{eq:divcenter10},   we have
   \begin{align}
 & \nonumber \| \int _ {t-A}^t\mathcal{U}(t,s) P_d    \mathcal{V}  _2^{D}
	h  ds\| _{L^2 ((A,T), L^{2,-\sigma }_x)}  \\&\nonumber  \le  \sum _j  \left \| \int _ {t-A}^t       \langle t-s \rangle ^{-\frac 32}     \|        \mathcal{V}  _2(\cdot +{\textbf{D} } ) e_j  \| _{  L^{2  }_x }     ds\right \| _{L^2   (A,T) }    \| 	 h \|   _{  L^\infty (  (0,T),L^{2  }_x) } \\& \le
	     	C C_0      \epsilon \gamma  \int _A^T \langle \textbf{D} (s) \rangle ^{-2} \|         \langle t-s \rangle ^{-\frac 32}   \| _{L^2 _{t}  (A,T) }  ds   \le   C^{\prime  } C_0 \epsilon   ^2   . \nonumber \end{align}
 We have by a version of \eqref{eq:duh462},
  for  $\mathcal{V}  _2^{\sigma}(x ) := \mathcal{V}  _2(x )  \langle x\rangle ^{ \sigma } $
\begin{align}
 & \nonumber \| \int _ {t-A}^t\mathcal{U}(t,s) F (|\im \nabla |\le N)   \mathcal{V}  _2^{D}
	h  ds\| _{L^2 ((A,T), L^{2,-\sigma }_x)}  \\&  \nonumber \le
		\left \|    \int _ {t-A}^t      \|  \langle x \rangle ^{-\sigma }       \mathcal{U}(t,s)  F (|\im \nabla |\le N)    \mathcal{V}  _2^{\sigma}(\cdot  +{\textbf{D} } )    \| _{   {2  }  \to 2 }      \| 	 h \|   _{   L^{6  }_x  }  ds\right \| _{L^2   (A,T) }
	\\& \nonumber  \le    C'(N)\gamma   \left \|    \int _ {t-A}^t     \left  \langle    \textbf{D}(s) +\int _s^t v(\tau ) d\tau  \right \rangle ^{-2}     \| 	 h \|   _{   L^{6  }_x  }  ds\right \| _{L^2   (A,T) } \\&   \le   2 C'(N)\gamma  \sqrt{A} \|        \langle     \textbf{D}(t)     \rangle ^{-2}  \|   _{L^1   (A,T) }       \| 	h \|   _{ L^2   ((0,T) , L^{6  }_x ) }  <  C (N) \sqrt{A}   C_0 \epsilon   ^{2} ,  \nonumber \end{align}
	where we used Minkowski inequality and,  for $t-A \le s \le t$ and $\sqrt{2} c(T)A<1$,
	\begin{equation*}
\begin{aligned} &     1+| \textbf{D}(s)|^2 <  2( 1+| \textbf{D}(s)|^2  - (c(T)A)^2)) \le 2  ( 1+| \textbf{D}(s)+\int _s^t v(\tau ) d\tau|^2   ).
\end{aligned}
\end{equation*}

	\noindent
  Finally we have
	\begin{align}
 & \nonumber \| \int _ {t-A}^t\mathcal{U}(t,s) F (|\im \nabla |\ge N)   \mathcal{V}  _2^{D}
	h  ds\| _{L^2 ((A,T), L^{2,-\sigma }_x)}  \\&  \label{eq:VD5} \lesssim    N ^{-1}
		\left \|    \int _ {t-A}^t      \|  \langle x\rangle ^{ -\sigma }      \mathcal{U}(t,s)        \mathcal{V}  _2^{\sigma}(y +{\textbf{D} } )  \| _{   {2  }  \to 2 }      \| 	h \|   _{   L^{6  }_x  }  ds\right \| _{L^2   (A,T) }
		\\&  \nonumber+   N ^{- \frac 12 }
			\left \|    \int _ {t-A}^t      \|  [ {\nabla}{\langle \im \nabla \rangle ^{-\frac 12}}  ,\langle x \rangle ^{-\sigma }]       \mathcal{U}(t,s)      \mathcal{V}  _2^{\sigma}(y +{\textbf{D} } )     \| _{   {2  }  \to 2 }      \| 	 h \|   _{   L^{6  }_x  }  ds\right \| _{L^2   (A,T) }  \\& \nonumber +    N ^{- \frac 12 }
			\left \|    \int _ {t-A}^t      \|     \langle x \rangle ^{-\sigma } {\nabla}{\langle \im \nabla \rangle ^{-\frac 12}}         \mathcal{U}(t,s)       \mathcal{V}  _2^{\sigma}(y +{\textbf{D} } )   \| _{   {2  }  \to 2 }      \| 	 h \|   _{   L^{6  }_x  }  ds\right \| _{L^2   (A,T) }
		 \end{align}
		By Young's inequality, the 2nd  and 3rd line of  \eqref{eq:VD5}   are  bounded by
		\begin{equation*}   \label{eq:VD7}  \begin{aligned}
 &   N ^{- \frac 12 }  \gamma  \left \|    \int _ {t-A}^t         \| 	 h \|   _{   L^{6  }_x  }  ds\right \| _{L^2   (A,T) }  \le A \gamma C_0N ^{- \frac 12 } \epsilon ,
 \end{aligned}
\end{equation*}
by   $\|   [  {\langle \im \nabla \rangle } ^{-\frac{1}{2}} \nabla , \langle x \rangle ^{-\sigma }] \| _{2\to 2} \lesssim 1$ and $\|   \mathcal{U}(t,s)      \mathcal{V}  _2^{\sigma}(y +{\textbf{D} } ) \| _{2\to 2} \lesssim \gamma$.
By the local smoothing effect  \eqref{eq:locsm} and Young's inequality  the
4th line  of  \eqref{eq:VD5} is  $\lesssim$
\begin{equation*}   \label{eq:VD8}  \begin{aligned}
 &       \gamma   N^{-\frac{1}{2}} \| 	h \|   _{ L^2   ((0,T) , L^{6  }_x ) }
  ( \sup _s \int h(t,s) dt + \sup _t \int h(t,s) ds)  <     C C_0  \gamma  \sqrt{A} N^{-\frac{1}{2}}      \epsilon \,  , \\&
 h(t,s):=    \chi _{[ A,T]} ( t)\ \ \chi _{[t-A,t]} ( s)  \ \  \|     \langle x \rangle ^{-\sigma } {\nabla}{\langle \im \nabla \rangle ^{-\frac 12}}        e^{\im (s-t) \mathcal{K}_0 }           \| _{   {2  }  \to 2 }  .  \end{aligned}
\end{equation*}
  Taking $N$ sufficiently large we obtain the desired bound.
		 \qed

\section{Proof of Lemma \ref{lem:scattf} }
\label{sec:scattf}
	 By \eqref{eq:eqh222}
		we can write    \begin{align} & \nonumber
    h(t) =  \mathcal{U} (t,0) h(0) -\im  \sigma _{3}\gamma    \int _0^t   \mathcal{U} (t, s)  P_c(\mathcal{K}_{\omega _0
 }  ) V (\cdot + \textbf{v}s+ D' +\resto ^{0,2} )
	 {h} (s) \,    ds  - \\&  \nonumber   \im    \int _0^t   \mathcal{U} (t, s)  P_c(\mathcal{K}_{\omega _0
 }  )   (  R_1 ^{\prime\prime } (s)+ R_2 ^{\prime\prime }  (s)  + \sum _{|\textbf{e}    \cdot(\mu-\nu)|>\omega   _0  }  z^\mu  (s)\overline{z}^\nu  (s)\mathbf{G}_{\mu \nu}(  s,\Pi (f)   )      )    ds
. \nonumber \end{align}
Taking $t_1<t_2$ we have 	
		\begin{align} & \nonumber  \|   \mathcal{U} ^{-1}(t_2,0)  h(t_2)- \mathcal{U} ^{-1}(t_1,0)  h(t_1)\| _{H^1}\lesssim \gamma    \|  V (x )\| _{W^{1,\frac 32  }_x}     \| h   \| _{L^2 ([t_1,t_2],W^{1, 6} )} \\&  + \| R_1^{\prime\prime } \| _{L^1 ([t_1,t_2],H^{ 1 } )} +    \| R_2^{\prime\prime } \| _{L^2 ([t_1,t_2],H^{ 1 ,s} )} +\sum _{|\textbf{e}    \cdot \mu |>\omega   _0  }    \| z^\mu  \| _{L^2 ( t_1,t_2  )}  .\label{eq:resolh1}
\end{align}
		Since  \eqref{Strichartzradiation}--\eqref{L^inftydiscrete}  on $I=[0,\infty )$ imply that
		the rhs of  \eqref{eq:resolh1}  goes to 0 as $t_1\nearrow \infty$,
	   then  $\displaystyle \lim _{t\nearrow \infty}\mathcal U^{-1}(t,0) h(t)=h_+$
	exists in $H^1$.
	Then, for $\textbf{A} \cdot\nabla = \sum _{a=1}^3\textbf{A}_a  \partial _{x_a}$,
\begin{equation}\label{eq:scattf1}
 \begin{aligned} &  \lim _{t\nearrow \infty} \| h(t) -h_0(t)  \| _{H^1} =0, \,   h_0(t):=
e^{  \im  t \sigma _3(  \Delta - \omega _0 )  -\int _0^t ( \textbf{A}(s) \cdot\nabla +
  \im  \textbf{A}_4 (s) \sigma _3) ds } h_+,
\end{aligned}
\end{equation}
	  with the coefficients $\textbf{A}_j$
	in  \eqref{eq:eqh222}.
	We have
	\begin{equation*}\label{eq:scattf2}\begin{aligned} & e^{-  J \frac{v_0 \cdot x}{2}}   M   h_0(t) = e^{- J \frac{v_0 \cdot x}{2}} e^{     t J(-\Delta + \omega _0 )} e^{ \int _0^t (
  J \textbf{A}_4 (s) -\textbf{A}(s) \cdot\nabla) ds } e^{  J \frac{v_0 \cdot x}{2}}     f_+
\end{aligned}
\end{equation*}
	 for $f_+:=e^{-\frac 12 J v_0 \cdot x}   Mh_+$  and by  $J=-M \im \sigma _3 M^{-1}$, see \eqref{eq:Homega2}.  By
\eqref{eq:comm1}--\eqref{eq:comm2}
\begin{equation*}\label{eq:scattf4}\begin{aligned} &
 e^{    -  t J \Delta  }
e^{ \frac 12 J v_0 \cdot x} u_0 =e^{ \frac 12 J v_0 \cdot x}
 e^{ t  \left ( J\frac{v_0^2}{4}+ v_0\cdot \nabla    \right )}
 e^{  -    t J \Delta  }
 u_0.
\end{aligned}
\end{equation*}
Recall $e^{   v_0\cdot \nabla }u_0(x)=u_0(x+v_0)$.
Hence we conclude
  \begin{equation}\label{eq:scattf3}\begin{aligned} & e^{-\frac 12 J v_0 \cdot x}   M   h_0(t) =   e^{   -  t J \Delta  }
	e^{J\varsigma (t) \cdot \Diamond}  f_+ \\&  	 e^{J\varsigma (t) \cdot \Diamond} := {e^{     J  \left ( t \left (\omega _0 +   \frac{v_0^2}{4}    \right )      + \int _0^t    {A}'_4 (s)   ds\right ) \Diamond _4}
	  e^{   J \sum _{a=1}^{3} \left ( - v_{0a} +\int _0^t  {A}_a'(s)   \right )   \Diamond _a  }}
\end{aligned}
\end{equation}
for $A'_a=\textbf{A}_a $  for $a\le 3$  and   $A'_4=\textbf{A}_4-\frac 12 \textbf{A}\cdot v_0$   the coefficients  in \eqref{eq:eqf22},  see under  \eqref{eq:eqh22}.
Then \eqref{eq:scattf1}   and $f(t) =e^{-\frac 12 J v_0 \cdot x}   M   h (t)$     yield Lemma \ref{lem:scattf}.\qed

\section{End of the proof of Theorem \ref{theorem-1.1}}
\label{sec:completion}

Recall that,   for $p'$ and $\Pi $ with limits as $t\nearrow\infty $, we have
  \begin{equation}\label{eq:compl0}
\begin{aligned}  u(t)=    e^{J \tau '(t)\cdot \Diamond} (  \Phi _{p'(t) } +P(p'(t))P(\Pi (t)) r'(t)) \text{ , } r'=e^{J \resto ^{0,2}\cdot \Diamond} (f +\mathbf{S}^{0,1}).      \end{aligned}
\end{equation}
By  Theorem \ref{thm:mainbounds}, we know  that $ e^{\im \Theta (t,x-D(t))}  r(t,x-D(t) )=A(t,x)+\widetilde{r}(t,x)$ with $A(t,x)$ and $\widetilde{r}(t,x)$ satisfying the conclusions of  Theorem \ref{theorem-1.1}.

To finish   Theorem \ref{theorem-1.1} we have yet to prove   \eqref{eq:scattering}.
 We know that $f$ satisfies Lemma \ref{lem:scattf}.
To prove \eqref{eq:scattering},
 it is sufficient to show that  there is a  $\zeta _0\in \R^4$ s.t.

 \begin{equation}\label{eq:compl2}
\begin{aligned} & \lim _{t \nearrow \infty}   (\tau '(t)+\varsigma (t)) =\zeta _0\end{aligned}
\end{equation}
for $\tau  '$ the phase   in \eqref{eq:compl0} and $\varsigma$
the function in Lemma \ref{lem:scattf}  written explicitly in \eqref{eq:scattf3}.
 It is easy to see  that  when we plug \eqref{eq:compl0}  in   $\dot u =J\nabla \mathbf{E} (u)$
  we get
\begin{equation}\label{eq:compl3}
\begin{aligned} &  \dot f=J(-\Delta -\dot \tau '\cdot \Diamond - \dot  \resto ^{0,2}\cdot \Diamond ) f+G_1(u),
\end{aligned}
\end{equation}
with $G_1(u)\in  C^0 (  H^1_x,   L^1_x)$.
 On the other hand,
$f$ satisfies also \eqref{eq:eqf22}, which is of the form
\begin{equation} \label{eq:compl4} \begin{aligned} &              \dot f = J(-\Delta -\lambda (p_0) \cdot \Diamond +  A' \cdot   \Diamond )
f+G_2(u)      \end{aligned}  \end{equation}
with  $G_2(u)\in  C^0 (  H^1_x,   L^1_x)$. Then, by the definition of $\varsigma$  we have \begin{equation} \label{eq:compl5} \begin{aligned} &
    ( \dot \tau '+ \dot  \resto ^{0,2} +  \dot \varsigma   )\cdot \Diamond f= G_1(u)
-G_2(u)      \end{aligned}  \end{equation}
In turn, all the coefficients $\dot \tau'$,  $\dot  \resto ^{0,2}$, $\dot \varsigma$   in  \eqref{eq:compl5} are in  $   C^0 (  H^1_x,   \R ^4)$.
By an argument in Lemma 13.8 \cite{Cu3} we  conclude
$\dot \tau  '+ \dot  \resto ^{0,2} +  \dot \varsigma =0$. Then, by
$\lim _{t \nearrow\infty}\resto ^{0,2} =0$, see the proof of \eqref{eq:partialmod21}, we obtain
\eqref{eq:compl2}.

\section*{Acknowledgments}   S.C. was partially funded  by a  grant  FRA 2009
from the University of Trieste and by the grant FIRB 2012 (Dinamiche Dispersive).
Part of this work was done when M.M. was visiting the Department of Mathematics and Geosciences, University
of Trieste and he wishes to thank its hospitality.
M.M. would like to thank Professor Yoshio Tsutsumi for covering the travel expenses to Trieste.
M.M. is supported by the Japan Society for the Promotion of Science (JSPS) with the Grant-in-Aid for Young Scientists (B) 24740081.

Department of Mathematics and Geosciences,  University
of Trieste, via Valerio  12/1  Trieste, 34127  Italy

{\it E-mail Address}: {\tt scuccagna@units.it}
\\

Department of Mathematics and Informatics,
Faculty of Science,
Chiba University,
Chiba 263-8522, Japan

{\it E-mail Address}: {\tt maeda@math.s.chiba-u.ac.jp}


\begin{thebibliography}{CP03}




\bibitem{ASFS}
  W.K.Abou Salem, J.Fr\"ohlich,  I.M.  Sigal,  {\em  Colliding solitons for the nonlinear Schr\"odinger
	equation\/},  Comm. Math. Phys. 291 (2009), 151--176.


\bibitem{AS}
  W.K.Abou Salem,      C. Sulem,  {\em  Resonant tunneling of fast solitons through large potential barriers} Canad. J. Math. 63 (2011),  1201--1219.

\bibitem{ALS}
  W.K.Abou Salem,    X. Liu, C. Sulem,  {\em
Numerical simulation of resonant tunneling of fast solitons for the nonlinear Schr\"odinger
	equation\/},  Discrete Contin. Dyn. Syst. 29 (2011),  1637--1649.



\bibitem {bambusi}
  D.Bambusi,  {\em   Asymptotic stability of ground states in some Hamiltonian PDEs with symmetry\/},   Comm. Math. Physics, 320 (2013),  499--542.

\bibitem{bambusicuccagna}
  D.Bambusi, S.Cuccagna, {\em On dispersion of
small energy solutions of the nonlinear Klein Gordon equation with a
potential},   Amer. Math. Jour., 133 (2011), 1421--1468.

\bibitem{beceanu}
  M. Beceanu  {\em New estimates for a time-dependent Schr\"odinger equation}, Duke Math. Jour.   159  (2011), 417--477.



\bibitem{boussaidcuccagna}
  N.Boussaid, S.Cuccagna,  {\em  On stability of standing waves of nonlinear Dirac equations\/},    Comm. in Partial Diff. Eq.  37  (2012),   1001--1056.


\bibitem{BP2}
V.Buslaev, G.Perelman, {\em On the stability of solitary waves for
nonlinear Schr\"odinger equations\/},   Nonlinear evolution
equations, editor N.N. Uraltseva, Transl. Ser. 2, 164, Amer. Math.
Soc.,
    75--98, Amer. Math. Soc., Providence (1995).

\bibitem{cai} K.Cai,
 {\em Fine properties of charge transfer models}, preprint   (2003),
 arXiv:math-ph/0311048v1.



\bibitem{CLC}
  C. Cote, S. Le Coz,  {\em High-speed excited multi-solitons in nonlinear
Schr\"odinger equations\/}, J. Math. Pures Appl.   96 (2011),   135--166.



\bibitem{Cu0}
  S.Cuccagna,  {\em On the Darboux and Birkhoff   steps in the asymptotic
  stability of    solitons\/},  Rend. Istit. Mat. Univ. Trieste   {\bf  44} (2012), 197--257.
	
		
\bibitem{Cu2}
  S.Cuccagna,  {\em The Hamiltonian structure of the nonlinear
Schr\"odinger equation and   the  asymptotic stability of its
ground states\/}, Comm. Math. Physics,
  305 (2011),  279-331.	
			
			
		
\bibitem{Cu3}
  S.Cuccagna,  {\em On asymptotic stability of moving  ground
 states    of the nonlinear Schr\"odinger equation\/}, to appear Trans. Amer. Math. Soc.	
			
			
			
			

 \bibitem {Cu6}
   S.Cuccagna, {\em Stabilization of solutions to  nonlinear
  Schr\"odinger equations},  Comm. Pure App. Math.  {54} (2001),
 1110--1145, erratum  Comm. Pure Appl. Math.  58  (2005),   147.


\bibitem{Cu7}   S.Cuccagna, {\em On asymptotic stability
of ground states of NLS},  Rev. Math. Phys.  {15} (2003),
877--903.



\bibitem{CM1} S.Cuccagna, M.Maeda, {\em On small energy stabilization in the NLS with a trapping potential},  arXiv:1309.0655v1 arXiv:1309.0655v1 .


\bibitem{CPV} S.Cuccagna, D.Pelinovsky, V.Vougalter, {\em Spectra of
positive and negative energies in the linearization of the NLS problem},  Comm.
Pure Appl. Math.  {58} (2005),   1--29. 	



		
			\bibitem{JFGS} J.Fr\"ohlich,B.Jonsson,S.Gustafson,I.M-Sigal
{\em Long time motion of NLS solitary waves in a confining potential}, Ann. H. Poinc.  7 (2006),  621--660.
		
			\bibitem{JFGS0} J.Fr\"ohlich,
 B.L.Jonsson,   S.Gustafson, I.M.Sigal
{\em  Solitary wave dynamics in an external potential} Commun. Math. Phys. 250(2004), 613--642.



	\bibitem {GSS1}   M.Grillakis,
J.Shatah, W.Strauss, {\em Stability of solitary waves in the presence of
symmetries, I },  Jour. Funct. An.   {74} (1987),   160--197.

\bibitem  {GSS2} M.Grillakis,
J.Shatah, W.Strauss, {\em Stability of solitary
    waves in the presence of symmetries, II}, Jour. Funct. An.   {94} (1990),
    308--348. 	
		
		
\bibitem{GHM} R.Goodman,   P.Holmes,  M.I.Weinstein  {\em Strong NLS soliton-defect interactions}, Phys. D 192 (2004),   215--248


\bibitem{GNT}
S.Gustafson,   K.Nakanishi, T.P.Tsai,  {\em Asymptotic stability and completeness in the energy space for nonlinear Schr\"odinger equations with small solitary waves}, Int. Math. Res. Not., 2004 (2004) no. 66, 3559--3584



		
			\bibitem {HMZ1}
		J.Holmer, J.Marzuola, M. Zworski,  {\em  Soliton splitting by external delta potentials} J. Nonlinear Sci. 17 (2007),   349--367.
		
		
		
		
			\bibitem {HMZ2}
		J.Holmer, J.Marzuola, M. Zworski,  {\em  Fast soliton scattering by delta impurities} Comm. Math. Phys. 274 (2007),   187--216.
		
		
		\bibitem {HZ1}
		J.Holmer,   M. Zworski,  {\em Slow soliton interaction with delta impurities}
J. Mod. Dyn. 1 (2007),   689--718.


\bibitem {HZ2} J.Holmer,   M. Zworski,  {\em
		Soliton interaction with slowly varying potentials},   Internat. Math. Res. Notices (2008), Art. ID runn026, 36 pp.
		
		
		
		
		\bibitem{JSS}
 J.L.Journe, A.Soffer, C.D.Sogge,
{\em Decay estimates for Schrodinger operators }, Comm.Pure Appl. Mat.
  44 (1991), pp. 573--604.





\bibitem  {lp} F.Linares, G. Ponce, {\em Introduction to nonlinear dispersive equations} Universitext. Springer, New York, 2009.


\bibitem  {MM}
Y.Martel, F.Merle, {\em Multi solitary waves for nonlinear Schr\"odinger equations},
 Ann. Inst. H. Poinc.  Anal. Non Lin., 23   (2006),   849--864.

\bibitem  {MMT}
Y.Martel, F.Merle, T.P.Tsai, {\em Stability in $H^1$ of the sum of K solitary waves for some nonlinear Schr\"odinger equations},
Duke Math. J., 133   (2006),   405--466

\bibitem  {M1}  C.Munoz, {\em Sharp inelastic character of slowly varying NLS solitons},   arXiv:1202.5807.


\bibitem  {M2}  C.Munoz, {\em On the soliton dynamics under slowly varying medium for Nonlinear Schr\"odinger equations},     Math. Ann., 353 (2012), 867--943.

 \bibitem {NPT}K.Nakanishi, T.V.Phan, T.P.Tsai, {\em Small solutions of nonlinear Schr\"odinger equations near first excited states}, Jour. Funct. Analysis,  263 (2012), 703--781.

\bibitem {NS}K.Nakanishi, W.Schlag, {\em Global dynamics above the ground state energy for the cubic NLS equation in 3D},  Cal.Var.Par.Diff.Eq.   44 (2012), 1--45.



\bibitem  {perelman1} G.Perelman, {\em  Two soliton collision for nonlinear Schr\"odinger equations in dimension 1}, Ann. Inst. H. Poinc.  Anal. Non Lin. 28 (2011), 357--384.


\bibitem  {perelman2} G.Perelman, {\em  A remark on soliton-potential interactions for nonlinear     Schr\"odinger equations},
Int. Math. Res. Not.    37 (2005), 2289--2313



\bibitem  {perelman3} G.Perelman, {\em  Asymptotic stability of multi-soliton solutions for nonlinear     Schr\"odinger equations},  Comm. Partial Diff.  29 (2004),  1051--1095.



\bibitem  {perelman4} G.Perelman, {\em
Some results on the scattering of weakly interacting solitons for nonlinear Schr\"odinger equations}, Spectral theory, microlocal analysis, singular manifolds, 78--137,
Math. Top., 14, Akademie Verlag, Berlin, 1997.





\bibitem  {RSS1} I.Rodnianski, W.Schlag, A.Soffer, {\em
Dispersive analysis of charge transfer models},
Comm. Pure   Appl. Math.   58(2005), 149--216.

\bibitem  {RSS2} I.Rodnianski, W.Schlag, A.Soffer, {\em
Asymptotic stability of N-soliton states of NLS}, preprint (2003),
	arXiv:math/0309114v1 .

 \bibitem  {schlag}W.Schlag,  {\em
 Stable manifolds for an orbitally unstable nonlinear Schr�dinger equation}, Annals of Math.  169 (2009),  139--227.



\bibitem {SW3}
A.Soffer, M.I.Weinstein, {\em Resonances, radiation damping and
instability in Hamiltonian nonlinear wave equations},  Inv.
Math., 136
 (1999),
  9--74.


 \bibitem {SW4}    A.Soffer, M.I.Weinstein,
 {\em Selection of the ground state for nonlinear Schr\"odinger
 equations}, Rev. Math. Phys.  16 (2004), pp.  977--1071.


\bibitem{strauss}
W.Strauss,  {\em Nonlinear wave equations},   CBMS Regional Conf.
Ser. Mat.   AMS  76    (1989).


 \bibitem {TY3}
 {  T.P.Tsai, H.T.Yau}, {\em Classification of asymptotic profiles
 for nonlinear Schr\"odinger equations with small initial data}, Adv.
 Theor. Math. Phys.  {6} (2002), pp.  107--139.

\bibitem{W2}
   M.I.Weinstein, {\em Lyapunov stability of ground
states of nonlinear dispersive equations},  Comm. Pure Appl. Math.
 39  (1986),    51--68.


\bibitem{Y1}   K.Yajima, The $W^{k,p}$-continuity of wave operators
for Schr\"{o}dinger operators, J. Math. Soc. Japan,  {47} (1995),
pp. 551--581.



	\end{thebibliography}
\end{document}